\documentclass[11pt]{sat}
\usepackage{amsmath,amsthm,amscd,amssymb}
\usepackage{latexsym}

\newcommand{\R}{{\mathbb{R}}}
\newcommand{\Z}{{\mathbb{Z}}}
\newcommand{\C}{{\mathbb{C}}}
\newcommand{\D}{{\mathbb{D}}}
\newcommand{\h}{{\mathcal{H}}}
\newcommand{\pol}{{\mathcal{P}}}
\newcommand{\calH}{{\mathcal H}}
\newcommand{\calL}{{\mathcal L}}
\newcommand{\calM}{{\mathcal M}}
\newcommand{\calV}{{\mathcal V}}

\newcommand{\bbC}{{\mathbb{C}}}
\newcommand{\bbD}{{\mathbb{D}}}
\newcommand{\bbL}{{\mathbb{L}}}
\newcommand{\bbR}{{\mathbb{R}}}
\newcommand{\bbS}{{\mathbb{S}}}
\newcommand{\bbT}{{\mathbb{T}}}
\newcommand{\bbU}{{\mathbb{U}}}
\newcommand{\bbZ}{{\mathbb{Z}}}

\allowdisplaybreaks
\numberwithin{equation}{section}

\newtheorem{theorem}{Theorem}[section]
\newtheorem{proposition}[theorem]{Proposition}
\newtheorem{lemma}[theorem]{Lemma}
\newtheorem{corollary}[theorem]{Corollary}
\theoremstyle{definition}

\newtheorem{conjecture}[theorem]{Conjecture}

\theoremstyle{remark}
\newtheorem*{remark}{Remark}
\newtheorem*{remarks}{Remarks}
\newtheorem*{definition}{Definition}

\newcommand{\abs}[1]{\lvert#1\rvert}

\newcommand{\norm}[1]{\lVert#1\rVert}
\newcommand{\snorm}[1]{{ \pmb\lvert}#1{ \pmb\rvert}}

\DeclareMathOperator{\Tr}{Tr}
\DeclareMathOperator{\diam}{diam}
\DeclareMathOperator{\Ima}{Im}
\DeclareMathOperator{\Span}{span}
\DeclareMathOperator{\Ran}{Ran}
\DeclareMathOperator{\lc}{linear\  combination\  of}

\newcommand{\ang}[1]{\langle\! \langle#1\rangle\! \rangle}

\def\jap#1{\langle #1 \rangle}

\newcommand{\lb}{\label}
\newcommand{\bdone}{{\boldsymbol{1}}}
\newcommand{\bdzero}{{\boldsymbol{0}}}
\newcommand{\bddot}{{\boldsymbol{\cdot}}}

\newcommand{\supp}{\text{\rm{supp}}}
\newcommand{\ess}{\text{\rm{ess}}}
\newcommand{\s}{\text{\rm{s}}}
\newcommand{\ol}{\overline}
\newcommand{\f}{\frac}
\newcommand{\dott}{\,\cdot\,}
\newcommand{\cvh}{\text{\rm{cvh}}}
\newcommand{\tr}{\text{\rm{tr}}}
\newcommand{\dist}{\text{\rm{dist}}}
\newcommand{\ti}{\tilde  }
\newcommand{\wti}{\widetilde  }

\let\det=\undefined\DeclareMathOperator{\det}{det}

\newcounter{smalllist}
\newenvironment{SL}{\begin{list}{{\rm\roman{smalllist})}}{%
\setlength{\topsep}{0mm}\setlength{\parsep}{0mm}\setlength{\itemsep}
{0mm}%
\setlength{\labelwidth}{2em}\setlength{\leftmargin}{2em}\usecounter
{smalllist}%
}}{\end{list}}

\title{The Analytic Theory of Matrix
Orthogonal Polynomials}
\def\shorttitle{Matrix Orthogonal Polynomials}

\author{David Damanik, Alexander Pushnitski, and Barry Simon}
\def\shortauthor{D.\ Damanik, A.\ Pushnitski, B.\ Simon }

\def\versiondate{January 30, 2008}

\def\abstracttext{We survey the analytic theory of matrix orthogonal polynomials.}

\def\MSCnumbers{42C05, 47B36, 30C10}

\def\keywords{orthogonal polynomials, matrix-valued measures, block
Jacobi matrices, block CMV matrices}

\def\startpagenumber{1}
\def\volumenumber{4}
\def\year{2008}

\newcommand{\beginddoc}{
\begin{document}
\maketitle
\begin{abstract}
\abstracttext
\vskip1pt MSC: \MSCnumbers
\ifx\keywords\empty\else\vskip1pt keywords: \keywords\fi
\end{abstract}
\insert\footins{\scriptsize
\medskip
\baselineskip 8pt
\leftline{Surveys in Approximation Theory}
\leftline{Volume \volumenumber, \year.
pp.~\thepage--\pageref{endpage}.}
\leftline{\copyright\ \year\ Surveys in Approximation Theory.}
\leftline{ISSN 1555-578X}
\leftline{All rights of reproduction in any form reserved.}
\smallskip
\par\allowbreak}
\tableofcontents}
\renewcommand\rightmark{\ifodd\thepage{\it \hfill\shorttitle\hfill}\else {\it \hfill\shortauthor\hfill}\fi}
\markboth{{\it \shortauthor}}{{\it \shorttitle}}
\markright{{\it \shorttitle}}
\def\endddoc{\label{endpage}\end{document}}
\date{{\small \versiondate}}
\setlength\oddsidemargin{0pc}
\setlength\evensidemargin{0pc}
\setlength\topmargin{0in}
\setlength\textwidth{6.5in}
\setlength\textheight{8.6in}
\beginddoc





\section{Introduction}

\subsection{Introduction and Overview} \lb{s1.1}

Orthogonal polynomials on the real line (OPRL) were developed in
the nineteenth century and orthogonal polynomials on the unit
circle (OPUC) were initially developed around 1920 by Szeg\H{o}.
Their matrix analogues are of much more recent vintage. They were
originally developed in the MOPUC case indirectly in the study of
prediction theory \cite{HL1,HL2, Kol,Krein1,Krein2,Lev,Wie} in the
period 1940--1960. The connection to OPUC in the scalar case was
discovered by Krein \cite{Krein1}. Much of the theory since is in
the electrical engineering literature \cite
{DG92,DGK,DGK3,DGK79,DGK2,DGK81,Kai74,Kai87,K6.1,KVM,YK}; see also
\cite{Ger81,Ge92,GC,GW07,MarRod}.

The corresponding real line theory (MOPRL) is still more recent:
Following early work of Krein \cite{Krein3} and Berezan'ski
\cite{Ber} on block Jacobi matrices, mainly as applied to
self-adjoint extensions, there was a seminal paper of
Aptekarev--Nikishin \cite{AN} and a flurry of papers since the
1990s \cite{Berg,Berg-Duran1,Berg-Duran2,CFMV,CMV05,CMV07,CG05,CG,CG07,CH,CGR,DGIM,DS,Dur1,Dur95,
Dur96,Dur6,Dur99,Dur7,Dur8,DuDa,DD,DDe,DG04,DG05a,DG05b,DG05c,DG05d,DG06,DG07a,DG07b,DdI07a,
DdI07b,DL,DL97a,DL97b, DL00,DL04,DL07,DLS,DP02,DP,DVA,Dyu,Fu,Ger82,Gru,Gru07,
GI07a,GI07b,GI03,GPT01,GPT02a,GPT02,GPT03,GPT04,GPT05,GT07,LR,Lop,Lop01,
MS,MY,MZ,M05a,M05b,M05c,PT04,PT07a,PT07b,PR,RT,Ros,SVI,Tirao,
Yak,YM02,YMP01,YMP02,Zyg}; see also \cite{BB83}.

There is very little on the subject in monographs
--- the more classical ones (e.g.,
\cite{Chi,FrB,Ger77,Szb}) predate most of the subject; see,
however, Atkinson \cite[Section~6.6]{Atk}. Ismail \cite {IsmBk} has no
discussion and Simon \cite{S,S2} has a single section! Because of
the use of MOPRL in \cite{DKSppt}, we became interested in the
subject and, in particular, we needed some basic results for that
paper which we couldn't find in the literature or which, at least,
weren't very accessible. Thus, we decided to produce this
comprehensive review that we hope others will find useful.

As with the scalar case, the subject breaks into two parts,
conveniently called the analytic theory (general structure
results) and the algebraic theory (the set of non-trivial
examples). This survey deals entirely with the analytic theory. We
note, however, that one of the striking developments in recent
years has been the discovery that there are rich classes of
genuinely new MOPRL, even at the classical level of Bochner's
theorem; see \cite{CG,DG04,DuLo4,DuLo5,Gru,GPT02,GPT03,GPT04,GPT05,GT07,PT07a,RT}
and the forthcoming monograph \cite{DGPT} for further
discussion of this algebraic side.

In this introduction, we will focus mainly on the MOPRL case. For
scalar OPRL, a key issue
is the passage from measure to monic OPRL, then to normalized OPRL,
and finally to
Jacobi parameters. There are no choices in going from measure to
monic OP, $P_n(x)$.
They are determined by
\begin{equation} \lb{1.1}
P_n(x) = x^n + \text{lower order}, \qquad \jap{x^j, P_n}=0 \quad j=1,
\dots, n-1.
\end{equation}

However, the basic condition on the orthonormal polynomials, namely,
\begin{equation} \lb{1.2}
\jap{p_n,p_m} =\delta_{nm}
\end{equation}
does not uniquely determine the $p_n(x)$. The standard choice is
$$
p_n(x) =  \f{P_n(x)}{\norm{P_n}}.
$$
However, if
$\theta_0, \theta_1, \dots$ are arbitrary real numbers,
then
\begin{equation} \lb{1.3}
\ti p_n(x) =  \f{e^{i\theta_n} P_n(x)}{\norm{P_n}}
\end{equation}
also obey \eqref{1.2}. If the recursion coefficients (aka Jacobi
parameters), are defined via
\begin{equation} \lb{1.4}
xp_n =a_{n+1} p_{n+1} + b_{n+1} p_n + a_n p_{n-1},
\end{equation}
then the choice \eqref{1.3} leads to
\begin{equation} \lb{1.5}
\ti b_n =b_n, \qquad \ti a_n =e^{i\theta_n} a_n e^{-i\theta_{n-1}}.
\end{equation}

The standard choice is, of course, most natural here; for example, if
\begin{equation} \lb{1.6}
p_n(x) = \kappa_n x^n + \text{lower order},
\end{equation}
then $a_n>0$ implies $\kappa_n >0$. It would be crazy to make any
other choice.

For MOPRL, these choices are less clear. As we will explain in
Section~\ref{s1.2}, there are
now two matrix-valued ``inner products" formally written as
\begin{align}
\ang{f,g}_R &= \int f(x)^\dagger\, d\mu(x) g(x), \lb{1.7} \\
\ang{f,g}_L &= \int g(x)\, d\mu(x) f (x)^\dagger, \lb{1.8}
\end{align}
where now $\mu$ is a matrix-valued measure and  ${}^\dagger$
denotes the adjoint, and corresponding two sets of monic OPRL:
$P_n^R(x)$ and $P_n^L(x)$. The orthonormal polynomials are
required to obey
\begin{equation} \lb{1.9}
\ang{p_n^R,p_m^R}_R = \delta_{nm}\bdone.
\end{equation}
The analogue of \eqref{1.3} is
\begin{equation} \lb{1.10z}
\ti p_n^R(x) = P_n^R(x) \ang{P_n^R,P_n^R}^{-1/2} \sigma_n
\end{equation}
for a unitary $\sigma_n$. For the immediately following, use $p_n^R$
to be
the choice $\sigma_n\equiv 1$. For any such choice, we have a
recursion relation,
\begin{equation} \lb{1.10x}
xp_n^R(x) =p_{n+1}^R(x) A_{n+1}^\dagger + p_n^R(x) B_{n+1} + p_{n-1}^R
(x) A_n
\end{equation}
with the analogue of \eqref{1.5} (comparing $\sigma_n\equiv
\bdone$ to general $\sigma_n$)
\begin{equation} \lb{1.10}
\ti B_n =\sigma_n^\dagger B_n \sigma_n \qquad \ti A_n =\sigma_{n-1}^
\dagger A_n\sigma_n.
\end{equation}

The obvious analogue of the scalar case is to pick $\sigma_n\equiv
\bdone$, which makes $\kappa_n$ in
\begin{equation} \lb{1.11x}
p_n^R(x) = \kappa_n x^n + \text{lower order}
\end{equation}
obey $\kappa_n >0$. Note that \eqref{1.10x} implies
\begin{equation} \lb{1.11}
\kappa_n = \kappa_{n+1} A_{n+1}^\dagger
\end{equation}
or, inductively,
\begin{equation} \lb{1.13}
\kappa_n = (A_n^\dagger \dots A_1^\dagger)^{-1}.
\end{equation}
In general, this choice does not lead to $A_n$ positive or even
Hermitian. Alternatively, one can pick $\sigma_n$ so $\ti A_n$ is
positive. Besides these two ``obvious" choices, $\kappa_n > 0$ or
$A_n >0$, there is a third that $A_n$ be lower triangular that, as
we will see in Section~\ref{s1.4}, is natural. Thus, in the study
of MOPRL one needs to talk about equivalent sets of $p_n^R$ and of
Jacobi parameters, and this is a major theme of Chapter~\ref{s2}.
Interestingly enough for MOPUC, the commonly picked choice
equivalent to $A_n>0$ (namely, $\rho_n >0$) seems to suffice for
applications. So we do not discuss equivalence classes for MOPUC.

Associated to a set of matrix Jacobi parameters is a block Jacobi
matrix, that is, a matrix which
when written in $l\times l$ blocks is tridiagonal; see \eqref{b6.1}
below.

In Chapter~\ref{s2}, we discuss the basics of MOPRL while Chapter~\ref
{s3} discusses MOPUC.
Chapter~\ref{s4} discusses the Szeg\H{o} mapping connection of MOPUC
and MOPRL. Finally,
Chapter~\ref{s5} discusses the extension of the theory of regular OPs
\cite{StT} to MOPRL.

While this is mainly a survey, it does have numerous new results, of
which we mention:
\begin{SL}
\item[(a)] The clarification of equivalent Jacobi parameters and
several new theorems (Theorems~\ref{l.b4} and \ref{l.b5}).
\item[(b)] A new result (Theorem~\ref{T2.18H}) on the order of
poles or zeros of $m(z)$ in terms of eigenvalues of $J$ and the
once stripped $J^{(1)}$. \item[(c)] Formulas for the resolvent in
the MOPRL (Theorem~\ref {T2.27}) and MOPUC (Theorem~\ref{T3.21})
cases. \item[(d)] A theorem on zeros of $\det(\Phi_n^R)$
(Theorem~\ref{T3.6A}) and eigenvalues of a cutoff CMV matrix
(Theorem~\ref{T3.6D}). \item[(e)] A new proof of the Geronimus
relations (Theorem~\ref{T4.2}). \item[(f)] Discussion of regular
MOPRL (Chapter~\ref{s5}).
\end{SL}

\smallskip
There are numerous open questions and conjectures in this paper, of
which we mention:
\begin{SL}
\item[(1)] We prove that type~1 and type~3 Jacobi parameters in
the Nevai class have $A_n\to\bdone$ but do not know if this is
true for type~2 and, if so, how to prove it.

\item[(2)] Determine which monic matrix polynomials, $\Phi$, can
occur as monic MOPUC.
We know $\det(\Phi(z))$ must have all of its zeros in the unit
disk in $\C$, but unlike the scalar
case where this is sufficient, we do not know necessary and
sufficient conditions.

\item[(3)] Generalize Khrushchev theory
\cite{Kh2000,Khr,KhGo} to MOPUC; see Section~\ref{s3.13}.

\item[(4)] Provide a proof of Geronimus relations for MOPUC that
uses the theory of canonical moments \cite{DS}; see the discussion
at the start of Chapter~\ref{s4}.

\item[(5)] Prove Conjecture~\ref{Con5.7} extending a result of
Stahl--Totik \cite{StT} from OPRL to MOPRL.
\end{SL}

\subsection{Matrix-Valued Measures} \lb{s1.2}

Let $\calM_l$ denote the ring of all $l\times l$ complex-valued
matrices; we denote by
$\alpha^\dagger$ the Hermitian conjugate of $\alpha\in \calM_l$.
(Because of the use of
$\,^*\,$ for Szeg\H{o} dual in the theory of OPUC, we do not use it
for adjoint.)
For $\alpha\in\calM_l$, we denote by $\norm{\alpha}$ its Euclidean
norm  (i.e., the norm of $\alpha$ as a linear operator on $\C^l$ with the
usual Euclidean norm).
Consider the
set $\pol$ of all polynomials in $z\in\C$ with coefficients from $
\calM_l$. The set $\pol$
can be considered either as a right or as a left module over $\calM_l
$; clearly, conjugation
makes the left and right structures isomorphic. For $n=0,1,\dots$, $
\pol_n$ will denote those
polynomials in $\pol$ of degree at most $n$. The set $\mathcal{V}$
denotes the set of all
polynomials in $z \in \C$ with coefficients from $\C^l$. The standard
inner product in
$\C^l$ is denoted by $\langle\cdot, \cdot\rangle_{\C^l}$.

A matrix-valued measure, $\mu$, on $\R$ (or $\C$) is the
assignment of a positive semi-definite $l\times l$ matrix $\mu(X)$
to every Borel set $X$ which is countably additive. We will
usually normalize it by requiring
\begin{equation} \lb{1.20}
\mu(\R)=\bdone
\end{equation}
(or $\mu(\C)=\bdone$) where $\bdone$ is the $l\times l$ identity
matrix. (We use $\bdone$ in general for an identity operator,
whether in $\calM_l$ or in the operators on some other Hilbert
space, and $\bdzero$ for the zero operator or matrix.) Normally,
our measures for MOPRL will have compact support and, of course,
our measures for MOPUC will be supported on all or part of
$\partial\D$ ($\D$ is the unit disk in $\C$).

Associated to any such measures is a scalar measure
\begin{equation} \lb{1.21}
\mu_\tr(X) =\Tr(\mu(X))
\end{equation}
(the trace normalized by $\Tr(\bdone)=l$). $\mu_\tr$ is normalized
by $\mu_\tr(\R)=l$.

Applying the Radon--Nikodym theorem to the matrix elements of $\mu$,
we see there is a
positive semi-definite matrix function $M_{ij}(x)$ so
\begin{equation} \lb{1.22}
d\mu_{ij}(x)=M_{ij}(x)\, d\mu_\tr (x).
\end{equation}
Clearly, by \eqref{1.21},
\begin{equation} \lb{1.23}
\Tr(M(x)) =1
\end{equation}
for $d\mu_\tr$-a.e.\ $x$. Conversely, any scalar measure with
$\mu_\tr (\R)=l$ and positive semi-definite matrix-valued function
$M$ obeying \eqref{1.23} define a matrix-valued measure normalized
by \eqref{1.20}.

Given $l\times l$ matrix-valued functions $f,g$, we define the $l
\times l$ matrix $\ang{f,g}_R$ by
\begin{equation} \lb{1.24}
\ang{f,g}_R=\int f(x)^\dagger M(x)g(x)\, d\mu_\tr (x),
\end{equation}
that is, its $(j,k)$ entry is
\begin{equation} \lb{1.25}
\sum_{nm} \int \ol{f_{nj}(x)}\, M_{nm}(x) g_{mk}(x)\, d\mu_\tr(x).
\end{equation}
Since $f^\dagger Mf\geq 0$, we see that
\begin{equation} \lb{1.25a}
\ang{f,f}_R \geq 0.
\end{equation}
One might be tempted to think of $\ang{f,f}_R^{1/2}$ as some kind of
norm, but that is
doubtful. Even if $\mu$ is supported at a single point, $x_0$, with
$M=l^{-1}\bdone$, this
``norm'' is essentially the absolute value of $A=f(x_0)$, which is
known not to obey the
triangle inequality! (See \cite[Sect.~I.1]{S-TI} for an example.)

However, if one looks at
\begin{equation} \lb{1.26}
\norm{f}_R=(\Tr\ang{f,f}_R)^{1/2},
\end{equation}
one does have a norm (or, at least, a semi-norm). Indeed,
\begin{equation} \lb{1.27}
\langle f,g\rangle_R=\Tr\ang{f,g}_R
\end{equation}
is a sesquilinear form which is positive semi-definite, so
\eqref{1.26} is the semi-norm corresponding to an inner product
and, of course, one has a Cauchy--Schwarz inequality
\begin{equation} \lb{1.28}
\abs{\Tr\ang{f,g}_R}\leq \norm{f}_R \norm{g}_R.
\end{equation}

We have not specified which $f$'s and $g$'s can be used in \eqref
{1.24}. We have in mind
mainly polynomials in $x$ in the real case and Laurent polynomials in
$z$ in the $\partial\bbD$
case although, obviously, continuous functions are okay. Indeed, it
suffices that $f$ (and $g$) be measurable and
obey
\begin{equation} \lb{1.28a}
\int \Tr(f^\dagger (x) f(x))\, d\mu_\tr (x) <\infty
\end{equation}
for the integrals in \eqref{1.25} to converge. The set of
equivalence classes under $f\sim g$ if $\norm{f-g}_R = 0$ defines
a Hilbert space, $\calH$, and $\langle f,g\rangle_R$ is
the inner product on this space.

Instead of \eqref{1.24}, we use the suggestive shorthand
\begin{equation} \lb{1.29}
\ang{f,g}_R = \int f(x)^\dagger \, d\mu(x) g(x).
\end{equation}
The use of $R$ here comes from ``right'' for if $\alpha\in\calM_l$,
\begin{align}
\ang{f,g\alpha}_R &= \ang{f,g}_R \alpha, \lb{1.30} \\
\ang{f\alpha,g}_R &= \alpha^\dagger\ang{f,g}_R, \lb{1.31}
\end{align}
but, in general, $\ang{f,\alpha g}_R$ is not related to $\ang{f,g}_R$.

While $(\Tr\ang{f,f}_R)^{1/2}$ is a natural analogue of the norm
in the scalar case, it will sometimes be useful to instead
consider
\begin{equation} \lb{1.31a}
[\det \ang{f,f}_R]^{1/2}.
\end{equation}
Indeed, this is a stronger ``norm'' in that $\det >0\Rightarrow
\Tr>0$ but not vice-versa.

When $d\mu$ is a ``direct sum,'' that is, each $M(x)$ is diagonal,
one can appreciate the difference.
In that case, $d\mu=d\mu_1\oplus\cdots\oplus d\mu_l$ and the MOPRL
are direct sums (i.e., diagonal
matrices) of scalar OPRL
\begin{equation} \lb{1.31b}
P_n^R (x,d\mu) =P_n(x,d\mu_1) \oplus\cdots\oplus P_n(x,d\mu_l).
\end{equation}
Then
\begin{equation} \lb{1.31c}
\norm{P_n^R}_R =\biggl(\, \sum_{j=1}^l \norm{P_n(\cdot, d\mu_j)}_{L^2
(d\mu_j)}^2\biggr)^{1/2},
\end{equation}
while
\begin{equation} \lb{1.31d}
(\det \ang{P_n^R,P_n^R}_R)^{1/2} = \prod_{j=1}^l \norm{P_n (\cdot, d
\mu_j)}_{L^2(d\mu_j)}.
\end{equation}
In particular, in terms of extending the theory of regular
measures \cite{StT}, $\norm{P_n^R}^{1/n}_R$ is only sensitive to
$\max\norm{P_n (\cdot, d\mu_j)}_{L^2(d\mu_j)}^{1/2} $ while
$(\det\ang{P_n^R,P_n^R}_R)^{1/2}$ is sensitive to them all. Thus,
$\det$ will be needed for that theory (see Chapter~\ref{s5}).

There will also be a left inner product and, correspondingly, two
sets of MOPRL and MOPUC. We discuss
this further in Sections~\ref{s2.1} and \ref{s3.1}.

Occasionally, for $\bbC^l$ vector-valued functions $f$ and $g$, we
will want to consider the scalar
\begin{equation} \lb{1.32}
\sum_{k,j}\int \ol{f_k(x)}\, M_{kj}(x) g_j (x)\, d\mu_\tr(x),
\end{equation}
which we will denote
\begin{equation} \lb{1.33}
\int d \langle f(x),\mu(x) g(x)\rangle_{\C^l}.
\end{equation}

We next turn to module Fourier expansions. A set $\{\varphi_j\}_{j=1}
^N$ in $\calH$ ($N$ may be infinite)
is called orthonormal if and only if
\begin{equation} \lb{1.33a}
\ang{\varphi_j,\varphi_k}_R =\delta_{jk}\bdone.
\end{equation}
This natural terminology is an abuse of notation since
\eqref{1.33a} implies orthogonality in $\langle\cdot,\cdot
\rangle_R$ but not normalization, and is much stronger than orthogonality in
$\langle \cdot,\cdot\rangle_R$.

Suppose for a moment that $N<\infty$. For any $a_1, \dots, a_N\in
\calM_l$, we can form $\sum_{j=1}^N
\varphi_j a_j$ and, by the right multiplication relations \eqref
{1.30}, \eqref{1.31}, and \eqref{1.33a},
we have
\begin{equation} \lb{1.33b}
\biggl< \!\!\!\biggl< \sum_{j=1}^N \varphi_j a_j, \sum_{j=1}^N
\varphi_j b_j \biggr>\!\! \!\biggr>_R =
\sum_{j=1}^N a_j^\dagger b_j.
\end{equation}
We will denote the set of all such $\sum_{j=1}^N \varphi_j a_j$ by $
\calH_{(\varphi_j)}$---it is a
vector subspace of $\calH$ of dimension (over $\bbC$) $Nl^2$.

Define for $f\in\calH$,
\begin{equation} \lb{1.33c}
\pi_{(\varphi_j)} (f)=\sum_{j=1}^N \varphi_j \ang{\varphi_j,f}_R.
\end{equation}
It is easy to see it is the orthogonal projection in the scalar
inner product $\langle \cdot,\cdot\rangle_R$ from $\calH$ to $\calH_{(\varphi_j)}$.

By the standard Hilbert space calculation, taking care to only
multiply on the right, one finds the
Pythagorean theorem,
\begin{equation} \lb{1.33d}
\ang{f,f}_R = \ang{f-\pi_{(\varphi_j)} f, f-\pi_{(\varphi_j)}f}_R +
\sum_{j=1}^N
\ang{\varphi_j,f}_R^\dagger \ang{\varphi_j,f}_R.
\end{equation}

As usual, this proves for infinite $N$ that
\begin{equation} \lb{1.33e}
\sum_{j=1}^N \ang{\varphi_j,f}_R^\dagger \ang{\varphi_j,f}_R \leq \ang
{f,f}_R
\end{equation}
and the convergence of
\begin{equation} \lb{1.33f}
\sum_{j=1}^N \varphi_j \ang{\varphi_j,f}_R \equiv \pi_{(\varphi_j)}(f)
\end{equation}
allowing the definition of $\pi_{(\varphi_j)}$ and of $\calH_
{(\varphi_j)}\equiv\Ran \pi_{(\varphi_j)}$
for $N=\infty$.

An orthonormal set is called complete if
$\calH_{(\varphi_j)}=\calH$. In that case, equality holds in
\eqref{1.33e} and $\pi_{(\varphi_j)}(f)=f$.

For orthonormal bases, we have the Parseval relation from
\eqref{1.33d}
\begin{equation} \lb{1.33g}
\ang{f,f}_R =\sum_{j=1}^\infty \ang{\varphi_j,f}_R^\dagger \ang
{\varphi_j,f}_R
\end{equation}
and
\begin{equation}\lb{1.33h}
\norm{f}_R^2 =\sum_{j=1}^\infty \Tr(\ang{\varphi_j,f}_R^\dagger
\ang {\varphi_j,f}_R).
\end{equation}

\subsection{Matrix M\"obius Transformations} \lb{s1.3}

Without an understanding of matrix M\"obius transformations, the
form of the MOPUC Geronimus theorem we will prove in
Section~\ref{s3.5} will seem strange-looking. To set the stage,
recall that scalar fractional linear transformations (FLT) are
associated to matrices $T=\left( \begin{smallmatrix} a&b \\ c&d
\end{smallmatrix}\right)$ with $\det T \not= 0$ via
\begin{equation} \lb{1.20x}
f_T(z) = \f{az+b}{cz+d}\, .
\end{equation}
Without loss, one can restrict to
\begin{equation} \lb{1.21x}
\det(T)=1.
\end{equation}
Indeed, $T\mapsto f_T$ is a $2$ to $1$ map of $\bbS\bbL(2,\bbC)$ to
maps of $\bbC\cup\{\infty\}$ to itself.
One advantage of the matrix formalism is that the map is a matrix
homomorphism, that is,
\begin{equation} \lb{1.22x}
f_{T\circ S} = f_T \circ f_S,
\end{equation}
which shows that the group of FLTs is $\bbS\bbL(2,\bbC)/\{\bdone, -
\bdone\}$.

While \eqref{1.22x} can be checked by direct calculation, a more
instructive way is to look at the
complex projective line. $u,v\in\bbC^2\setminus\{0\}$ are called
equivalent if there is $\lambda\in
\bbC\setminus\{0\}$ so that $u=\lambda v$. Let $[\cdot]$ denote
equivalence classes. Except for
$[\binom{1}{0}]$, every equivalence class contains exactly one point
of the form $\binom{z}{1}$ with
$z\in\bbC$. If $[\binom{1}{0}]$ is associated with $\infty$, the set
of equivalence classes is naturally
associated with $\bbC\cup\{\infty\}$. $f_T$ then obeys
\begin{equation} \lb{1.23x}
\biggl[ T\binom{z}{1}\biggr] = \biggl[ \binom{f_T(z)}{1}\biggr]
\end{equation}
from which \eqref{1.22x} is immediate.

By M\"obius transformations we will mean those FLTs that map
$\bbD$ onto itself. Let
\begin{equation} \lb{1.24x}
J=\left( \begin{array}{rr}
1 & 0 \\
0 & -1
\end{array} \right).
\end{equation}
Then $[u]=[\binom{z}{1}]$ with $\abs{z}=1$ (resp.\ $\abs{z}<1$) if
and only if $\jap{u,Ju}=0$
(resp.\ $\jap{u,Ju}<0$). From this, it is not hard to show that if $
\det(T)=1$, then $f_T$ maps
$\bbD$ invertibly onto $\bbD$ if and only if
\begin{equation} \lb{1.25x}
T^\dagger \! JT=J.
\end{equation}
If $T$ has the form $\left(\begin{smallmatrix} a&b \\ c&d \end
{smallmatrix}\right)$, this is equivalent to
\begin{equation} \lb{1.26x}
\abs{a}^2 - \abs{c}^2 =1, \qquad \abs{b}^2-\abs{d}^2 =-1, \qquad \bar
ab-\bar cd =0.
\end{equation}
The set of $T$'s obeying $\det(T)=1$ and \eqref{1.25x} is called $\bbS
\bbU(1,1)$. It is
studied extensively in \cite[Sect.~10.4]{S2}.

The self-adjoint elements of $\bbS\bbU(1,1)$ are parametrized by $
\alpha\in\bbD$ via $\rho =
(1-\abs{\alpha}^2)^{1/2}$,
\begin{equation} \lb{1.27x}
T_\alpha = \f{1}{\rho}\begin{pmatrix}
1 & \alpha \\
\bar\alpha & 1 \end{pmatrix}
\end{equation}
associated to
\begin{equation} \lb{1.28x}
f_{T_\alpha}(z) = \f{z+\alpha}{1+\bar\alpha z}\, .
\end{equation}
Notice that
\begin{equation} \lb{1.29x}
T_\alpha^{-1} = T_{-\alpha}
\end{equation}
and that
\[
\forall z\in\bbD,\, \exists \, ! \, \alpha \text{ such that } T_
\alpha (0)=z,
\]
namely, $\alpha =z$.

It is a basic theorem that every holomorphic bijection of $\bbD$ to $
\bbD$ is an $f_T$ for some
$T$ in $\bbS\bbU(1,1)$ (unique up to $\pm\bdone$).

With this in place, we can turn to the matrix case. Let $\calM_l$ be
the space of $l\times l$
complex matrices with the Euclidean norm induced by the vector norm
$\jap{\cdot, \cdot}_{\bbC^l}^{1/2}$. Let
\begin{equation} \lb{1.30x}
\bbD_l =\{A\in\calM_l\, \colon \norm{A} <1\}.
\end{equation}
We are interested in holomorphic bijections of $\bbD_l$ to itself,
especially via a suitable notion of
FLT. There is a huge (and diffuse) literature on the subject,
starting with its use in analytic number
theory. It has also been studied in connection with electrical
engineering filters and indefinite
matrix Hilbert spaces. Among the huge literature, we mention \cite
{ALK,Alp,Dym,GK,Helg,ScZ}. Especially
relevant to MOPUC is the book of Bakonyi--Constantinescu \cite{BakCon}.

Consider $\calM_l \oplus \calM_l = \calM_l [2]$ as a right module
over $\calM_l$. The $\calM_l$-projective line is defined by saying
$\left[\begin{smallmatrix} X \\ Y
\end{smallmatrix}\right]\sim \left[\begin{smallmatrix} X' \\
Y'\end{smallmatrix}\right]$, both in $\calM_l[2]
\setminus\{\bdzero\}$, if and only if there exists
$\Lambda\in\calM_l$, $\Lambda$ invertible so that
\begin{equation} \lb{1.31x}
X= X'\Lambda, \qquad Y=Y'\Lambda.
\end{equation}
Let $T$ be a map of $\calM_l[2]$ of the form
\begin{equation} \lb{1.32x}
T=\begin{pmatrix} A&B \\ C&D
\end{pmatrix}
\end{equation}
acting on $\calM_l[2]$ by
\begin{equation} \lb{1.33x}
T\begin{bmatrix} X \\ Y \end{bmatrix} =
\begin{bmatrix} AX + BY \\ CX + DY \end{bmatrix}.
\end{equation}
Because this acts on the left and $\Lambda$ equivalence on the right,
$T$ maps equivalence classes
to themselves. In particular, if $CX+D$ is invertible, $T$ maps the
equivalence class of
$\left[\begin{smallmatrix}X \\ \bdone
\end{smallmatrix}\right]$ to the equivalence class of
$\left[\begin{smallmatrix}f_T[X] \\ \bdone \end{smallmatrix}\right]$,
where
\begin{equation} \lb{1.34}
f_T[X]=(AX+B)(CX+D)^{-1}.
\end{equation}

So long as $CX+D$ remains invertible, \eqref{1.22x} remains true. Let
$J$ be the $2l\times 2l$ matrix
in $l\times l$ block form
\begin{equation} \lb{1.35}
J = \left( \begin{array}{rr}
\bdone  & \bdzero \\ \bdzero & -\bdone
\end{array} \right).
\end{equation}
Note that (with $\left[\begin{smallmatrix} X \\ \bdone \end
{smallmatrix}\right]^\dagger = [X^\dagger \bdone]$)
\begin{equation} \lb{1.36}
\begin{bmatrix} X \\ \bdone \end{bmatrix}^\dagger J \begin{bmatrix} X
\\ \bdone \end{bmatrix} \leq
\bdzero \Leftrightarrow X^\dagger X\leq \bdone \Leftrightarrow \norm{X}
\leq \bdone.
\end{equation}
Therefore, if we define $\bbS\bbU(l,l)$ to be those $T$'s with
$\det T=1$ and
\begin{equation} \lb{1.37}
T^\dagger \! JT =J,
\end{equation}
then
\begin{equation} \lb{1.38}
T\in\bbS\bbU(l,l) \Rightarrow f_T [\bbD_l]=\bbD_l \text{ as a
bijection}.
\end{equation}
If $T$ has the form \eqref{1.32x}, then \eqref{1.37} is equivalent to
\begin{align}
& A^\dagger A - C^\dagger C = D^\dagger D - B^\dagger B =\bdone, \lb
{1.61a} \\
& A^\dagger B = C^\dagger D \lb{1.61b}
\end{align}
(the fourth relation $B^\dagger A = D^\dagger C$ is equivalent to
\eqref{1.61b}).

This depends on

\begin{proposition}\lb{P1.3.1} If $T=\left(\begin{smallmatrix} A&B \\
C&D\end{smallmatrix}\right)$ obeys \eqref{1.37} and $\norm{X} <1$,
then $CX+D$ is invertible.
\end{proposition}

\begin{proof} \eqref{1.37} implies that
\begin{align}
T^{-1} &= JT^\dagger \! J \lb{1.37a} \\
&= \left( \begin{array}{rr}
A^\dagger & -C^\dagger \\
-B^\dagger & D^\dagger
\end{array} \right). \lb{1.37b}
\end{align}
Clearly, \eqref{1.37} also implies $T^{-1}\in\bbS\bbU(l,l)$. Thus,
by \eqref{1.61a} for $T^{-1}$,
\begin{equation} \lb{1.37c}
DD^\dagger -CC^\dagger =\bdone.
\end{equation}
This implies first that $DD^\dagger\geq \bdone$, so $D$ is
invertible, and second that
\begin{equation} \lb{1.37d}
\norm{D^{-1} C}\leq 1.
\end{equation}
Thus, $\norm{X}<1$ implies $\norm{D^{-1}CX}<1$ so $\bdone+D^{-1}CX$
is invertible, and thus
so is $D(\bdone+D^{-1}CX)$.
\end{proof}

It is a basic result of Cartan \cite{Cartan} (see Helgason
\cite{Helg} and the discussion therein) that

\begin{theorem}\lb{T1.3.2} A holomorphic bijection, $g$, of $\bbD_l$
to itself is either of
the form
\begin{equation} \lb{1.38x}
g(X)=f_T(X)
\end{equation}
for some $T\in\bbS\bbU(l,l)$ or
\begin{equation} \lb{1.39}
g(X)=f_T(X^t).
\end{equation}
\end{theorem}

Given $\alpha\in\calM_l$ with $\norm{\alpha}<1$, define
\begin{equation} \lb{1.40}
\rho^L = (\bdone-\alpha^\dagger \alpha)^{1/2}, \qquad
\rho^R =(\bdone-\alpha\alpha^\dagger)^{1/2}.
\end{equation}

\begin{lemma} \lb{L1.3.3} We have
\begin{alignat}{2}
\alpha\rho^L &= \rho^R\alpha, \qquad & \alpha^\dagger \rho^R &=
\rho^L \alpha^\dagger, \lb{1.41} \\
\alpha(\rho^L)^{-1} &= (\rho^R)^{-1}\alpha, \qquad & \alpha^\dagger
(\rho^R)^{-1} &= (\rho^L)^{-1}
\alpha^\dagger. \lb{1.42}
\end{alignat}
\end{lemma}

\begin{proof} Let $f$ be analytic in $\bbD$ with $f(z)=\sum_{n=0}^
\infty c_n z^n$ its Taylor series at $z=0$.
Since $\norm{\alpha^\dagger\alpha}<1$, we have
\begin{equation} \lb{1.43}
f(\alpha^\dagger \alpha) = \sum_{n=0}^\infty c_n (\alpha^\dagger
\alpha)^n
\end{equation}
norm convergent, so $\alpha (\alpha^\dagger\alpha)^n = (\alpha\alpha^
\dagger)^n\alpha$ implies
\begin{equation} \lb{1.44}
\alpha f(\alpha^\dagger\alpha) = f(\alpha\alpha^\dagger)\alpha,
\end{equation}
which implies the first halves of \eqref{1.41} and \eqref{1.42}. The
other halves follow by
taking adjoints.
\end{proof}

\begin{theorem}\lb{T1.4} There is a one-one correspondence between $
\alpha$'s in $\calM_l$ obeying
$\norm{\alpha}<1$ and positive self-adjoint elements of $\bbS\bbU(l,l)
$ via
\begin{equation} \lb{1.45}
T_\alpha = \begin{pmatrix}
(\rho^R)^{-1} & (\rho^R)^{-1}\alpha \\
(\rho^L)^{-1} \alpha^\dagger & (\rho^L)^{-1}
\end{pmatrix}.
\end{equation}
\end{theorem}

\begin{proof} A straightforward calculation using Lemma~\ref{L1.3.3}
proves that $T_\alpha$ is self-adjoint and $T_\alpha^\dagger
JT_\alpha =J$. Conversely, if $T$ is self-adjoint, $T=\left(
\begin{smallmatrix} A&B\\C&D\end{smallmatrix} \right)$ and in
$\bbS\bbU(l,l)$, then $T^\dagger = T\Rightarrow A^\dagger = A$,
$B^\dagger =C$, so \eqref{1.61a} becomes
\begin{equation} \lb{1.46}
AA^\dagger - BB^\dagger =\bdone,
\end{equation}
so if
\begin{equation} \lb{1.47}
\alpha = A^{-1} B,
\end{equation}
then \eqref{1.46} becomes
\begin{equation} \lb{1.48}
A^{-1} (A^{-1})^\dagger + \alpha\alpha^\dagger =\bdone.
\end{equation}
Since $T\geq 0$, $A\geq 0$ so \eqref{1.48} implies $A=(\rho^R)^{-1}$,
and then \eqref{1.47} implies
$B=(\rho^R)^{-1} \alpha$.

By Lemma~\ref{L1.3.3},
\begin{equation} \lb{1.49}
C=B^\dagger = \alpha^\dagger (\rho^R)^{-1} = (\rho^L)^{-1} \alpha^
\dagger
\end{equation}
and then (by $D=D^\dagger$, $C^\dagger =B$, and \eqref{1.61a})
$D D^\dagger - C C^\dagger =\bdone$ plus $D>0$ implies
$D= (\rho^L)^{-1}$.
\end{proof}

\begin{corollary}\lb{C1.3.5} For each $\alpha\in\bbD_l$, the map
\begin{equation} \lb{1.50}
f_{T_\alpha}(X) = (\rho^R)^{-1} (X+\alpha) (\bdone+\alpha^\dagger X)^
{-1} (\rho^L)
\end{equation}
takes $\bbD_l$ to $\bbD_l$. Its inverse is given by
\begin{equation} \lb{1.51}
f_{T_\alpha}^{-1}(X) = f_{T_{-\alpha}}(X) = (\rho^R)^{-1} (X-\alpha)
(\bdone-\alpha^\dagger X)^{-1} (\rho^L).
\end{equation}
\end{corollary}

There is an alternate form for the right side of \eqref{1.50}.

\begin{proposition} \lb{P1.6} The following identity holds true for any $X$,
$\norm{X}\leq 1$:
\begin{equation} \lb{1.82a}
\rho^R (1+X\alpha^\dagger)^{-1} (X+\alpha)(\rho^L)^{-1} =
(\rho^R)^{-1} (X+\alpha)(1+\alpha^\dagger X)^{-1}\rho^L.
\end{equation}
\end{proposition}

\begin{proof}
By the definition of $\rho^L$ and $\rho^R$, we have
$$
X (\rho^L)^{-2}(1-\alpha^\dagger\alpha)=(\rho^R)^{-2}(1-\alpha\alpha^\dagger)X.
$$
Expanding, using \eqref{1.42} and rearranging, we get
$$
X(\rho^L)^{-2}+\alpha(\rho^L)^{-2}\alpha^\dagger X
=
(\rho^R)^{-2}X+X\alpha^\dagger(\rho^R)^{-2}\alpha.
$$
Adding $\alpha(\rho^L)^{-2}+X(\rho^L)^{-2}\alpha^\dagger X$ to both sides and
using \eqref{1.42} again, we obtain
\begin{multline*}
X(\rho^L)^{-2}+\alpha(\rho^L)^{-2}+X(\rho^L)^{-2}\alpha^\dagger X
+\alpha(\rho^L)^{-2}\alpha^\dagger X
\\
=
(\rho^R)^{-2}X+(\rho^R)^{-2}\alpha+X\alpha^\dagger(\rho^R)^{-2}X
+X\alpha^\dagger(\rho^R)^{-2}\alpha,
\end{multline*}
which is the same as
$$
(X+\alpha)(\rho^L)^{-2}(1+\alpha^\dagger X)
=
(1+X\alpha^\dagger)(\rho^R)^{-2}(X+\alpha).
$$
Multiplying by $(1+X\alpha^\dagger)^{-1}$ and $(1+\alpha^\dagger X)^{-1}$, we get
$$
(1+X\alpha^\dagger)^{-1}(X+\alpha)(\rho^L)^{-2}
=
(\rho^R)^{-2}(X+\alpha)(1+\alpha^\dagger X)^{-1}
$$
and the statement follows.
\end{proof}

\subsection{Applications and Examples} \lb{s1.4}

There are a number of simple examples which show that beyond their
intrinsic mathematical interest, MOPRL
and MOPUC have wide application.

\subsubsection*{{\rm (a)} Jacobi matrices on a strip} Let $\Lambda
\subset\bbZ^\nu$ be a subset (perhaps infinite) of the
$\nu$-dimensional lattice $\bbZ^\nu$ and let $\ell^2 (\Lambda)$ be
square summable sequences indexed by $\Lambda$. Suppose a real
symmetric matrix $\alpha_{ij}$ is given for all $i,j\in\Lambda$
with $\alpha_{ij}=0$ unless $\abs{i-j}=1$ (nearest neighbors). Let
$\beta_i$ be a real sequence indexed by $i\in\Lambda$. Suppose
\begin{equation} \lb{1.40x}
\sup_{i,j}\, \abs{\alpha_{ij}} + \sup_i\, \abs{\beta_i} <\infty.
\end{equation}
Define a bounded operator, $J$, on $\ell^2(\Lambda)$ by
\begin{equation} \lb{1.41x}
(Ju)_i = \sum_j \alpha_{ij} u_j + \beta_i u_i.
\end{equation}
The sum is finite with at most $2\nu$ elements.

The special case $\Lambda =\{1,2,\dots\}$ with $b_i =\beta_i$, $a_i =
\alpha_{i,i+1}>0$
corresponds precisely to classical semi-infinite tridiagonal Jacobi
matrices.

Now consider the situation where $\Lambda' \subset\bbZ^{\nu-1}$ is a
finite set with $l$ elements
and
\begin{equation} \lb{1.42x}
\Lambda = \{j\in\bbZ^\nu \, \colon j_1\in \{1,2,\dots\};\, (j_2,
\dots j_\nu)\in\Lambda'\},
\end{equation}
a ``strip'' with cross-section $\Lambda'$. $J$ then has a block $l
\times l$ matrix Jacobi form
where $(\gamma,\delta\in\Lambda'$)
\begin{alignat}{2}
(B_i)_{\gamma\delta} &= b_{(i,\gamma)}, && \qquad (\gamma=\delta), \lb
{1.42y} \\
&= a_{(i,\gamma)(i,\delta)}, && \qquad (\gamma\neq\delta), \lb{1.42z} \\
(A_i)_{\gamma\delta} &= a_{(i,\gamma)(i+1,\delta)}. \lb{1.43x}
\end{alignat}
The nearest neighbor condition says $(A_i)_{\gamma\delta}=0$ if $
\gamma\neq\delta$. If
\begin{equation} \lb{1.44x}
a_{(i,\gamma)(i+1,\gamma)} >0
\end{equation}
for all $i,\gamma$, then $A_i$ is invertible and we have a block
Jacobi matrix of the
kind described in Section~\ref{s2.2} below.

By allowing general $A_i,B_i$, we obtain an obvious generalization of
this model---an interpretation of
general MOPRL.

Schr\"odinger operators on strips have been studied in part as
approximations to $\bbZ^\nu$; see \cite{CrS,GKM,KS88,Lac,MV,SB04}.
{}From this point of view, it is also natural to allow periodic
boundary conditions in the vertical directions. Furthermore, there
is closely related work on Schr\"odinger (and other) operators
with matrix-valued potentials; see, for example, \cite{BGMS03,
CG01,CG02,CG03,CG04, CGHL00, CGZ07,GS03, GS00, SB07}.

\subsubsection*{{\rm (b)} Two-sided Jacobi matrices} This example goes
back at least to Nikishin \cite{Nik}. Consider the
case $\nu=2$, $\Lambda'=\{0,1\}
\subset\bbZ$, and $\Lambda$ as above. Suppose \eqref{1.44x} holds,
and in addition,
\begin{align}
a_{(1,0)(1,1)} &>0, \lb{1.45x} \\
a_{(i,0)(i,1)} &=0,  \quad i=2,3,\dots . \lb{1.46x}
\end{align}
Then there are no links between the rungs of the ``ladder,'' $\{1,2,
\dots\}\times \{0,1\}$ except at the
end and the ladder can be unfolded to $\bbZ$! Thus, a two-sided
Jacobi matrix can be viewed as a special
kind of one-sided $2\times 2$ matrix Jacobi operator.

It is known that for two-sided Jacobi matrices, the spectral theory
is determined by the $2\times 2$
matrix
\begin{equation} \lb{1.47x}
d\mu = \begin{pmatrix}
d\mu_{00} & d\mu_{01} \\
d\mu_{10} & d\mu_{11}
\end{pmatrix},
\end{equation}
where $d\mu_{kl}$ is the measure with
\begin{equation} \lb{1.48x}
\jap{\delta_k, (J-\lambda)^{-1}\delta_l} = \int \f{d\mu_{kl}(x)}{x-
\lambda},
\end{equation}
but also that it is very difficult to specify exactly which $d\mu$
correspond to two-sided Jacobi
matrices.

This difficulty is illuminated by the theory of MOPRL. By Favard's
theorem (see Theorem~\ref{favard}),
every such $d\mu$ (given by \eqref{1.47x} and positive definite and
non-trivial in a sense we will
describe in Lemma~\ref{nondeg}) yields a unique block Jacobi matrix
with $A_j>0$ (positive definite).
This $d\mu$ comes from a two-sided Jacobi matrix if and only if
\begin{SL}
\item[(a)] $B_j$ is diagonal for $j=2,3,\dots$.
\item[(b)] $A_j$ is diagonal for $j=1,2,\dots$.
\item[(c)] $B_j$ has strictly positive off-diagonal elements.
\end{SL}
These are very complicated indirect conditions on $d\mu$!

\subsubsection*{{\rm (c)} Banded matrices} Classical Jacobi matrices
are semi-infinite symmetric
tridiagonal matrices, that is,
\begin{equation} \lb{1.49x}
J_{km}=0 \quad\text{if } \abs{k-m}>1
\end{equation}
with
\begin{equation} \lb{1.50x}
J_{km}>0 \quad\text{if } \abs{k-m}=1.
\end{equation}

A natural generalization are $(2l+1)$-diagonal symmetric matrices,
that is,
\begin{alignat}{2}
J_{km} &= 0 &  \quad \text{if } \abs{k-m} &>l \lb{1.51x},  \\
J_{km} &> 0 & \quad \text{if } \abs{k-m} &=l. \lb{1.52x}
\end{alignat}

Such a matrix can be partitioned into $l\times l$ blocks, which is
tridiagonal in block. The
conditions \eqref{1.51x} and \eqref{1.52x} are equivalent to $A_k\in
\calL$, the set of lower
triangular matrices; and conversely, $A_k\in\calL$, with $A_k,B_k$
real (and $B_k$ symmetric)
correspond precisely to such banded matrices. This is why we
introduce type~3 MOPRL.

Banded matrices correspond to certain higher-order difference
equations. Unlike the second-order
equation (which leads to tridiagonal matrices) where every equation
with positive coefficients is
equivalent via a variable rescaling to a symmetric matrix, only
certain higher-order difference
equations correspond to symmetric block Jacobi matrices.

\subsubsection*{{\rm (d)} Magic formula} In \cite{DKSppt}, Damanik,
Killip, and Simon studied
perturbations of Jacobi and CMV matrices with periodic Jacobi
parameters (or Verblunsky coefficients).
They proved that if $\Delta$ is the discriminant of a two-sided
periodic $J_0$, then a bounded
two-sided $J$ has $\Delta(J)=S^p + S^{-p}$ ($(Su)_n\equiv u_{n+1}$)
if and only if $J$ lies
in the isospectral torus of $J_0$. They call this the magic formula.

This allows the study of perturbations of the isospectral torus by
studying $\Delta(J)$ which is
a polynomial in $J$ of degree $p$, and so a $2p+1$ banded matrix.
Thus, the study of perturbations
of periodic problems is connected to perturbations of $S^p + S^{-p}$
as block Jacobi matrices.
Indeed, it was this connection that stimulated our interest in MOPRL,
and \cite{DKSppt} uses some
of our results here.

\subsubsection*{{\rm (e)} Vector-valued prediction theory} As noted
in Section~\ref{s1.1},
both prediction theory and filtering theory use OPUC and have natural
MOPUC settings that
motivated much of the MOPUC literature.

\section{Matrix Orthogonal Polynomials on the Real Line} \lb{s2}

\subsection{Preliminaries} \lb{s2.1}

OPRL are the most basic and developed of orthogonal polynomials,
and so this chapter on the matrix analogue is the most important
of this survey. We present the basic formulas, assuming enough
familiarity with the scalar case (see
\cite{Chi,FrB,S,Rice,Szb,Teschl}) that we do not need to explain
why the objects we define are important.

\subsubsection{Polynomials, Inner Products, Norms} \lb{s2.1.1}

Let $d\mu$ be an $l\times l$ matrix-valued Hermitian positive
semi-definite finite measure on $\R$
normalized by $\mu(\R)={\boldsymbol 1}\in\calM_l $.
We assume for simplicity that $\mu$ has a compact support.
However, many of the results below do not need the latter restriction and in fact can be
found in the literature for matrix-valued measures with unbounded support.

Define (as in
\eqref{1.24})
\begin{alignat*}{3}
\ang{f , g}_R &= \int f(x)^\dagger \, d\mu(x) \, g(x), \quad & \norm
{f}_R & =(\Tr \ang{f,f}_R)^{1/2},
\quad & f,g & \in\pol, \\
\ang{f , g}_L &= \int g(x)\, d\mu(x) \, f(x)^\dagger,
\quad & \norm{f}_L &=(\Tr \ang{f,f}_L)^{1/2}, \quad & f,g &\in \pol.
\end{alignat*}
Clearly, we have
\begin{alignat}{2}
\ang{f,g}_R^\dagger &=\ang{g,f}_R, \qquad
& \ang{f,g}_L^\dagger & =\ang{g,f}_L, \label{b1} \\
\ang{f,g}_L &=\ang{g^\dagger,f^\dagger}_R, \qquad
& \norm{f}_L &=\norm{f^\dagger}_R. \label{b2}
\end{alignat}

As noted in Section~\ref{s1.2}, we have the left and right analogues
of the Cauchy inequality
$$
\abs{\Tr \ang{f,g}_R}\leq \norm{f}_R\norm{g}_R,
\quad
\abs{\Tr \ang{f,g}_L}\leq \norm{f}_L\, \norm{g}_L.
$$
Thus, $\norm{\cdot}_R$ and $\norm{\cdot}_L$ are semi-norms in $\pol$.
Indeed, as noted
in Section~\ref{s1.2}, they are associated to an inner product. The
sets $\{ f \,\colon
\norm{f}_R = 0\}$ and $\{ f \, \colon \norm{f}_L = 0\}$ are linear
subspaces. Let $\pol_R$ be
the completion of $\pol / \{ f \, \colon \norm{f}_R = 0\}$ (viewed as
a right module over $\calM_l$)
with respect to the norm $\norm{\cdot}_R$. Similarly, let $\pol_L$ be
the completion of $\pol
/ \{ f \, \colon \norm{f}_L = 0\}$ (viewed as a left module) with
respect to the norm
$\norm{\cdot}_L$.

The set $\mathcal{V}$ defined in Section~\ref{s1.2} is a linear space.
Let us introduce a semi-norm in $\mathcal{V}$ by
\begin{equation}\label{f.vsn}
\snorm{f}=\biggl\{\int d\langle f(x), \mu(x) f(x)\rangle_{\C^l}\biggr
\}^{1/2}.
\end{equation}
Let $\mathcal{V}_0\subset \mathcal{V}$ be the linear subspace of all
polynomials such that $\snorm{f}=0$ and let $\mathcal{V}_\infty$ be
the completion of the quotient space $\mathcal{V}/ \mathcal{V}_0$
with respect to the norm $\snorm{\cdot}$.

\begin{lemma}\label{nondeg}
The following are equivalent:
\begin{SL}
\item[{\rm{(1)}}] $\norm{f}_R > 0$ for every non-zero $f\in \pol$.
\item[{\rm{(2)}}] For all $n$, the dimension in $\pol_R$ of the set
of all polynomials of degree at most $n$ is
$(n+1)l^2$.
\item[{\rm{(3)}}] $\norm{f}_L > 0$ for every non-zero $f\in \pol$.
\item[{\rm{(4)}}] For all $n$, the dimension in $\pol_L$ of the set
of all polynomials of degree at most $n$ is
$(n+1)l^2$.
\item[{\rm{(5)}}] For every non-zero $v \in \mathcal{V}$, we have that
$\snorm{v} \not= 0$.
\item[{\rm{(6)}}] For all $n$, the dimension in $\mathcal{V}_\infty$
of all vector-valued polynomials
of degree at most $n$ is $(n+1)l$.
\end{SL}
The measure $d\mu$ is called non-trivial if these equivalent
conditions hold.
\end{lemma}

\begin{remark} If $l=1$, these are equivalent to the usual non-
triviality
condition, that is, $\supp(\mu)$ is infinite. For
$l>1$, we cannot define triviality in this simple way, as can be
seen by looking at the direct sum of a trivial and non-trivial
measure. In that case, the measure is not non-trivial in the above
sense but its support is infinite.
\end{remark}

\begin{proof}
The equivalences  (1) $\Leftrightarrow$ (2),
(3) $\Leftrightarrow$ (4), and  (5) $\Leftrightarrow$ (6)
are immediate. The
equivalence (1) $\Leftrightarrow$ (3) follows from \eqref{b2}.
Let us prove the equivalence (1) $\Leftrightarrow$ (5).
Assume that (1) holds and let $v \in\mathcal{V}$ be non-zero.
Let $f\in\calM_l$ denote the matrix that has $v$ as its leftmost
column and that has zero columns otherwise. Then, $0 \not= \norm{f}^2_R = \Tr \ang{f,f}_R
=\snorm{v}^2$ and hence (5) holds. Now assume that (1) fails
and let $f \in \pol$ be non-zero with $\norm{f}_R = 0$. Then, at
least one of the column
vectors of $f$ is non-zero.
Suppose for simplicity that this is the first column and
denote this column vector by $v$.
Let $t\in\calM_l$ be the matrix $t_{ij}=\delta_{i1}\delta_{j1}$;
then we have
$$
\norm{f}_R=0 \Rightarrow \ang{f,f}_R=0 \Rightarrow
0=\Tr( t^*\ang{f,f}_R t) =\snorm{v}^2
$$
and hence (5) fails.
\end{proof}

Throughout the rest of this chapter, we assume the measure $d\mu$ to
be non-trivial.

\subsubsection{Monic Orthogonal Polynomials}

\begin{lemma}\label{l.b1} Let $d\mu$ be a non-trivial measure.
\begin{SL}
\item[{\rm{(i)}}] There exists a unique monic polynomial $P_n^R$  of
degree $n$, which minimizes
the norm $\norm{P_n^R}_R$.

\item[{\rm{(ii)}}] The polynomial $P_n^R$ can be equivalently defined
as the monic polynomial of
degree $n$ which satisfies
\begin{equation}
\ang{P_n^R, f}_R = {\boldsymbol 0} \quad \text{ for any } f\in \pol,
\quad \deg f<n.
\label{b2.2}
\end{equation}

\item[{\rm{(iii)}}] There exists a unique monic polynomial $P_n^L$
of degree $n$, which minimizes
the norm $\norm{P_n^L}_L$.

\item[{\rm{(iv)}}] The polynomial $P_n^L$ can be equivalently defined
as the monic polynomial of
degree $n$ which satisfies
\begin{equation}
\ang{P_n^L, f}_L = {\boldsymbol 0} \quad \text{ for any } f\in \pol,
\quad \deg f<n.
\label{b2.3}
\end{equation}

\item[{\rm{(v)}}] One has $P_n^L(x) = P_n^R(x)^\dagger$ for all  $x
\in \R$ and
\begin{equation}
\label{b2.4} \ang{P_n^R,P_n^R}_R = \ang{P_n^L,P_n^L}_L.
\end{equation}
\end{SL}
\end{lemma}

\begin{proof} As noted, $\pol$ has an inner product $\langle\cdot,
\cdot\rangle_R$, so
there is an orthogonal projection $\pi_n^{(R)}$ onto $\pol_n$
discussed in Section~\ref{s1.2}. Then
\begin{equation}\lb{2.6}
P_n^R(x)=x^n -\pi_{n-1}^{(R)}(x^n).
\end{equation}
As usual, in inner product spaces, this uniquely minimizes $x^n-Q$
over all $Q\in\pol_{n-1}$.
It clearly obeys
\begin{equation}\lb{2.6a}
\Tr(\ang{P_n^R,f}_R)=0
\end{equation}
for all $f\in\pol_{n-1}$. But then for any matrix $\alpha$,
\[
\Tr(\ang{P_n^R,f}_R\alpha)=\Tr(\ang{P_n^R,f\alpha}_R)=0
\]
so \eqref{b2.2} holds.

This proves (i) and (ii). (iii) and (iv) are similar. To prove (v),
note that
$P_n^L(x)=P_n^R(x)^\dagger$ follows from the criteria \eqref{b2.2},
\eqref{b2.3}.
The identity \eqref{b2.4} follows from \eqref{b2}.
\end{proof}

\begin{lemma}\label{l.b2}
Let $\mu$ be non-trivial. For any monic polynomial $P$, we have $\det
\ang{P,P}_R \not= 0$ and
$\det\ang{P,P}_L \not= 0$.
\end{lemma}

\begin{proof}
Let $P$ be a monic polynomial of degree $n$ such that $\ang{P,P}_R$
has a non-trivial
kernel. Then one can find $\alpha\in\calM_l$, $\alpha\not=
{\boldsymbol 0}$, such that
$\alpha^\dagger \ang{P,P}_R\alpha = {\boldsymbol 0}$. It follows that
$\norm{P\alpha }_R
= 0$. But since $P$ is monic, the leading coefficient of $P\alpha$ is
$\alpha$, so
$P\alpha \not= {\boldsymbol 0}$, which contradicts the non-triviality
assumption.
A similar argument works for $\ang{P,P}_L$.
\end{proof}

By the orthogonality of $Q_n-P_n^R$ to $P_n^R$ for any monic
polynomial $Q_n$ of degree $n$, we have
\begin{equation}\lb{2.7a}
\ang{Q_n,Q_n}_R = \ang{Q-P_n^R, Q-P_n^R}_R + \ang{P_n^R, P_n^R}_R
\end{equation}
and, in particular,
\begin{equation}\lb{2.7b}
\ang{P_n^R,P_n^R}_R \leq \ang{Q_n,Q_n}_R
\end{equation}
with (by non-triviality) equality if and only if $Q_n=P_n^R$.
Since $\Tr$ and $\det$ are strictly monotone on strictly positive
matrices, we have the following variational principles
(\eqref{2.7c} restates (i) of Lemma~\ref{l.b1}):

\begin{theorem}\lb{T2.3AA} For any monic $Q_n$ of degree $n$, we have
\begin{gather}
\norm{Q_n}_R \geq \norm{P_n^R}_R, \lb{2.7c} \\
\det\ang{Q_n,Q_n}_R \geq \det\ang{P_n^R,P_n^R}_R \lb{2.7d}
\end{gather}
with equality if and only if $P_n^R=Q_n$.
\end{theorem}

\subsubsection{Expansion} \lb{s2.1.3}

\begin{theorem}\lb{T2.3A} Let $d\mu$ be non-trivial.
\begin{SL}
\item[{\rm{(i)}}] We have
\begin{equation}\lb{2.6A}
\ang{P_k^R,P_n^R}_R = \gamma_n \delta_{kn}
\end{equation}
for some positive invertible matrices $\gamma_n$.

\item[{\rm{(ii)}}] $\{P_k^R\}_{k=0}^n$ are a right-module basis for $
\pol_n$; indeed, any $f\in\pol_n$ has
a unique expansion,
\begin{equation}\lb{2.6B}
f=\sum_{j=0}^n P_j^R f_j^R.
\end{equation}
Indeed, essentially by \eqref{1.33c},
\begin{equation}\lb{2.6C}
f_j^R = \gamma_j^{-1} \ang{P_j^R,f}_R.
\end{equation}
\end{SL}
\end{theorem}

\begin{remark} There are similar formulas for $\ang{\cdot,\cdot}_L$.
By \eqref{b2.4},
\begin{equation}\lb{2.6D}
\ang{P_k^L, P_n^L}_L =\gamma_n \delta_{kn}
\end{equation}
(same $\gamma_n$, which is why we use $\gamma_n$ and not $\gamma_n^R$).
\end{remark}

\begin{proof}  (i) \eqref{2.6A} for $n<k$ is immediate from \eqref
{b2.3} and for $n>k$ by symmetry.
$\gamma_n\geq 0$ follows from \eqref{1.25a}. By Lemma~\ref{l.b2}, $
\det(\gamma_n)\neq 0$, so
$\gamma_n$ is invertible.

\smallskip
(ii) Map $(\calM_l)^{n+1}$ to $\pol_n$ by
\[
\langle \alpha_0,\dots,\alpha_n\rangle \mapsto \sum_{j=0}^n P_j^R
\alpha_j
\equiv X (\alpha_0, \dots, \alpha_n).
\]
By \eqref{2.6A},
\[
\alpha_j=\gamma_j^{-1} \ang{P_j^R,X(\alpha_0,\dots,\alpha_n)}
\]
so that map is one-one. By dimension counting, it is onto.
\end{proof}

\subsubsection{Recurrence Relations for Monic Orthogonal Polynomials}
Denote by $\zeta_n^R$ (resp.\ $\zeta_n^L$) the coefficient of
$x^{n-1}$ in $P_n^R(x)$ (resp.\ $P_n^L(x)$), that is,
\begin{align*}
P_n^R(x)&=x^n \bdone +\zeta_n^R x^{n-1}+\text{lower order terms},
\\
P_n^L(x)&=x^n \bdone +\zeta_n^L x^{n-1}+\text{lower order terms}.
\end{align*}
Since $P_n^R(x)^\dagger=P_n^L(x)$, we have $(\zeta_n^R)^\dagger=
\zeta_n^L$.
Using the parameters $\gamma_n$ of \eqref{2.6A} and $\zeta_n^R$, $
\zeta_n^L$ one can write down
recurrence relations for $P_n^R(x)$, $P_n^L(x)$.

\begin{lemma}
\begin{SL}
\item[{\rm{(i)}}] We have a commutation relation
\begin{equation}
\gamma_{n-1}(\zeta_n^R-\zeta_{n-1}^R)
=
(\zeta_n^L-\zeta_{n-1}^L)\gamma_{n-1}.
\label{b.2.7.1}
\end{equation}

\item[{\rm{(ii)}}] We have the recurrence relations
\begin{align}
xP_n^R(x)&=P_{n+1}^R(x)+P_n^R(x)(\zeta_n^R-\zeta_{n+1}^R)+P_{n-1}(x)
\gamma_{n-1}^{-1}\gamma_n,
\label{b.2.7.2}
\\
xP_n^L(x)&=P_{n+1}^L(x)+(\zeta_n^L-\zeta_{n+1}^L)P_n^L(x)+\gamma_n
\gamma_{n-1}^{-1}P_{n-1}^L(x).
\label{b.2.7.3}
\end{align}
\end{SL}
\end{lemma}

\begin{proof}
(i) We have
$$
P_n^R(x)-xP_{n-1}^R(x)=(\zeta_n^R-\zeta_{n-1}^R)x^{n-1}+\text{lower
order terms}
$$
and so
\begin{align*}
(\zeta_n^L-\zeta_{n-1}^L)\gamma_{n-1}
&= (\zeta_n^R-\zeta_{n-1}^R)^\dagger \ang{P_{n-1}^R, P_{n-1}^R}_R \\
& = (\zeta_n^R-\zeta_{n-1}^R)^\dagger \ang{x^{n-1}, P_{n-1}^R}_R \\
& = \ang{P_n^R-xP_{n-1}^R,P_{n-1}^R}_R \\
& = \ang{P_n^R,P_{n-1}^R}_R - \ang{xP_{n-1}^R,P_{n-1}^R}_R \\
& = -\ang{xP_{n-1}^R,P_{n-1}^R}_R \\
& = -\ang{P_{n-1}^R,xP_{n-1}^R}_R \\
& = \ang{P_{n-1}^R,P_n^R-xP_{n-1}^R}_R \\
& = \ang{P_{n-1}^R,x^{n-1}(\zeta_n^R-\zeta_{n-1}^R)}_R \\
& = \ang{P_{n-1}^R,x^{n-1}}_R (\zeta_n^R-\zeta_{n-1}^R) \\
& = \gamma_{n-1} (\zeta_n^R-\zeta_{n-1}^R).
\end{align*}

\smallskip
(ii)
By Theorem~\ref{T2.3A},
$$
xP_n^R (x)=P_{n+1}^R(x)C_{n+1}+P_n^R(x)C_n+P_{n-1}^R(x)C_{n-1}+\cdots
+P_0^R C_0
$$
with some matrices $C_0, \dots, C_{n+1}$. It is straightforward
that $C_{n+1}={\boldsymbol 1}$ and $C_n=\zeta_n^R-\zeta_{n+1}^R$. By
the orthogonality property \eqref{b2.2}, we find
$C_0=\dots=C_{n-2} = {\boldsymbol 0}$. Finally, it is easy to
calculate $C_{n-1}$:
\begin{align*}
\gamma_n & = \ang{P_n^R,xP_{n-1}^R}_R = \ang{xP_n^R,P_{n-1}^R}_R \\
& = \ang{P_{n+1}^R+P_n^R(\zeta_n^R-\zeta_{n+1}^R)+P_{n-1}^RC_{n-1},P_
{n-1}^R}_R \\
& = C_{n-1}^\dagger \gamma_{n-1}
\end{align*}
and so, taking adjoints and using self-adjointness of $\gamma_j$, $C_
{n-1} =
\gamma_{n-1}^{-1}\gamma_n$. This proves \eqref{b.2.7.2}; the other
relation \eqref{b.2.7.3}
is obtained by conjugation.
\end{proof}

\subsubsection{Normalized Orthogonal Polynomials} We call  $p_n^R\in
\pol$ a right orthonormal
polynomial if $\deg p_n^R\leq n$ and
\begin{gather}
\ang{p_n^R,f}_R= \bdzero \text{ for every $f\in \pol$ with } \deg
f<n, \label{b2.8}
\\
\ang{p_n^R,p_n^R}_R={\boldsymbol 1}.
\label{b2.9}
\end{gather}
Similarly, we call  $p_n^L\in \pol$ a left orthonormal polynomial if $
\deg
p_n^L\leq n$ and
\begin{gather}
\ang{p_n^L,f}_L=\bdzero \text{ for every $f\in \pol$ with }\deg f<n,
\label{b2.10}
\\
\ang{p_n^L,p_n^L}_L={\boldsymbol 1}.
\label{b2.11}
\end{gather}

\begin{lemma}\label{l.b3}
Any orthonormal polynomial has the form
\begin{equation}
p_n^R(x)=P_n^R(x)\ang{P_n^R,P_n^R}_R^{-1/2}\sigma_n,
\qquad
p_n^L(x)=\tau_n\ang{P_n^L,P_n^L}_L^{-1/2}P_n^L(x)
\label{b2.11.1}
\end{equation}
where $\sigma_n,\tau_n\in\calM_l$ are unitaries. In particular,
$\deg p_n^R=\deg p_n^L=n$.
\end{lemma}

\begin{proof}
Let $K_n$ be the coefficient of $x^n$ in $p_n^R$.
Consider the polynomial $q(x)=P_n^R(x) K_n-p_n^R(x)$,
where $P_n^R$ is the monic orthogonal polynomial from Lemma~\ref{l.b1}.
Then $\deg q<n$ and so from \eqref{b2.2} and \eqref{b2.8}, it follows
that
$\ang{q,q}_R=0$ and so $q(x)$
vanishes identically. Thus, we have
\begin{equation}
{\boldsymbol 1}=\ang{p_n^R,p_n^R}_R
=
K_n^\dagger \ang{P_n^R,P_n^R}_R K_n
\label{b.2.12}
\end{equation}
and so $\det(K_n)\not=0$.
{}From \eqref{b.2.12} we get
$(K_n^\dagger)^{-1}K_n^{-1}=\ang{P_n^R,P_n^R}_R$,
and so $K_n K_n^\dagger=\ang{P_n^R,P_n^R}_R^{-1}$.
{}From here we get $K_n=\ang{P_n^R,P_n^R}_R^{-1/2}\sigma_n$
with a unitary $\sigma_n$.
The proof for $p_n^L$ is similar.
\end{proof}

By Theorem~\ref{T2.3A}, the polynomials $p_n^R$ form a right
orthonormal module basis in
$\pol_R$. Thus, for any $f\in \pol_R$, we have
\begin{equation}
f(x)=\sum_{m=0}^\infty p_m^R f_m,
\qquad f_m=\ang{p_m^R, f}_R
\label{b5}
\end{equation}
and the Parseval identity
\begin{equation}
\sum_{m=0}^\infty \Tr(f_m f_m^\dagger)=\norm{f}^2_R
\label{b6}
\end{equation}
holds true.
Obviously, since $f$ is a polynomial, there are only finitely many
non-zero
terms in \eqref{b5} and \eqref{b6}.

\subsection{Block Jacobi Matrices} \lb{s2.2}

The study of block Jacobi matrices goes back at least to Krein \cite
{Krein3}.

\subsubsection{Block Jacobi Matrices as Matrix Representations}

Suppose that a sequence of unitary matrices ${\boldsymbol 1} =
\sigma_0, \sigma_1,
\sigma_2, \ldots$ is fixed, and $p_n^R$ are defined according to
\eqref{b2.11.1}. As
noted above, $p_n^R$ form a right orthonormal basis in $\pol_R$.

The map $f(x)\mapsto xf(x)$ can be considered as a right homomorphism
in $\pol_R$.
Consider the matrix $J_{nm}$ of this homomorphism with respect to the
basis $p_n^R$, that
is,
\begin{equation}\lb{2.22a}
J_{nm}=\ang{p_{n-1}^R,xp_{m-1}^R}_R.
\end{equation}
Following Killip--Simon \cite{KS} and Simon \cite{S,S2,Rice}, our
Jacobi matrices are
indexed with $n=1,2,\dots$ but, of course, $p_n$ has $n=0,1,2,\dots$.
That is why
\eqref{2.22a} has $n-1$ and $m-1$.

As in the scalar case, using the orthogonality properties of $p_n^R$,
we get that $J_{nm}
= {\boldsymbol 0}$ if $\abs{n-m}>1$. Denote
$$
B_n=J_{nn}=\ang{p_{n-1}^R,x p_{n-1}^R}_R
$$
and
$$
A_n=J_{n,n+1}=J_{n+1,n}^\dagger=\ang{p_{n-1}^R,x p_{n}^R}_R.
$$
Then we have
\begin{equation}
J=
\begin{pmatrix}
B_1 & A_1 & {\boldsymbol 0} & \cdots \\
A_1^\dagger & B_2 & A_2 & \cdots \\
{\boldsymbol 0} & A_2^\dagger & B_3 &  \cdots \\
\vdots & \vdots & \vdots & \ddots
\end{pmatrix}
\label{b6.1}
\end{equation}

Applying \eqref{b5} to $f(x)=xp_n^R(x)$, we get the recurrence relation
\begin{equation}
xp_n^R(x)=p_{n+1}^R(x)A^\dagger_{n+1}+p_{n}^R(x)B_{n+1}+p_{n-1}^R(x)A_
{n},
\quad
n=1,2,\dots .
\label{b7}
\end{equation}
If we set $p_{-1}^R(x) = {\boldsymbol 0}$ and $A_0={\boldsymbol 1}$,
the relation
\eqref{b7} also holds for $n=0$. By \eqref{b2}, we can always pick
$p_n^L$
so that for $x$ real, $p_n^L(x) = p_n^R(x)^\dagger$, and thus for
complex $z$,
\begin{equation} \lb{2.30a}
p_n^L(z) = p_n^R(\bar z)^\dagger
\end{equation}
by analytic continuation. By conjugating \eqref{b7}, we get
\begin{equation}
xp_n^L(x)=
A_{n+1} p_{n+1}^L(x) + B_{n+1} p_{n}^L(x) +A_{n}^\dagger p_{n-1}^L(x),
\quad
n=0,1,2,\dots .
\label{b8}
\end{equation}
Comparing this with the recurrence relations \eqref{b.2.7.2}, \eqref
{b.2.7.3},
we get
\begin{equation}
A_n=\sigma_{n-1}^\dagger\gamma_{n-1}^{-1/2}\gamma_n^{1/2}\sigma_n,
\qquad
B_n=\sigma_{n-1}^\dagger\gamma_{n-1}^{1/2}(\zeta_{n-1}^R-\zeta_n^R)
\gamma_{n-1}^{-1/2}\sigma_{n-1}.
\label{b8.00}
\end{equation}
In particular, $\det A_n\not=0$ for all $n$.

Notice that since $\sigma_n$ is unitary, $\abs{\det(\sigma_n)}=1$, so
\eqref{b8.00} implies $\det(\gamma_n^{1/2}) = \det(\gamma_{n-1}^{1/2})
\abs{\det(A_n)}$ which, by induction, implies that
\begin{equation} \lb{2.32a}
\det \ang{P_n^R,P_n^R} = \abs{\det(A_1\cdots A_n)}^2.
\end{equation}

Any block matrix of the form \eqref{b6.1} with $B_n = B_n^\dagger$
and $\det A_n \not= 0$
for all $n$ will be called a block Jacobi matrix corresponding to the
Jacobi parameters
$A_n$ and $B_n$.

\subsubsection{Basic Properties of Block Jacobi Matrices}

Suppose we are given a block Jacobi matrix $J$ corresponding to
Jacobi parameters $A_n$
and $B_n$, where $B_n = B_n^\dagger$ and $\det A_n \not= 0$ for each
$n$.

Consider the Hilbert space $\h_v = \ell^2(\Z_+, \C^l)$ (here
$\bbZ_+ = \{1,2,3,\dots\}$) with inner product
$$
\langle f,g \rangle_{\h_v} = \sum_{n=1}^\infty \langle f_n , g_n
\rangle_{\C^l}
$$
and orthonormal basis $\{ e_{k,j} \}_{k \in \Z_+, 1 \le j \le l}$, where
$$
(e_{k,j})_n = \delta_{k,n} v_j
$$
and $\{ v_j \}_{1 \le j \le l}$ is the standard basis of $\C^l$. $J$
acts on $\h_v$ via
\begin{equation}\label{jacobiop}
(J f)_n = A_{n-1}^\dagger f_{n-1} + B_n f_n + A_n f_{n+1}, \quad f\in
\h_v
\end{equation}
(with $f_0 = 0$) and defines a symmetric operator on this space. Note
that using
invertibility of the $A_n$'s, induction shows
\begin{equation}\label{cyclicity}
\Span \{ e_{k,j}\, \colon 1 \le k \le n, \, 1 \le j \le l \} = \Span \{
J^{k-1} e_{1,j}\, \colon 1 \le k \le n, \, 1 \le j \le l \}
\end{equation}
for every $n \ge 1$. We want to emphasize that elements of $\calH_v$
and $\calH$ are
vector-valued and matrix-valued, respectively. For this reason, we
will be interested
in both matrix- and vector-valued solutions of the basic difference
equations.

We will consider only bounded block Jacobi matrices, that is, those
corresponding to
Jacobi parameters satisfying
\begin{equation}\label{boundedbjm}
\sup_{n} \, \Tr (A_n^\dagger A_n+B_n^\dagger B_n) < \infty.
\end{equation}
Equivalently,
\begin{equation} \lb{2.29a}
\sup_n (\norm{A_n} + \norm{B_n}) <\infty.
\end{equation}
In this case, $J$ is a bounded self-adjoint operator. This is
equivalent to $\mu$ having
compact support.

We call two Jacobi matrices $J$ and $\tilde J$ equivalent if there
exists a sequence of
unitaries $u_n\in\calM_l$, $n\geq 1$, with $u_1={\boldsymbol 1}$ such
that $\tilde
J_{nm}=u_n^\dagger J_{nm} u_m$.
{}From Lemma~\ref{l.b3} it is clear that if $p_n^R$, $\tilde p_n^R$
are two sequences of normalized orthogonal polynomials, corresponding
to the same
measure (but having different normalization), then the Jacobi matrices
$J_{nm}=\ang{p_{n-1}^R, xp_{m-1}^R}_R$ and $\tilde J_{nm}=\ang{\tilde
p_{n-1}^R,
x\tilde p_{m-1}^R}_R$ are equivalent ($u_n=\sigma_{n-1}^\dagger \ti
\sigma_{n-1}$).
Thus,
\begin{equation} \lb{2.35a}
\ti B_n = u_n^\dagger B_n u_n, \qquad \ti A_n = u_n^\dagger A_n u_{n+1}.
\end{equation}

Therefore, we have a map
\begin{equation}
\begin{split}
\Phi \,\colon  \mu\mapsto &\{ J \,\colon J_{mn}=\ang{p_{n-1}^R,x p_
{m-1}^R}_R, \, p_n^R \\
&\quad \text{ correspond to $d\mu$ for some normalization}\} \label{Phi}
\end{split}
\end{equation}
from the set of all Hermitian positive semi-definite non-trivial
compactly supported measures to the set of all equivalence classes
of bounded block Jacobi matrices. Below, we will see how to invert
this map.

\subsubsection{Special Representatives of the Equivalence Classes}

Let $J$ be a block Jacobi matrix with the Jacobi parameters $A_n$,
$B_n$. We say that $J$
is:
\begin{itemize}
\item
of type~1, if $A_n>\bdzero$ for all $n$;
\item
of type~2, if $A_1 A_2 \cdots A_n> \bdzero$ for all $n$;
\item
of type~3, if $A_n\in {\mathcal L}$ for all $n$.
\end{itemize}
Here, $\mathcal L$ is the class of all lower triangular matrices
with strictly positive elements on the diagonal. Type~3 is of
interest because they correspond precisely to bounded Hermitian
matrices with $2l+1$ non-vanishing diagonals with the extreme
diagonals strictly positive; see Section~\ref{s1.4}(c). Type~2 is
the case where the leading coefficients of $p_n^R$ are strictly
positive definite.

\begin{theorem} \label{l.b4}
\begin{SL}
\item[{\rm{(i)}}] Each equivalence class of block Jacobi matrices
contains exactly one element
each of type~1, type~2, or type~3.

\item[{\rm{(ii)}}] Let $J$ be a block Jacobi matrix corresponding to
a sequence of polynomials
$p_n^R$ as in \eqref{b2.11.1}. Then $J$ is of type~2 if and only if $
\sigma_n =
{\boldsymbol 1}$ for all $n$.
\end{SL}
\end{theorem}

\begin{proof}
The proof is based on the following two well-known facts:
\begin{SL}
\item[(a)] For any $t\in\calM_l$ with $\det(t)\not=0$, there
exists a unique unitary $u\in \calM_l$ such that $tu$ is Hermitian
positive semi-definite: $tu\geq 0$. \item[(b)] For any
$t\in\calM_l$ with $\det(t)\not=0$, there exists a unique unitary
$u\in\calM_l$ such that $tu\in \calL$.
\end{SL}

We first prove that every equivalence class of block Jacobi matrices
contains at least
one element of type~1.
For a given sequence $A_n$, let us construct a sequence $u_1=
\bdone,u_2,u_3,\dots$
of unitaries such that $u_n^\dagger A_n u_{n+1}\geq\bdzero$.
By the existence part of (a), we find $u_2$ such that $A_1 u_2\geq
\bdzero$,
then find $u_3$ such that $u_2^\dagger A_2 u_3\geq\bdzero$, etc.
This, together with \eqref{2.35a}, proves the statement.
In order to  prove the uniqueness part, suppose we have $A_n\geq\bdzero$
and $u_n^\dagger A_n u_{n+1}\geq\bdzero$ for all $n$.
Then, by the uniqueness part of (a), $A_1\geq\bdzero$ and $A_1 u_2
\geq \bdzero$
imply $u_2=\bdone$; next, $A_2\geq\bdzero$ and $u_2^\dagger A_2
u_3=A_2 u_3\geq\bdzero$
imply $u_3=\bdone$, etc.

The statement (i) concerning type~$3$ can be proven in the same way,
using (b) instead of
(a).

The statement (i) concerning type~$2$ can be proven similarly.
Existence: find $u_2$
such that $A_1 u_2\geq \bdzero$, then $u_3$ such that
$(A_1 u_2)(u_2^\dagger A_2 u_3)=A_1A_2u_3\geq \bdzero$, etc.
Uniqueness: if
$A_1\dots A_n\geq\bdzero$  and $A_1\cdots A_n u_{n+1}\geq\bdzero$,
then $u_{n+1}=\bdone$.

By \eqref{b8.00}, we have $A_1 A_2 \cdots A_n = \gamma_n^{1/2}\sigma_n
$ and the statement
(ii) follows from the positivity of $\gamma_n$.
\end{proof}

We say that a block Jacobi matrix $J$ belongs to the Nevai class if
$$
B_n\to \bdzero \text{ and } A_n^\dagger A_n\to {\boldsymbol 1} \text
{ as }n\to\infty.
$$
It is clear that $J$ is in the Nevai class if and only if all
equivalent Jacobi matrices
belong to the Nevai class.

\begin{theorem}\label{l.b5}
If $J$ belongs to the Nevai class and is of type~1 or type~3, then
$A_n\to 1$ as $n\to
\infty$.
\end{theorem}

\begin{proof}
If $J$ is of type~$1$, then $A_n^\dagger A_n=A_n^2\to {\boldsymbol 1}
$ clearly
implies $A_n\to {\boldsymbol 1}$ since square root is continuous on
positive Hermitian
matrices.

Suppose $J$ is of type~$3$. We shall prove that $A_n\to \bdone$ by
considering the rows
of the matrix $A_n$ one by one, starting from the $l$th row. Denote $
(A_n)_{jk} =
a_{j,k}^{(n)}$. We have
$$
(A_n^\dagger A_n)_{ll}=(a_{l,l}^{(n)})^2\to1,
\text{ and so } a_{l,l}^{(n)} \to1.
$$
Then, for any $k<l$, we have
$$
(A_n^\dagger A_n)_{lk}=a^{(n)}_{l,l}a^{(n)}_{l,k}\to 0,
\text{ and so } a_{l,k}^{(n)} \to 0.
$$
Next, consider the $(l-1)$st row. We have
$$
(A_n^\dagger A_n)_{l-1,l-1}=(a_{l-1,l-1}^{(n)})^2+\abs{a_{l,l-1}^
{(n)}}^2 \to 1
$$
and so, using the previous step, $a_{l-1,l-1}^{(n)}\to1$ as $n\to
\infty$.
Then for all $k<l-1$, we have
$$
(A_n^\dagger A_n)_{l-1,k}=\overline{a_{l-1,l-1}^{(n)}}\, a^{(n)}_{l-1,k}
+\overline{a_{l,l-1}^{(n)}} \, a^{(n)}_{l,k}\to 0
$$
and so, using the previous steps, $a_{l-1,k}\to 0$.
Continuing this way, we get $a^{(n)}_{j,k}\to \delta_{j,k}$ as required.
\end{proof}

It is an interesting open question if this result also applies to the
type~2 case.

\subsubsection{Favard's Theorem}

Here we construct an inverse of the mapping $\Phi$ (defined by \eqref
{Phi}).
Thus, $\Phi$ sets up a bijection between non-trivial measures of compact
support and equivalence classes of bounded block Jacobi matrices.

Before proceeding to do this, let us prove:

\begin{lemma}\label{l.b6}
The mapping $\Phi$ is injective.
\end{lemma}

\begin{proof}
Let $\mu$ and $\tilde \mu$ be two Hermitian positive semi-definite
non-trivial compactly supported measures. Suppose that
$\Phi(\mu)=\Phi(\tilde \mu)$.

Let $p_n^R$ and $\tilde p_n^R$ be normalized orthogonal polynomials
corresponding to $\mu$ and $\tilde\mu$. Suppose that the normalization
both for $p_n^R$ and for $\tilde p_n^R$ has been chosen such that
$\sigma_n={\boldsymbol 1}$ (see \eqref{b2.11.1}), that is, type~2.
{}From Lemma~\ref{l.b4} and the assumption $\Phi(\mu)=\Phi(\tilde \mu)
$ it
follows that the corresponding Jacobi matrices coincide, that is,
$\ang{p_n^R,xp_m^R}_R=\ang{\tilde p_n^R,x\tilde p_m^R}_R$
for all $n$ and $m$.
Together with the recurrence relation \eqref{b7} this yields $p_n^R=
\tilde p_n^R$
for all $n$.

For any $n \geq 0$, we can represent $x^n$ as
$$
x^n = \sum_{k=0}^{n} p_k^R(x) C^{(n)}_k =\sum_{k=0}^{n} \tilde p^R_k
(x) \tilde C^{(n)}_k.
$$
The coefficients $C^{(n)}_k$ and $\tilde C^{(n)}_k$ are completely
determined by the
coefficients of the polynomials $p_n^R$ and $\tilde p_n^R$ and so
$C^{(n)}_k=\tilde C^{(n)}_k$ for all $n$ and $k$.

For the moments of the measure $\mu$, we have
$$
\int x^n d\mu(x)= \ang{{\boldsymbol 1}, x^n}_R
=
\sum_{k=0}^n \ang{ {\boldsymbol 1}, p^R_k C_k^{(n)}}_R
=
\ang{ {\boldsymbol 1}, {\boldsymbol 1}}_R \, C_0^{(n)} =C_0^{(n)}.
$$
Since the same calculation is valid for the measure $\tilde \mu$, we get
$$
\int x^n d\mu(x)=\int x^n d\tilde \mu(x)
$$
for all $n$. It follows that
$$
\int f(x) d\mu(x)g(x)=\int f(x) d\tilde \mu(x) g(x)
$$
for all matrix-valued polynomials $f$ and $g$, and so the measures $
\mu$ and $\tilde \mu$
coincide.
\end{proof}

We can now construct the inverse of the map $\Phi$. Let a block
Jacobi matrix $J$ be given. By a version of the spectral theorem
for self-adjoint operators with finite multiplicity (see, e.g.,
\cite[Sect.~72]{AG}), there exists a matrix-valued measure $d\mu$
with
\begin{equation}
\langle e_{1,j} , f(J) e_{1,k} \rangle_{\h_v} = \int f(x) \, d\mu_
{j,k}(x)
\label{defmu}
\end{equation}
and an isometry
$$
R \, \colon \h_v \to L^2(\R,d\mu;\C^l)
$$
such that (recall that $\{v_j\}$ is the standard basis in $\C^l$)
\begin{equation}\label{constant}
[R e_{1,j}](x) = v_j, \quad 1 \le j \le l,
\end{equation}
and, for any $g\in \h_v$, we have
\begin{equation}\label{diagonal}
(RJg)(x) = x (Rg)(x).
\end{equation}
If the Jacobi matrices $J$ and $\tilde J$ are equivalent, then we
have $\tilde J = U^*\! J
U$ for some $U = \oplus_{n=1}^\infty u_n$, $u_1={\boldsymbol 1}$. Thus,
$$
\langle e_{1,j} , f(\tilde J) e_{1,k} \rangle_{\h_v} = \langle U e_
{1,j} , f(J) Ue_{1,k}
\rangle_{\h_v} = \langle  e_{1,j} , f(J) e_{1,k} \rangle_{\h_v}
$$
and so the measures corresponding to $J$ and $\tilde J$ coincide.
Thus, we have a map
\begin{equation}
\Psi \colon \{\tilde J \colon \tilde J \text{ is equivalent to } J\}
\mapsto \mu
\label{Psi}
\end{equation}
from the set of all equivalence classes of bounded block Jacobi
matrices to the set of all Hermitian positive semi-definite
compactly supported measures.

\begin{theorem}\label{favard}
\begin{SL}
\item[{\rm{(i)}}] All measures in the image of the map $\Psi$ are non-
degenerate.
\item[{\rm{(ii)}}] $\Phi\circ \Psi=\text{id}$.
\item[{\rm{(iii)}}] $\Psi\circ \Phi=\text{id}$.
\end{SL}
\end{theorem}

\begin{proof}
(i) To put things in the right context, we first recall that
$\norm{\cdot}_{\h_v}$ is a norm (rather than a semi-norm), whereas
$\snorm{\cdot}$  on $\mathcal{V}$ (cf.~\eqref{f.vsn}) is, in
general, a semi-norm. Using the assumption that $\det(A_k)\not=0$
for all $k$ (which is included in our definition of a Jacobi
matrix), we will prove that $\snorm{\cdot}$ is in fact a norm.
More precisely, we will prove that  $\snorm{p}>0$ for any
non-zero polynomial $p\in\mathcal{V}$; by Lemma~\ref{nondeg} this will
imply that $\mu$ is non-degenerate.

Let $p\in\mathcal{V}$ be a non-zero polynomial, $\deg p=n$. Notice
that \eqref{constant} and \eqref{diagonal} give
\begin{equation}\label{conjugation}
[R J^k e_{1,j}](x) = x^k v_j
\end{equation}
for every $k \ge 0$ and $1 \le j \le l$. This shows that $p$ can
be represented as $p=Rg$, where $g=\sum_{k=0}^n J^k f_k$, and
$f_0, \dots, f_n$ are vectors in $\calH_v$ such that $\langle f_i,
e_{j,k}\rangle_{\calH_v}=0$ for all $i=0,\dots,n$, $j\geq2$,
$k=1,\dots, l$ (i.e., the only non-zero components of $f_j$ are in
the first $\bbC^l$ in $\calH_v$). Assumption $\deg p=n$ means
$f_n\not=0$.

Since $R$ is isometric, we have $\snorm{p}=\norm{g}_{\h_v}$, and
so we have to prove that $g\not=0$.  Indeed, suppose that $g=0$.
Using the assumption $\det(A_k)\not=0$ and the tri-diagonal nature
of $J$, we  see that $\sum_{k=0}^n J^k f_k=0$ yields $f_n=0$,
contrary to our assumption.

\smallskip
(ii) Consider the elements $R e_{n,k}\in L^2(\R,d\mu;\C^l)$. First
note that, by
\eqref{cyclicity} and \eqref{conjugation}, $R e_{n,k}$ is a
polynomial of degree at most
$n-1$. Next, by the unitarity of $R$, we have
\begin{equation}
\langle R e_{n,k}, R e_{m,j}\rangle_{L^2(\R,d\mu;\C^l)} = \delta_
{m,n} \delta_{k,j}.
\label{ortho}
\end{equation}
Let us construct matrix-valued polynomials $q_n(x)$, using $R e_{n,
1}, R e_{n,2}, \dots,
R e_{n,l}$ as columns of $q_{n-1}(x)$:
$$
[q_{n-1}(x)]_{j,k} =  [R e_{n,k}(x)]_j.
$$
We have $\deg q_n \leq n$ and $\ang{q_m,q_n}_R = \delta_{m,n}
{\boldsymbol 1}$; the last
relation is just a reformulation of \eqref{ortho}. Hence the $q_n$'s
are right
normalized orthogonal polynomials with respect to the measure $d\mu$.
We find
\begin{align*}
J_{nm} & = [ \langle e_{n,j}, J e_{m,k} \rangle_{\h_v} ]_{1 \le j,k
\le l} \\
& = [ \langle R e_{n,j}, R J e_{m,k} \rangle_{L^2(\R,d\mu;\C^l)} ]_{1
\le j,k \le l} \\
& = [ \langle R e_{n,j}, x R e_{m,k} \rangle_{L^2(\R,d\mu;\C^l)} ]_{1
\le j,k \le l} \\
& =  [ \langle [q_{n-1}(x)]_{\bddot,j}, x
[q_{m-1}(x)]_{\bddot,k} \rangle_{L^2(\R,d\mu;\C^l)}]_{1 \le j,k \le
l} \\
& =  \ang{q_{n-1} , x q_{m-1}}_R
\end{align*}
as required.

\smallskip
(iii)
Follows from (ii) and from Lemma~\ref{l.b6}.
\end{proof}

\subsection{The $m$-Function}

\subsubsection{The Definition of the $m$-Function}

We denote the Borel transform of $d\mu$ by $m$:
\begin{equation}
m(z) = \int \frac{d\mu(x)}{x-z}\,  ,\qquad \Ima z > 0. \label{b.8.0}
\end{equation}
It is a matrix-valued Herglotz function, that is, it is analytic
and obeys $\Ima m(z) > 0$. For information on matrix-valued
Herglotz functions, see \cite {GT} and references therein.
Extensions to operator-valued Herglotz functions can be found in
\cite{GKMT}.

\begin{lemma}
Suppose $d\mu$ is given, $p_n^R$ are right normalized orthogonal
polynomials, and $J$ is
the associated block Jacobi matrix. Then,
\begin{equation}\label{mfitop}
m(z) = \ang{p_0^R,(x-z)^{-1}p_0^R}_R
\end{equation}
and
\begin{equation}\label{mfitoj}
m(z) = \langle e_{1,\bddot}, (J-z)^{-1} e_{1,\bddot} \rangle_{\h_v}.
\end{equation}
\end{lemma}

\begin{proof}
Since $p_0^R = {\boldsymbol 1}$, \eqref{mfitop} is just a way of
rewriting the definition
of $m$. The second identity, \eqref{mfitoj}, is a consequence of
\eqref{defmu} and
Theorem~\ref{favard}(iii).
\end{proof}

\subsubsection{Coefficient Stripping}

If $J$ is a block Jacobi matrix corresponding to invertible $A_n$'s
and Hermitian
$B_n$'s, we denote the $k$-times stripped block Jacobi matrix,
corresponding to $\{
A_{k+n} , B_{k+n} \}_{n \ge 1}$, by $J^{(k)}$. That is,
$$
J^{(k)} =
\begin{pmatrix}
B_{k+1} & A_{k+1} & {\boldsymbol 0} & \cdots \\
A_{k+1}^\dagger & B_{k+2} & A_{k+2} & \cdots \\
{\boldsymbol 0} & A_{k+2}^\dagger & B_{k+3} &  \cdots \\
\vdots & \vdots & \vdots & \ddots
\end{pmatrix}.
$$
The $m$-function corresponding to $J^{(k)}$ will be denoted by $m^
{(k)}$. Note that, in
particular, $J^{(0)} = J$ and $m^{(0)} = m$.

\begin{proposition}
Let $J$ be a block Jacobi matrix with $\sigma_\mathrm{ess}(J)
\subseteq [a,b]$. Then, for
every $\varepsilon > 0$, there is $k_0 \ge 0$ such that for $k \ge k_0
$, we have that
$\sigma(J^{(k)}) \subseteq [a-\varepsilon,b+\varepsilon]$.
\end{proposition}

\begin{proof}
This is an immediate consequence of (the proof of) \cite[Lemma~1]{Den}.
\end{proof}

\begin{proposition}[Due to Aptekarev--Nikishin \cite{AN}]\label{strip}
We have that
$$
m^{(k)}(z)^{-1}  = B_{k+1} - z - A_{k+1} m^{(k+1)}(z) A_{k+1}^\dagger
$$
for $\Ima z > 0$ and $k \ge 0$.
\end{proposition}

\begin{proof} It suffices to handle the case $k=0$.
Given \eqref{mfitoj}, this is a special case of a general formula for
$2\times 2$ block
operator matrices, due to Schur \cite{Schur}, that reads
$$
\begin{pmatrix}
A & B \\
C & D
\end{pmatrix}^{-1} = \begin{pmatrix}
(A-BD^{-1}C)^{-1} & -A^{-1} B (D-CA^{-1}B)^{-1} \\
-D^{-1} C (A-BD^{-1}C)^{-1} & (D-CA^{-1}B)^{-1}
\end{pmatrix}
$$
and which is readily verified. Here $A=B_1-z$, $B=A_1$, $C=A_1^\dagger
$, and $D=J^{(1)}-z$.
\end{proof}

\subsection{Second Kind Polynomials} Define the second kind
polynomials by $q_{-1}^R(z) =
-{\boldsymbol 1}$,
$$
q_n^R(z)=\int_\R d\mu(x) \frac{p_n^R(z)-p_n^R(x)}{z-x},
\quad n=0,1,2,\dots .
$$
As is easy to see, for $n\geq 1$, $q_n^R$ is a polynomial of degree
$n-1$. For future reference,
let us display the first several polynomials $p_n^R$ and $q_n^R$:
\begin{alignat}{3}
p_{-1}^R(x) &= {\boldsymbol 0}, \quad & p_{0}^R(x) &= {\boldsymbol
1}, \quad & p_{1}^R(x) &=
(x-B_1)A_1^{-1}, \label{b10} \\
q_{-1}^R(x) &= -{\boldsymbol 1}, \quad & q_{0}^R(x) &= {\boldsymbol
0}, \quad & q_{1}^R(x) &=
A_1^{-1}. \label{b11}
\end{alignat}
The polynomials $q_n^R$ satisfy the equation (same form as \eqref{b7})
\begin{equation}
xq_n^R(x)=q_{n+1}^R(x)A_{n+1}^\dagger+q_{n}^R(x)B_{n+1}+q_{n-1}^R(x)A_
{n},
\quad
n=0,1,2,\dots .
\label{b12}
\end{equation}
For $n=0$, this can be checked by a direct substitution of
\eqref{b11}. For $n\geq1$, as in the scalar case, this can be
checked by taking \eqref{b7} for $x$ and for $z$, subtracting,
dividing by $x-z$, integrating over $d\mu$, and taking into
account the orthogonality relation
$$
\int d\mu(x) p_n^R(x) = {\boldsymbol 0}, \quad n\geq1.
$$
Finally, let us define
$$
\psi_n^R(z)=q_n^R(z)+m(z)p_n^R(z).
$$
According to the definition of $q_n^R$, we have
$$
\psi_n^R(z)=\ang{f_z, p_n^R}_R,
\quad
f_z(x)=(x-\bar{z})^{-1}.
$$
By the Parseval identity, this shows that for all $\Ima z>0$, the
sequence $\psi_n^R(z)$ is in $\ell^2$, that is,
\begin{equation}\lb{2.49ax}
\sum_{n=0}^\infty \Tr (\psi_n^R(z)^\dagger \psi_n^R(z)) <\infty .
\end{equation}

In the same way, we define $q_{-1}^L(z)=-{\boldsymbol 1}$,
$$
q_n^L(z)=\int_\R  \frac{p_n^L(z)-p_n^L(x)}{z-x}d\mu(x)\, ,
\quad n=0,1,2,\dots
$$
and $\psi_n^L(z)=q_n^L(z)+p_n^L(z)m(z)$.

\subsection{Solutions to the Difference Equations} For $\Ima z>0$,
consider the solutions
to the equations
\begin{align}
z u_n(z)&=
\sum_{m=1}^\infty u_m(z) J_{mn},
\quad
n=2,3,\dots,
\label{b15}
\\
z v_n(z)&=
\sum_{m=1}^\infty  J_{nm} v_m(z) ,
\quad
n=2,3,\dots .
\label{b16}
\end{align}
Clearly, $u_n(z)$ solves \eqref{b15} if and only if $v_n(z)=(u_n(\bar
z))^\dagger$
solves \eqref{b16}. In the above, we normally intend $z$ as a fixed
parameter
but it then can be used as a variable. That is, $z$ is fixed and $u_n
(z)$ is a fixed sequence, not a $z$-dependent function. A statement
like $v_n
(z)=
(u_n(\bar z))^\dagger$ means if $u_n$ is a sequence solving \eqref
{b15} for
$z=\bar z_0$, then $v_n$ obeys \eqref{b16} for $z=z_0$. Of course, if
$u_n(z)$
is a function of $z$ in a region, we can apply our estimates to all $z
$ in the
region. For any solution $\{u_n(z)\}_{n=1}^\infty$ of \eqref{b15},
let us define
\begin{equation}\lb{2.49a}
u_0(z)=zu_1(z)-u_1(z) B_1-u_2(z)A_1^\dagger.
\end{equation}
With this definition, the equation \eqref{b15} for $n=1$ is
equivalent to
$u_0(z)=0$. In the same way, we set
$$
v_0(z)=zv_1(z)-B_1 v_1(z) -A_1 v_2(z).
$$

\begin{lemma}
Let $\Ima z>0$ and suppose $\{u_n(z)\}_{n=0}^\infty$ solves
\eqref{b15} {\rm{(}}for $n\geq 2${\rm{)}} and \eqref{2.49a} and
belongs to $\ell^2$.
Then
\begin{equation}
(\Ima z) \sum_{n=1}^\infty \Tr (u_n(z)^\dagger u_n(z))
=
-\Ima \Tr (u_1(z) u_0(z)^\dagger).
\label{b17}
\end{equation}
In particular, $u_n(z) = \alpha p_{n-1}^R(z)$ is in $\ell^2$ only if
$\alpha=\bdzero$.
\end{lemma}

\begin{proof}
Denote $s_n=\Tr (u_n(z) A_{n-1}^\dagger u_{n-1}(z)^\dagger)$. Here
$A_0=\bdone$. Multiplying \eqref{b15} for $n\geq 2$ and
\eqref{2.49a} for $n=1$ by $u_n(z)^\dagger$ on the right, taking
traces, and summing over $n$, we get
$$
z\sum_{n=1}^N \Tr(u_n(z) u_n(z)^\dagger)
=
\sum_{n=1}^N s_{n+1}
+ \sum_{n=1}^N \Tr(u_n(z) B_{n+1} u_n(z)^\dagger)
+ \sum_{n=1}^N \overline{s_n}\, .
$$
Taking imaginary parts and letting $N\to \infty$, we obtain
\eqref{b17} since the middle sum is real and the outer sums cancel
up to boundary terms. Applying \eqref{b17} to $u_n(z)=\alpha p_{n-1}^R
(z)$,
we get zero in the right-hand side:
$$
(\Ima z)   \sum_{n=1}^\infty \Tr (\alpha p_{n-1}^R(z) p_{n-1}^R(z)^
\dagger\alpha^\dagger)
=0
$$
and hence $\alpha=\bdzero$ since $p_{0}^R=\bdone$.
\end{proof}

\begin{theorem}\label{t.b1}
Let $\Ima z>0$.
\begin{SL}
\item[{\rm{(i)}}] Any  solution $\{u_n(z)\}_{n=0}^\infty$ of \eqref
{b15} {\rm{(}}for $n\geq 2${\rm{)}}
can be represented as
\begin{equation}
u_n(z)=a p_{n-1}^R(z) + b q_{n-1}^R(z)
\label{b18}
\end{equation}
for suitable $a,b\in\calM_l$. In fact, $a=u_1(z)$ and $b=-u_0(z)$.

\item[{\rm{(ii)}}] A sequence  \eqref{b18} satisfies \eqref{b15} for
all $n\geq 1$ if and only if $b=0$.

\item[{\rm{(iii)}}] A sequence \eqref{b18} belongs to $\ell^2$ if and
only if $u_n(z)=c\psi_{n-1}^R(z)$ for some
$c\in\calM_l$. Equivalently, a sequence \eqref{b18} belongs to $\ell^2
$ if and only if $u_1(z)+u_0(z)m(z)=0$.
\end{SL}
\end{theorem}

\begin{proof}
(i) Let $u_n(z)$ be a solution to \eqref{b15}. Consider
$$
\tilde u_n(z)
=
u_n(z) - u_1(z) p_{n-1}^R(z) + u_0(z)q_{n-1}^R(z).
$$
Then $\tilde u_n(z)$ also solves \eqref{b15} and
$\tilde u_0(z)=\tilde u_1(z)=0$.
It follows that $\tilde u_n(z)=0$ for all $n$. This proves (i).

\smallskip
(ii) A direct substitution of \eqref{b18} into \eqref{b15} for $n=1$
yields the statement.

\smallskip
(iii)
We already know that $c\psi_{n-1}^R$ is an $\ell^2$ solution.
Suppose that $u_n(z)$ is an $\ell^2$ solution to \eqref{b15}. Rewrite
\eqref{b18} as
$$
u_n(z)=(a-bm(z))p_{n-1}^R(z)+b \psi_{n-1}^R(z).
$$
Since $\psi_n^R$ is in $\ell^2$ and $cp_n^R$ is not in $\ell^2$,
we get $a=bm(z)$, which is equivalent to $u_1(z)+u_0(z)m(z)=0$ or
to $u_n(z)=b\psi_{n-1}^R(z)$.
\end{proof}

By conjugation, we obtain:

\begin{theorem}\label{t.b2}
Let $\Ima z>0$.
\begin{SL}
\item[{\rm{(i)}}] Any  solution $\{v_n(z)\}_{n=0}^\infty$ of \eqref
{b16} {\rm{(}}for $n\geq 2${\rm{)}}
can be represented as
\begin{equation}
v_n(z)= p_{n-1}^L(z)a + q_{n-1}^L(z)b .
\label{b19}
\end{equation}
In fact, $a=v_1(z)$ and $b=-v_0(z)$.

\item[{\rm{(ii)}}] A sequence  \eqref{b19} satisfies \eqref{b16} for
all $n\geq1$ if and only if $b=0$.

\item[{\rm{(iii)}}] A sequence \eqref{b19} belongs to $\ell^2$ if
and only if $v_n(z)=\psi_{n-1}^L(z)c$ for some $c\in\calM_l$.
Equivalently, a sequence \eqref{b19} belongs to $\ell^2$ if and
only if $v_1(z)+m(z) v_0(z)=0$.
\end{SL}
\end{theorem}

\subsection{Wronskians and the Christoffel--Darboux Formula}

For any two $\calM_l$-valued sequences $u_n$, $v_n$, define the
Wronskian by
\begin{equation}
W_n(u,v)=u_{n} A_n v_{n+1} - u_{n+1} A_n^\dagger v_{n}. \label{b8.1}
\end{equation}

Note that $W_{n}(u, v)= - W_{n}(v^\dagger, u^\dagger)^\dagger$. If
$u_n(z)$ and $v_n(z)$
are solutions to \eqref{b15} and \eqref{b16}, then by a direct
calculation, we see that
$W_n(u(z),v(z))$ is independent of $n$. Put differently, if both $u_n
(z)$ and $v_n(z)$
are solutions to \eqref{b15},
then  $W_n(u(z),v(\bar z)^\dagger)$ is independent of $n$.
Or, if both $u_n(z)$ and $v_n(z)$ are solutions to \eqref{b16},
then  $W_n(u(\bar z)^\dagger, v(z))$ is independent of $n$.
In particular, by a direct evaluation for $n=0$, we get
\begin{align*}
W_n(p_{\bddot\, -1}^R(z), p_{\bddot\, -1}^R(\bar{z})^\dagger)&= W_n(q_
{\bddot\, -1}^R(z),
q_{\bddot\, -1}^R(\bar{z})^\dagger)= \bdzero, \\
W_n(p_{\bddot\, -1}^L(\bar{z})^\dagger, p_{\bddot \, -1}^L(z))&= W_n
(q_{\bddot\, -1}^L(\bar{z})^\dagger,
q_{\bddot\, -1}^L(z))= \bdzero, \\
W_n(p_{\bddot\, -1}^R(z), q_{\bddot\, -1}^R(\bar{z})^\dagger)&= W_n(p_
{\bddot\, -1}^L(\bar{z})^\dagger,
q_{\bddot\, -1}^L(z))={\boldsymbol 1}.
\end{align*}
Let both $u(z)$ and $v(z)$ be solutions to \eqref{b15} of the type
\eqref{b18}, namely,
$$
u_n(z)=a p_{n-1}^R(z) + b q_{n-1}^R(z), \quad v_n(z)=c p_{n-1}^R(z) +
d q_{n-1}^R(z).
$$
Then the above calculation implies
$$
W_n(u(z), v(\bar{z})^\dagger)=ad^\dagger - b c^\dagger.
$$

\begin{theorem}[CD Formula]\lb{T2.18}
For any $x,y\in\C$ and $n\geq1$, one has
\begin{equation}
(x-y)\sum_{k=0}^n p_k^R(x) p_k^L(y) = -W_{n+1}(p_{\bddot\, -1}^R(x),
p_{\bddot\, -1}^L(y)). \label{b8.2}
\end{equation}
\end{theorem}
\begin{proof}
Multiplying \eqref{b7} by $p_n^L(y)$ on the right and \eqref{b8}
(with $y$ in place of $x$) by $p_n^R(x)$ on the left and
subtracting, we get
$$
(x-y) p_n^R(x) p_n^R(y) = W_n(p_{\bddot\, -1}^R(x), p_{\bddot\, -1}^L
(y))
-W_{n+1}(p_{\bddot\, -1}^R(x), p_{\bddot\, -1}^L(y)).
$$
Summing over $n$ and noting that $W_0(p^R(x), p^L(y))=0$, we get the
required statement.
\end{proof}

\subsection{The CD Kernel}\lb{2.7}

The CD kernel is defined for $z,w\in\bbC$ by
\begin{align}
K_n(z,w) &= \sum_{k=0}^n p_k^R(z) p_k^R(\bar w)^\dagger \lb{2.54a} \\
&= \sum_{k=0}^n p_k^L (\bar z)^\dagger p_k^L(w). \lb{2.54b}
\end{align}

\eqref{2.54b} follows from \eqref{2.54a} and \eqref{2.30a}.
Notice that $K$ is independent of the choices $\sigma_n,\tau_n$ in
\eqref{b2.11.1} and that \eqref{b8.2} can be written
\begin{equation} \lb{2.54c}
(z-\bar w)K_n(z,w) = -W_{n+1} (p_{\bddot\, -1}^R(z), p_{\bddot\, -1}^R
(\bar w)^\dagger).
\end{equation}

The independence of $K_n$ of $\sigma,\tau$ can be understood by
noting that if $f_m$ is given by \eqref{b5}, then
\begin{equation} \lb{2.54cA}
\int K_n(z,w)\, d\mu(w) f(w) =\sum_{m=0}^n p_m^R (z) f_m
\end{equation}
so $K$ is the kernel of the orthogonal projection onto polynomials
of degree up to $n$, and so intrinsic. Similarly, if
$f_m^{(L)}=\ang{f, p_m^L}_L$, so
\begin{equation} \lb{2.54d}
f(x)=\sum_{m=0}^\infty f_m^{(L)} p_m^L(x),
\end{equation}
then, by \eqref{2.54b},
\begin{equation} \lb{2.54e}
\int f(z)\, d\mu(z) K_n(z,w) =\sum_{m=0}^n f_m^{(L)} p_m^L (w).
\end{equation}

One has
\begin{equation} \lb{2.54f}
\int K_n(z,w)\, d\mu(w) K_n(w,\zeta) =K_n (z,\zeta)
\end{equation}
as can be checked directly and which is an expression of the fact
that the map in \eqref{2.54cA} is a projection, and so its own
square.

We will let $\pi_n$ be the map of $L^2 (d\mu)$ to itself given by
\eqref{2.54cA} (or by \eqref{2.54e}). \eqref{2.54c} can then be
viewed (for $z,w\in\bbR$) as an expression of the integral kernel
of the commutator $[J,\pi_n]$, which leads to another proof of it
\cite{Sim-cd}.

Let $J_{n;F}$ be the finite $nl\times nl$ matrix obtained from $J$
by taking the top leftmost $n2$ blocks. It is the matrix of
$\pi_{n-1} M_x\pi_{n-1}$ where $M_x$ is multiplication by $x$ in
the $\{p_j^R\}_{j=0}^{n-1}$ basis. For $y\in\bbC$ and
$\gamma\in\bbC^l$, let $\varphi_{n,\gamma}(y)$ be the vector whose
components are
\begin{equation} \lb{2.54g}
(\varphi_{n,\gamma}(y))_j =p_{j-1}^L(y)\gamma
\end{equation}
for $j=1,2,\dots, n$. We claim that
\begin{equation} \lb{2.54h}
[(J_{n;F}-y)\varphi_{n,\gamma}(y)]_j = -\delta_{jn} A_{n} p_n^L(y)\gamma
\end{equation}
as follows immediately from \eqref{b8}.

This is intimately related to \eqref{b8.2} and \eqref{2.54c}. For
recalling $J$ is the matrix in
$p^R$ basis, $\varphi_{n,\gamma}(y)$ corresponds to the function
\[
\sum_{j=0}^{n-1} p_j^R(x) (\varphi_{n,\gamma}(y))_{j-1} = K_n (x,y)
\gamma .
\]
As we will see in the next two sections, \eqref{2.54h} has important
consequences.

\subsection{Christoffel Variational Principle} \lb{s2.7A}

There is a matrix version of the Christoffel variational principle
(see Nevai \cite{NevFr} for a
discussion of uses in the scalar case; this matrix case is discussed
by Duran--Polo \cite{DP}):

\begin{theorem}\lb{T2.18A} For any non-trivial $l\times l$ matrix-
valued measure, $d\mu$, on
$\bbR$, we have that for any $n$, any $x_0\in\bbR$, and matrix
polynomials $Q_n(x)$ of degree
at most $n$ with
\begin{equation}\lb{2.66a}
Q_n(x_0)=\bdone,
\end{equation}
we have that
\begin{equation}\lb{2.66b}
\ang{Q_n, Q_n}_R \geq K_n(x_0,x_0)^{-1}
\end{equation}
with equality if and only if
\begin{equation}\lb{2.66c}
Q_n(x) =K_n(x,x_0) K_n(x_0,x_0)^{-1}.
\end{equation}
\end{theorem}

\begin{remark} \eqref{2.66b} also holds for $\ang{\cdot,\cdot}_L$ but
the minimizer is then
$K_n(x_0,x_0)^{-1} K_n(x,x_0)$.
\end{remark}

\begin{proof} Let $Q_n^{(0)}$ denote the right-hand side of \eqref
{2.66c}. Then for any
polynomial $R_n$ of degree at most $n$, we have
\begin{equation}\lb{2.66d}
\ang{Q_n^{(0)}, R_n}_R = K_n(x_0,x_0)^{-1} R_n(x_0)
\end{equation}
because of \eqref{2.54cA}. Since $Q_n(x_0)=Q_n^{(0)}(x_0)=\bdone$, we
conclude
\begin{equation}\lb{2.66e}
\ang{Q_n-Q_n^{(0)},Q_n-Q_n^{(0)}}_R =\ang{Q_n,Q_n}_R - K_n(x_0,x_0)^{-1}
\end{equation}
from which \eqref{2.66b} is immediate and, given the supposed non-
triviality, the uniqueness of
minimizer.
\end{proof}

With this, one easily gets an extension of a result of M\'at\'e--Nevai \cite{MN} to MOPRL.
(They had it for scalar OPUC. For OPUC, it is in M\'at\'e--Nevai--Totik \cite{MNT88}
on $[-1,1]$ and in Totik \cite{Tot} for general OPRL. The result
below can be proven using
polynomial mappings \`a la Totik \cite{Tot-acta} or Jost solutions
\`a la Simon \cite{Weak-cd}.)

\begin{theorem} \lb{T2.18B} Let $d\mu$ be a non-trivial $l\times l$
matrix-valued measure on
$\bbR$ with compact support, $E$. Let $I=(a,b)$ be an open interval
with $I\subset E$. Then for
Lebesgue a.e.\ $x\in I$,
\begin{equation}\lb{2.66f}
\limsup (n+1) K_n(x,x)^{-1} \leq w(x).
\end{equation}
\end{theorem}

\begin{remark} This is intended in terms of expectations in any fixed
vector.
\end{remark}

We state this explicitly since we will need it in Section~\ref
{s5.3A}, but we note that
the more detailed results of M\'at\'e--Nevai--Totik \cite{MNT88},
Lubinsky \cite{Lub},
Simon \cite{2ext}, and Totik \cite{Tot-prep} also extend.

\subsection{Zeros} \lb{s2.8}

We next look at zeros of $\det(P_n^L(z))$, which we will prove soon
is also $\det(P_n^R(z))$. Following
\cite{DL,Sinap}, we will identify these zeros with eigenvalues of $J_
{n;F}$. It is not obvious a priori
that these zeros are real and, unlike the scalar situation, where the
classical arguments on the
zeros rely on orthogonality, we do not know how to get reality just
from that (but see the remark
after Theorem~\ref{T2.18E}).

\begin{lemma}\lb{L2.18A} Let $C(z)$ be an $l\times l$ matrix-valued
function analytic near $z=0$. Let
\begin{equation} \lb{2.54i}
k=\dim(\ker(C(0))).
\end{equation}
Then $\det(C(z))$ has a zero at $z = 0$ of order at least $k$.
\end{lemma}

\begin{remarks} 1. Even in the $1\times 1$ case, where $k=1$, the
zeros can clearly be of higher
order than $k$ since $c_{11}(z)$ can have a zero of any order!

\smallskip
2. The temptation to take products of eigenvalues will lead at best
to a complicated proof
as the cases $C(z)=\left( \begin{smallmatrix} 0&z\\ 1&0 \end
{smallmatrix}\right)$ and $C(z)
=\left(\begin{smallmatrix} 0&z^2 \\ 1&0 \end{smallmatrix}\right)$
illustrate.
\end{remarks}

\begin{proof} Let $e_1, \dots , e_l$ be an orthonormal basis with $e_1,
\dots , e_k\in \ker(C(0))$.
By Hadamard's inequality (see Bhatia \cite{Bha}),
\begin{align*}
\abs{\det(C(z))} &\leq \norm{C(z) e_1} \cdots \norm{C(z) e_l} \\
&\leq C\abs{z}^k
\end{align*}
since $\norm{C(z)e_j}\leq C\abs{z}$ if $j=1, \dots, k$ and $\norm{C
(z) e_j}\leq d$ for $j=k+1, \dots, l$.
\end{proof}

The following goes back at least to \cite{Dean-Martin}; see also
\cite{DL,SB07,Sinap,SinapVanAsche}.
\begin{theorem}\lb{T2.18BB} We have that
\begin{equation} \lb{2.54j}
\det_{\bbC^l}(P_n^L(z)) =\det_{\bbC^{nl}}(z-J_{n;F}).
\end{equation}
\end{theorem}

\begin{proof} By \eqref{2.54h}, if $\gamma$ is such that $p_n^L(y)
\gamma =0$, then $\varphi_{n,\gamma}(y)$ is an eigenvector for
$J_{n;F}$ with eigenvalue $y$. Conversely, if $ \varphi$ is an
eigenvector and $\gamma$ is defined as that vector in $\bbC^l$
whose components are the first $l$ components of $\varphi$, then a
simple inductive argument shows $\varphi=\varphi_{n,\gamma}(y)$
and then, by \eqref{2.54h} and the fact that $A_{n}$ is
invertible, we see that $p_n^L(y)\gamma=0$. This shows that for
any $y$,
\begin{equation} \lb{2.54k}
\dim\ker (P_n^L(y))=\dim\ker (J_{n;F}-y).
\end{equation}

By Lemma~\ref{L2.18A}, if $y$ is any eigenvalue of $J_{n;F}$ of
multiplicity $k$, then $\det(P_n^L(z))$ has a zero of order at
least $k$ at $y$. Now let us consider the polynomials in $z$ on
the left and right in \eqref{2.54j}. Since $J_{n;F}$ is Hermitian,
the sum of the multiplicities of the zeros on the right is $nl$.
Since the polynomial on the left is of degree $nl$, by a counting
argument it has the same set of zeros with the same multiplicities
as the polynomial on the right. Since both polynomials are monic,
they are equal.
\end{proof}

\begin{corollary}\lb{C2.18C} All the zeros of $\det(P_n^L(z))$ are
real. Moreover,
\begin{equation} \lb{2.54L}
\det(P_n^R(z))=\det(P_n^L(z)).
\end{equation}
\end{corollary}

\begin{proof} Since $J_{n;F}$ is Hermitian, all its eigenvalues are
real, so \eqref{2.54j} implies
the zeros of $\det(P_n^L(z))$ are real. Thus, since the polynomial is
monic,
\begin{equation} \lb{2.54m}
\ol{\det(P_n^L(\bar z))} =\det(P_n^L(z)).
\end{equation}
By Lemma~\ref{l.b1}(v), we have
\begin{equation} \lb{2.54n}
P_n^R(z)=P_n^L(\bar z)^\dagger
\end{equation}
since both sides are analytic and agree if $z$ is real. Thus,
\[
\det(P_n^R(z)) =\det(P_n^L(\bar z)^\dagger) =\ol{\det(P_n^L(\bar z))}
\]
proving \eqref{2.54L}.
\end{proof}

The following appeared before in \cite{SinapVanAsche}; see also
\cite {SB07}.
\begin{corollary} \lb{C2.18D} Let $\{x_{n,j}\}_{j=1}^{nl}$ be the
zeros of $\det(P_n^L(x))$ counting
multiplicity ordered by
\begin{equation} \lb{2.54o}
x_{n,1} \leq x_{n,2} \leq \cdots\leq x_{n,nl}.
\end{equation}
Then
\begin{equation} \lb{2.54p}
x_{n+1,j}\leq x_{n,j} \leq x_{n+1,j+l}.
\end{equation}
\end{corollary}

\begin{remarks} 1. This is interlacing if $l=1$ and replaces it for
general $l$.

\smallskip
2. Using $A_n$ invertible, one can show the inequalities in \eqref
{2.54p} are strict.
\end{remarks}

\begin{proof} The min-max principle \cite{RS4} says that
\begin{equation} \lb{2.84a}
x_{n,j} = \max_{\substack{ L\subset\bbC^{nl} \\ \dim(L)\leq j-1}}\,
\min_{\substack{f\in L^\perp \\ \norm{f}=1}}
\jap{f,J_{n;F}f}_{\bbC^{nl}}.
\end{equation}
If $P\colon \bbC^{(n+1)l} \to \bbC^{nl}$ is the natural projection, then
$\jap{Pf,J_{n+1;F}Pf}_{\bbC^{(n+1)l}} = \jap{Pf,J_{n;F} Pf}_{\bbC^{nl}}$
and as $L$ runs through all subspaces of $\bbC^{(n+1)l}$ dimension at most
$j-1$, $P[L]$ runs through all subspaces of dimension at most $j-1$ in $\bbC^{nl}$,
so \eqref{2.84a} implies $x_{n+1,j}\leq x_{n,j}$.
Using the same argument on $-J_{n;F}$ and
$-J_{n+1;F}$ shows $x_j(-J_{n;F}) \geq x_j (-J_{n+1;F})$. But $x_j (-
J_{n;F})=-x_{nl+1-j}
(J_{n;F})$ and $x_j(-J_{n+1;F})=-x_{(n+1)l+1-j} (J_{n+1;F})$. That
yields the other inequality.
\end{proof}

\subsection{Lower Bounds on $p$ and the Stieltjes--Weyl
Formula for $m$}\lb{s2.9}

Next, we want to exploit \eqref{2.54h} to prove uniform lower bounds
on $\norm{p_n^L(y)\gamma}$
when $y\notin \cvh(\sigma(J))$, the convex hull of the support of $J
$, and thereby uniform bounds
$\norm{p_n^L(y)^{-1}}$. We will then use that to prove that for $z
\notin\sigma(J)$, we have
\begin{equation} \lb{2.54q}
m(z)=\lim_{n\to\infty} -p_n^L(z)^{-1} q_n^L(z)
\end{equation}
the matrix analogue of a formula that spectral theorists associate
with Weyl's definition of the
$m$-function \cite{Weyl}, although for the discrete case, it goes
back at least to Stieltjes
\cite{Stie}.

We begin by mimicking an argument from \cite{EqMC}. Let $H=\cvh(\sigma
(J))=[c-D,c+D]$ with
\begin{equation} \lb{2.54r}
D=\tfrac12\, \diam(H).
\end{equation}
By the definition of $A_n$,
\begin{equation} \lb{2.54s}
\norm{A_n} =\norm{\ang{p_{n-1}^R, (x-c)p_n^R}_R}\leq D.
\end{equation}
Suppose $y\notin H$ and let
\begin{equation} \lb{2.54t}
d=\dist(y,H).
\end{equation}
By the spectral theorem, for any vector $\varphi\in\h_v$,
\begin{equation} \lb{2.54u}
\abs{\jap{\varphi, (J-y)\varphi}_{\h_v}}\geq d\norm{\varphi}^2.
\end{equation}

By \eqref{2.54h}, with $\varphi=\varphi_{n,\gamma}(y)$,
\begin{equation} \lb{2.54v}
\abs{\jap{\varphi, (J-y)\varphi}} \leq \norm{A_{n}}\, \norm{p_n^L(y)
\gamma}\,
\norm{p_{n-1}^L(y)\gamma}
\end{equation}
while
\begin{equation} \lb{2.54w}
\norm{\varphi}^2 =\sum_{j=0}^{n-1} \, \norm{p_j^L(y)\gamma}^2.
\end{equation}
As in \cite[Prop.~2.2]{EqMC}, we get:

\begin{theorem}\lb{T2.18E} If $y\notin H$, for any $\gamma$,
\begin{equation} \lb{2.54x}
\norm{p_n^L(y)\gamma}\geq \bigl(\tfrac{d}{D}\bigr)
\bigl(1+\bigl(\tfrac{d}{D}\bigr)^2\bigr)^{(n-1)/2}\norm{\gamma}.
\end{equation}
In particular,
\begin{equation} \lb{2.54y}
\norm{p_n^L(y)^{-1}} \leq \tfrac{D}{d}\, .
\end{equation}
\end{theorem}

\begin{remark} \eqref{2.54x} implies $\det(p_n^L(y))\neq 0$ if $\Ima
y >0$, providing
another proof that its zeros are real.
\end{remark}

By the discussion after \eqref{b12}, if $\Ima z>0$, $q_n^L(z) + p_n^L
(z)m(z)$ is in $\ell^2$,
so goes to zero. Since $p_n^L(z)^{-1}$ is bounded, we obtain:

\begin{corollary}\lb{C2.18F} For any $z\in\bbC_+=\{z : \Ima
z>0\}$,
\begin{equation} \lb{2.54z}
m(z)=\lim_{n\to\infty} -p_n^L(z)^{-1} q_n^L(z).
\end{equation}
\end{corollary}

\begin{remark} This holds for $z\notin H$.
\end{remark}

Taking adjoints using \eqref{2.54n} and $m(z)^\dagger=m(\bar z)$, we
see that
\begin{equation} \lb{2.54aa}
m(z)=\lim_{n\to\infty}  -q_n^R(z) p_n^R(z)^{-1}.
\end{equation}

\eqref{2.54z} and \eqref{2.54aa} are due to \cite{Dur96}, which uses the
proof based on the Gauss--Jacobi quadrature formula.

\subsection{Wronskians of Vector-Valued Solutions}\lb{s2.10}

Let $\alpha,\beta$ be two vector-valued solutions $(\bbC^l)$ of
\eqref {b16} for $n=2,3,\dots$. Define their scalar Wronskian as
(Euclidean inner product on $\bbC^l $)
\begin{equation} \lb{2.54bb}
W_n(\alpha,\beta) =\jap{\alpha_{n}, A_{n}\beta_{n+1}} - \jap{A_{n}
\alpha_{n+1}, \beta_{n}}
\end{equation}
for $n=2,3,\dots$. One can obtain two matrix solutions by using $
\alpha$ or $\beta$ for one column
and $0$ for the other columns. The scalar Wronskian is just a matrix
element of the resulting
matrix Wronskian, so $W_n$ is constant (as can also be seen by direct
calculation). Here is an
application:

\begin{theorem}\lb{T2.18F} Let $z_0 \in \R \setminus \sigma_\ess(J)$.
For $k = 0,1$, let $q_k$ be the
multiplicity of $z_0$ as an eigenvalue of $J^{(k)}$. Then, $q_0 + q_1
\le l$.
\end{theorem}

\begin{proof} If $\ti\beta$ is an eigenfunction for $J^{(1)}$ and we
define $\beta$ by
\begin{equation} \lb{2.54cc}
\beta_n =\begin{cases} 0, & n=1, \\
\ti\beta_{n-1}, & n\geq 2,
\end{cases}
\end{equation}
then $\beta$ solves \eqref{b16} for $n\geq2$. If $\alpha$ is an
eigenfunction for $J=J^{(0)}$, it also solves \eqref{b16}. Since
$\alpha_n\to 0$, $\beta_n\to 0$, and $A_n$ is bounded,
$W_n(\alpha,\beta)\to 0$ as $n\to\infty$ and so it is identically
zero. But since $\beta_1 =0 $,
\begin{equation} \lb{2.54dd}
\bdzero=W_1(\alpha,\beta) =\jap{\alpha_1, A_1\beta_2} = \jap
{\alpha_1, A_1\ti\beta_1}.
\end{equation}

Let $V^{(k)}$ be the set of values of eigenfunctions of $J^{(k)}$ at
$n=1$. \eqref{2.54dd}
says
\begin{equation} \lb{2.54ee}
V^{(0)}\subset [A_1 V^{(1)}]^\perp .
\end{equation}
Since $q_k=\dim(V^{(k)})$ and $A_1$ is non-singular, \eqref{2.54ee}
implies that $q_0\leq l-q_1$.
\end{proof}

\subsection{The Order of Zeros/Poles of $m(z)$} \lb{s2.11}

\begin{theorem}\lb{T2.18H} Let $z_0 \in \R \setminus \sigma_\ess(J)$.
For $k = 0,1$, let $q_k$ be the
multiplicity of $z_0$ as an eigenvalue of $J^{(k)}$. If $q_1 - q_0
\geq0$, then $\det(m(z))$ has a
zero at $z=z_0$ of order $q_1 - q_0$. If $q_1-q_0<0$, then  $\det(m
(z))$ has a pole at $z=z_0$
of order $q_0 - q_1$.
\end{theorem}

\begin{remarks} 1. To say $\det(m(z))$ has a zero at $z=z_0$ of order
$0$ means it is finite
and non-vanishing at $z_0$!

\smallskip
2. Where $d\mu$ is a direct sum of scalar measures, so is $m(z)$, and
$\det(m(z))$ is then a product.
In the scalar case, $m(z)$ has a pole at $z_0$ if $J^{(0)}$ has $z_0$
as eigenvalue and a zero at
$z_0$ if $J^{(1)}$ has $z_0$ as eigenvalue. In the direct sum case,
we see there can be
cancellations, which helps explain why $q_1-q_0$ occurs.

\smallskip
3. Formally, one can understand this theorem as follows. Cramer's
rule suggests $\det(m(z))=
\det (J^{(1)}-z)/\det(J^{(0)}-z)$. Even though $\det(J^{(k)}-z)$ is
not well-defined in the
infinite case, we expect a cancellation of zeros of order $q_1$ and
$q_0$. For $z_0\notin H$,
the convex hull of $\sigma(J^{(0)})$, one can use \eqref{2.54z} to
prove the theorem following
this intuition. Spurious zeros in gaps of $\sigma (J^{(0)})$ make
this strategy difficult in gaps.

\smallskip
4. Unlike in Lemma~\ref{L2.18A}, we can write $m$ as a product of
eigenvalues and analyze that directly because $m(x)$ is
self-adjoint for $x$ real, which sidesteps some of the problems
associated with non-trivial Jordan normal forms.

\smallskip
5. This proof gives another demonstration that $q_0 + q_1\leq l$.
\end{remarks}

\begin{proof} $m(z)$ has a simple pole at $z=z_0$ with a residue
which is rank $q_0$ and
strictly negative definite on its range. Let
\[
f(z)=(z-z_0) m(z).
\]
$f$ is analytic near $z_0$ and self-adjoint for $z$ real and near
$z_0$. Thus, by eigenvalue perturbation theory \cite{Kato,RS4},
$f(z)$ has $l$ eigenvalues $\rho_1(z), \dots, \rho_l(z)$ analytic
near $z_0$ with $\rho_1, \rho_2, \dots, \rho_{q_0}$ non-zero at
$z_0$ and $\rho_{q_0+1}, \dots, \rho_l$ zero at $z_0$.

Thus, $m(z)$ has $l$ eigenvalues near $z_0$,
$\lambda_j(z)=\rho_j(z)/(z-z_0)$, where $\lambda_1, \dots,
\lambda_{q_0}$ have simple poles and the others are regular.

By Proposition~\ref{strip}, $m(z)^{-1}$ has a simple pole at $z_0$
with residue of rank $q_1$ (because $A_1$ is non-singular),
so $m(z)^{-1}$ has
$q_1$ eigenvalues with poles. That means $q_1$ of
$\lambda_{q_0+1}, \dots, \lambda_l$ have simple zeros at $z_0$ and
the others are non-zero. Thus, $\det(m(z))=\prod_{j=1}^l
\lambda_j(z)$ has a pole/zero of order $q_0-q_1$.
\end{proof}

\subsection{Resolvent of the Jacobi Matrix} Consider the matrix
$$
G_{nm}(z)= \ang{p_{n-1}^R, (x-z)^{-1}p_{m-1}^R}_R.
$$

\begin{theorem} \lb{T2.27}
One has
\begin{equation}
G_{nm}(z)=
\begin{cases}
\psi_{n-1}^L(z) p_{m-1}^R(z), & \text{ if $n\geq m$},
\\
p_{n-1}^L(z) \psi_{m-1}^R(z), & \text{ if $n\leq m$}.
\end{cases}
\end{equation}
\end{theorem}
\begin{proof}
We have
\begin{alignat*}{2}
\sum_{m=1}^\infty G_{km}(z) J_{mn} &= z G_{kn}(z),
\quad & n &\not= k,
\\
\sum_{m=1}^\infty  J_{nm} G_{mk}(z) &= z G_{nk}(z),
\quad &  n& \not = k.
\end{alignat*}
Fix $k\geq0$ and let $u_m(z)=G_{km}(z)$. Then $u_m(z)$ satisfies the
equation \eqref{b15}
for $n\not =k$, and so we can use Theorem~\ref{t.b1} to describe this
solution.
First suppose $k>1$. As $u_m$ is an $\ell^2$ solution and $u_m$
satisfies \eqref{b15} for
$n=1$, we have
\begin{equation}
G_{km}(z)=
\begin{cases}
a_k(z) p_{m-1}^R(z), & m\leq k,
\\
b_k(z) \psi_{m-1}^R(z), & m\geq k.
\end{cases}
\label{b20}
\end{equation}
If $k=1$, \eqref{b20} also holds true. For $m\geq k$, this follows by
the same
argument, and for $m=k=1$, this is a trivial statement.
Next, similarly, let us consider $v_m(z)=G_{mk}(z)$.
Then $v_m(z)$ solves \eqref{b16} and so, using Theorem~\ref{t.b2}, we
obtain
\begin{equation}
G_{mk}(z)=
\begin{cases}
p_{m-1}^L(z) c_k(z), & m\leq k,
\\
\psi_{m-1}^L(z) d_k(z) , & m\geq k.
\end{cases}
\label{b21}
\end{equation}
Comparing \eqref{b20} and \eqref{b21}, we find
\begin{align*}
a_k(z) p_{m-1}^R(z) & =  \psi_{k-1}^L(z) d_m(z),
\\
b_k(z) \psi_{m-1}^R(z) & =  p_{k-1}^L(z) c_m(z).
\end{align*}
As $p_0^R=p_0^L={\boldsymbol 1}$,
it follows that
$$
a_k(z)=\psi_{k-1}^L(z) d_1(z),
\qquad
c_m(z)=b_1(z) \psi_{m-1}^R(z)
$$
and so we obtain
\begin{equation}
G_{nm}(z)=
\begin{cases}
\psi_{n-1}^L(z)d_1(z) p_{m-1}^R(z) & \text{ if $n\geq m$},
\\
p_{n-1}^L(z) b_1(z) \psi_{m-1}^R(z) & \text{ if $n\leq m$}.
\end{cases}
\label{b22}
\end{equation}
It remains to prove that
\begin{equation}
b_1(z)=d_1(z)={\boldsymbol 1}.
\label{b23}
\end{equation}
Consider the case $m=n=1$.
By the definition of the resolvent,
$$
G_{11}(z)
=
\ang{p_0^R,(J-z)^{-1}p_0^R}_R
=
\int \frac{d\mu(x)}{x-z}
=
m(z).
$$
On the other hand, by \eqref{b22},
\begin{align*}
G_{11}(z)&=\psi_0^L(z)d_1(z) p_0^R(z)=m(z) d_1(z),
\\
G_{11}(z)&=p_0^L(z) b_1(z) \psi_0^R(z)=b_1(z) m(z),
\end{align*}
which proves \eqref{b23}.
\end{proof}

\section{Matrix Orthogonal Polynomials on the Unit Circle} \lb{s3}

\subsection{Definition of MOPUC} \lb{s3.1}

In this chapter, $\mu$ is an $l\times l$ matrix-valued measure on $
\partial\bbD$.
$\ang{\cdot,\cdot}_R$ and $\ang{\cdot,\cdot}_L$ are defined as in
the MOPRL case. Non-triviality is defined as for MOPRL. We will
always assume $\mu$ is non-trivial.  We define monic matrix
polynomials $\Phi_n^R, \Phi_n^L$ by applying Gram--Schmidt to
$\{\bdone, z\bdone, \dots\}$, that is, $\Phi_n^R$ is the unique
matrix polynomial $z^n \bdone +$ lower order with
\begin{equation} \lb{3.1}
\ang{z^k\bdone, \Phi_n^R}_R =0 \quad k=0,1,\dots, n-1.
\end{equation}
We will define the normalized MOPUC shortly. We will only consider
the analogue of what we called type~1 for MOPRL because only those
appear to be useful. Unlike in the scalar case, the monic
polynomials do not appear much because it is for the normalized,
but not monic, polynomials that the left and right Verblunsky
coefficients are the same.

\subsection{The Szeg\H{o} Recursion}

Szeg\H{o} \cite{Szb} included the scalar Szeg\H{o} recursion for the
first time. It seems likely
that Geronimus had it independently shortly after Szeg\H{o}. Not
knowing of the connection
with this work, Levinson \cite{Lev} rederived the recursion but with
matrix coefficients!
So the results of this section go back to 1947.

For a matrix polynomial $P_n$ of degree $n$, we define the reversed
polynomial $P_n^*$ by
\begin{equation} \lb{3.1a}
P_n^*(z) = z^n P_n(1/\bar{z})^\dagger.
\end{equation}
Notice that
\begin{equation} \lb{3.2x}
(P_n^*)^* = P_n
\end{equation}
and for any $\alpha\in\calM_l$,
\begin{equation} \lb{3.2}
(\alpha P_n)^* = P_n^* \alpha^\dagger, \qquad (P_n\alpha)^* = \alpha^
\dagger P_n^*.
\end{equation}

\begin{lemma}
We have
\begin{equation} \lb{3.1ax}
\ang{f , g}_L = \ang{g , f}_L^\dagger \, ,
\qquad
\ang{f , g}_R = \ang{g , f}_R^\dagger
\end{equation}
and
\begin{equation} \lb{3.1b}
\ang{f^* , g^*}_L = \ang{f , g}_R^\dagger \, ,
\quad
\ang{f ^*, g^*}_R = \ang{f ,g}_L^\dagger \, .
\end{equation}
\end{lemma}

\begin{proof}
The first and second identities  follow immediately from the
definition. The third identity is
derived as follows:
\begin{align*}
\ang{f^* , g^*}_L & = \int e^{in\theta} g(\theta)^\dagger \,
d\mu(\theta) \, (e^{in\theta} f(\theta)^\dagger)^\dagger \\
& = \int e^{in\theta} g(\theta)^\dagger \, d\mu(\theta) \, e^{-in
\theta} f(\theta) \\
& = \int g(\theta)^\dagger \, d\mu(\theta) \, f(\theta) \\
& = \ang{g , f}_R \\
& = \ang{f , g}_R^\dagger .
\end{align*}
The proof of the last identity is analogous.
\end{proof}

\begin{lemma}\label{l.c1}
If $P_n$ has degree $n$ and is left-orthogonal with respect to $z
\bdone, \ldots , z^n \bdone$, then
$P_n = c (\Phi_n^R)^*$ for some suitable matrix $c$.
\end{lemma}

\begin{proof}
By assumption,
$$
\bdzero = \ang{P_n , z^j \bdone}_L = \ang{(z^j \bdone)^* , P_n^*}_R =
\ang{z^{n-j} \bdone , P_n^*}_R \; \text{ for } 1
\le j \le n.
$$
Thus, $P_n^*$ is right-orthogonal with respect to $\bdone, z\bdone,
\ldots, z^{n-1} \bdone$ and hence it
is a right-multiple of $\Phi_n^R$. Consequently, $P_n$ is a left-
multiple of
$(\Phi_n^R)^*$.
\end{proof}

Let us define normalized orthogonal matrix polynomials by
$$
\varphi_0^L=\varphi_0^R=\bdone,\quad
\varphi_n^L = \kappa_n^L \Phi_n^L \quad \text{ and } \quad \varphi_n^R =
\Phi_n^R\kappa_n^R
$$
where the $\kappa$'s are defined according to the normalization
condition
$$
\ang{\varphi_n^R,\varphi_m^R}_R=\delta_{nm}{\boldsymbol 1}
\quad
\ang{\varphi_n^L,\varphi_m^L}_L=\delta_{nm}{\boldsymbol 1}
$$
along with (a type~1 condition)
\begin{equation}\lb{3.4a}
\kappa_{n+1}^L (\kappa_n^L)^{-1} > \bdzero \quad \text{and} \quad
(\kappa_n^R)^{-1}
\kappa_{n+1}^R > \bdzero.
\end{equation}
Notice that $\kappa_{n}^L$ are determined by the normalization
condition up to multiplication on the left by unitaries; these
unitaries can always be uniquely chosen so as to satisfy
\eqref{3.4a}.

Now define
$$
\rho_n^L = \kappa_n^L (\kappa_{n+1}^L)^{-1} \quad \text{and} \quad
\rho_n^R =
(\kappa_{n+1}^R)^{-1} \kappa_n^R.
$$
Notice that as inverses of positives matrices, $\rho_n^L >0$ and $
\rho_n^R >0$.
In particular, we have that
$$
\kappa_n^L = (\rho_{n-1}^L \dots \rho_0^L)^{-1} \quad \text{and}
\quad \kappa_n^R =
(\rho_0^R \dots \rho_{n-1}^R)^{-1}.
$$

\begin{theorem}[Szeg\H{o} Recursion]\label{szethm}
\begin{SL}
\item[{\rm{(a)}}] For suitable matrices $\alpha_n^{L,R}$, one has
\begin{align}
z \varphi_n^L - \rho_n^L \varphi_{n+1}^L & = (\alpha_n^L)^\dagger
\varphi_n^{R,*},  \lb{3.5a}\\
z \varphi_n^R - \varphi_{n+1}^R \rho_n^R & = \varphi_n^{L,*}
(\alpha_n^R)^\dagger. \lb{3.5b}
\end{align}

\item[{\rm{(b)}}] The matrices $\alpha_n^L$ and $\alpha_n^R$ are
equal and will henceforth be
denoted by $\alpha_n$.\\

\item[{\rm{(c)}}] $\rho_n^L = (\bdone - \alpha_n^\dagger \alpha_n)^
{1/2}$ and $\rho_n^R = (\bdone -
\alpha_n \alpha_n^\dagger)^{1/2}$.
\end{SL}
\end{theorem}

\begin{proof}
(a) The matrix polynomial $z \varphi_n^L$ has leading term $z^{n+1}
\kappa_n^L$. On the
other hand, the matrix polynomial $\rho_n^L \varphi_{n+1}^L$ has
leading term $z^{n+1}
\rho_n^L \kappa_{n+1}^L$. By definition of $\rho_n^L$, these terms
are equal.
Consequently, $z \varphi_n^L - \rho_n^L \varphi_{n+1}^L$ is a matrix
polynomial of degree
at most $n$. Notice that it is left-orthogonal with respect to $z
\bdone, \ldots , z^n \bdone$ since
$$
\ang{z \varphi_n^L - \rho_n^L \varphi_{n+1}^L , z^j \bdone}_L =
\ang{\varphi_n^L , z^{j-1} \bdone}_L - \ang{\rho_n^L \varphi_{n+1}
^L , z^j \bdone}_L
= \bdzero - \bdzero = \bdzero.
$$
Now apply Lemma~\ref{l.c1}. The other claim is proved in the same way.\\

\smallskip
(b) By part (a) and identities established earlier,
\begin{align*}
(\alpha_n^L)^\dagger & = \bdzero + (\alpha_n^L)^\dagger \bdone \\
& = \ang{\varphi_n^{R,*} , \rho_n^L \varphi_{n+1}^L}_L +
(\alpha_n^L)^\dagger \ang{\varphi_n^R ,  \varphi_n^R} _R \\
& = \ang{\varphi_n^{R,*} , \rho_n^L \varphi_{n+1}^L}_L +
(\alpha_n^L)^\dagger \ang{\varphi_n^{R,*} ,\varphi_n^{R,*}}_L \quad
\text{by \eqref{3.1b}} \\
& = \ang{\varphi_n^{R,*} , \rho_n^L \varphi_{n+1}^L +
(\alpha_n^L)^\dagger \varphi_n^{R,*}}_L \\
& = \ang{\varphi_n^{R,*} , z \varphi_n^L}_L \\
& = \ang{z \varphi_n^R , \varphi_n^{L,*}}_R^\dagger  \quad \text
{(using the $(n+1)$-degree *)}\\
& = \ang{\varphi_{n+1}^R \rho_n^R + \varphi_n^{L,*} (\alpha_n^R)^
\dagger ,
\varphi_n^{L,*}}_R^\dagger \\
& = \ang{\varphi_{n+1}^R \rho_n^R , \varphi_n^{L,*}}_R^\dagger
+ \ang{\varphi_n^{L,*} (\alpha_n^R)^\dagger , \varphi_n^{L,*}}_R^
\dagger \\
& = \bdzero + \ang{\varphi_n^{L,*} , \varphi_n^{L,*}(\alpha_n^R)^
\dagger}_R \\
& = \ang{\varphi_n^{L,*} , \varphi_n^{L,*}}_R \, (\alpha_n^R)^\dagger \\
& = \ang{\varphi_n^L , \varphi_n^L}_L \,
(\alpha_n^R)^\dagger \\
& = (\alpha_n^R)^\dagger.
\end{align*}
(c) Using parts (a) and (b), we see that
\begin{align*}
\bdone & = \ang{z \varphi_n^L , z \varphi_n^L}_L \\
& = \ang{\rho_n^L \varphi_{n+1}^L + \alpha_n^\dagger \varphi_n^{R,*} ,
\rho_n^L \varphi_{n+1}^L + \alpha_n^\dagger \varphi_n^{R,*}}_L \\
& = \rho_n^L \ang{\varphi_{n+1}^L , \varphi_{n+1}^L}_L
(\rho_n^L)^\dagger + \alpha_n^\dagger \ang{\varphi_n^{R,*} ,
\varphi_n^{R,*}}_L \, \alpha_n \\
& = (\rho_n^L)2 + \alpha_n^\dagger \ang{\varphi_n^R , \varphi_n^R}_R
\, \alpha_n \\
& = (\rho_n^L)2 + \alpha_n^\dagger \alpha_n.
\end{align*}
A similar calculation yields the other claim.
\end{proof}

The matrices $\alpha_n$ will henceforth be called the Verblunsky
coefficients associated
with the measure $d\mu$. Since $\rho_n^L$ is invertible, we have
\begin{equation} \lb{3.3a}
\norm{\alpha_n}<1.
\end{equation}
We will eventually see (Theorem~\ref{T3.10B}) that any set of $
\alpha_n$'s obeying \eqref{3.3a}
occurs as the set of Verblunsky coefficients for a unique non-trivial
measure.

Note that the Szeg\H{o} recursion for the monic orthogonal
polynomials is
\begin{equation} \label{szmonic}
\begin{aligned}
z \Phi_n^L -  \Phi_{n+1}^L & = (\kappa_n^L)^{-1} \alpha_n^\dagger
(\kappa_n^R)^\dagger\Phi_n^{R,*}, \\
z \Phi_n^R - \Phi_{n+1}^R  & = \Phi_n^{L,*} (\kappa_n^L)^\dagger
\alpha_n^\dagger (\kappa_n^R)^{-1},
\end{aligned}
\end{equation}
so the coefficients in the $L$ and $R$ equations are not equal and
depend on all the
$\alpha_j$, $j=1,\dots, n$.

Let us write the Szeg\H{o} recursion in matrix form, starting with
left-orthogonal polynomials. By Theorem \ref{szethm},
\begin{align*}
\varphi_{n+1}^L & = (\rho_n^L)^{-1} [z \varphi_n^L - \alpha_n^\dagger
\varphi_n^{R,*} ], \\
\varphi_{n+1}^R & = [ z \varphi_n^R - \varphi_n^{L,*} \alpha_n^\dagger ]
( \rho_n^R )^{-1},
\end{align*}
which implies that
\begin{align}
\varphi_{n+1}^L & = z (\rho_n^L)^{-1} \varphi_n^L - (\rho_n^L)^{-1}
\alpha_n^\dagger \varphi_n^{R,*}, \lb{3.7a} \\
\varphi_{n+1}^{R,*} & = (\rho_n^R)^{-1} \varphi_n^{R,*} - z (\rho_n^R)
^{-1}
\alpha_n \varphi_n^L. \lb{3.7b}
\end{align}
In other words,
\begin{equation}
\begin{pmatrix} \varphi_{n+1}^L \\ \varphi_{n+1}^{R,*} \end{pmatrix}
= A^L(\alpha_n,z) \begin{pmatrix} \varphi_n^L \\ \varphi_n^{R,*} \end
{pmatrix}
\label{*1}
\end{equation}
where
$$
A^L(\alpha,z) = \begin{pmatrix} z (\rho^L)^{-1} & - (\rho^L)^{-1}
\alpha^\dagger \\
- z (\rho^R)^{-1} \alpha & (\rho^R)^{-1}
\end{pmatrix}
$$
and $\rho^L = (\bdone - \alpha^\dagger \alpha)^{1/2}$, $\rho^R =
(\bdone - \alpha
\alpha^\dagger)^{1/2}$. Note that, for $z \not= 0$, the inverse of
$A^L(\alpha,z)$ is
given by
$$
A^L(\alpha,z)^{-1} =
\begin{pmatrix}
z^{-1} (\rho^L)^{-1} & z^{-1} (\rho^L)^{-1} \alpha^\dagger \\
(\rho^R)^{-1} \alpha & (\rho^R)^{-1}
\end{pmatrix}
$$
which gives rise to the inverse Szeg\H{o} recursion (first
emphasized in the scalar and matrix cases by Delsarte el al.\
\cite{DGK})
\begin{align*}
\varphi_n^L & =  z^{-1} (\rho_n^L)^{-1} \varphi_{n+1}^L + z^{-1}
(\rho_n^L)^{-1} \alpha_n^\dagger \varphi_{n+1}^{R,*},
\\
\varphi_n^{R,*} & =  (\rho_n^R)^{-1} \alpha_n\varphi_{n+1}^L +
(\rho_n^R)^{-1} \varphi_{n+1}^{R,*}.
\end{align*}

For right-orthogonal polynomials, we find the following matrix
formulas. By Theorem~\ref{szethm},
\begin{align}
\varphi_{n+1}^R & = z \varphi_n^R (\rho_n^R)^{-1} - \varphi_n^{L,*}
\alpha_n^\dagger (\rho_n^R)^{-1},  \label{szright1}
\\
\varphi_{n+1}^{L,*} & = \varphi_n^{L,*} (\rho_n^L)^{-1} - z \varphi_n^R
\alpha_n (\rho_n^L)^{-1}.  \label{szright2}
\end{align}
In other words,
$$
\begin{pmatrix} \varphi_{n+1}^R & \varphi_{n+1}^{L,*} \end{pmatrix}  =
\begin{pmatrix} \varphi_n^R & \varphi_n^{L,*} \end{pmatrix} A^R
(\alpha_n,z)
$$
where
$$
A^R(\alpha,z) = \begin{pmatrix}
z (\rho^R)^{-1} & - z \alpha (\rho^L)^{-1} \\
- \alpha^\dagger (\rho^R)^{-1} & (\rho^L)^{-1}
\end{pmatrix}.
$$
For $z \not= 0$, the inverse of $A^R(\alpha,z)$ is given by
$$
A^R(\alpha,z)^{-1} = \begin{pmatrix}
z^{-1} (\rho^R)^{-1} & (\rho^R)^{-1} \alpha \\
z^{-1} (\rho^L)^{-1} \alpha^\dagger & (\rho^L)^{-1}
\end{pmatrix}
$$
and hence
\begin{align}
\varphi_n^R & =  z^{-1} \varphi_{n+1}^R (\rho_n^R)^{-1} + z^{-1}
\varphi_{n+1}^{L,*} (\rho_n^L)^{-1} \alpha_n^\dagger,  \label{szright3}
\\
\varphi_n^{L,*} & =  \varphi_{n+1}^R (\rho_n^R)^{-1} \alpha_n +
\varphi_{n+1}^{L,*} (\rho_n^L)^{-1}. \label{szright4}
\end{align}

\subsection{Second Kind Polynomials} \lb{s3.3}

In the scalar case, second kind polynomials go back to Geronimus \cite
{Ger44,GBk1,GBk}.
For $n\geq1$, let us introduce the second kind polynomials $\psi_n^
{L,R}$ by
\begin{align}
\psi_n^L (z)&= \int \frac{e^{i\theta}+z}{e^{i\theta}-z}\,
(\varphi_n^L ( e^{i\theta} )
-\varphi_n^L (z) ) \, d\mu (\theta), \label{a1}
\\
\psi_n^R (z)&= \int \frac{e^{i\theta}+z}{e^{i\theta}-z} \, d\mu
(\theta) \, (\varphi_n^R
( e^{i\theta} ) -\varphi_n^R (z) ) . \label{a2}
\end{align}
For $n=0$, let us set $\psi_0^L (z)=\psi_0^R (z)=\bdone$. For
future reference, let us display the first polynomials of each
series:
\begin{align}
\varphi_1^L(z)&=(\rho_0^L)^{-1}(z-\alpha_0^\dagger),
\qquad
\varphi_1^R(z)=(z-\alpha_0^\dagger)(\rho_0^R)^{-1},
\label{firstpoly1}
\\
\psi_1^L(z)&=(\rho_0^L)^{-1}(z+\alpha_0^\dagger),
\qquad
\psi_1^R(z)=(z+\alpha_0^\dagger)(\rho_0^R)^{-1}.
\label{firstpoly2}
\end{align}
We will also need formulas for
$\psi_n^{L,*}$ and $\psi_n^{R,*}$, $n\geq1$. These formulas follow
directly from the
above definition and from
$$
\overline{\left(\frac{e^{i\theta}+1/\bar z}{e^{i\theta}-1/\bar z}
\right)} =
-\frac{e^{i\theta}+z}{e^{i\theta}-z}\, .
$$
Indeed, we have
\begin{align}
\psi_n^{L,*} (z)&= z^n \int \frac{e^{i\theta}+z}{e^{i\theta}-z} \, d
\mu (\theta) \,
(\varphi_n^L (1/\bar z )^\dagger -\varphi_n^L (e^{i\theta})^\dagger),
\label{a2a}
\\
\psi_n^{R,*} (z)&= z^n \int \frac{e^{i\theta}+z}{e^{i\theta}-z}\,
(\varphi_n^R
(1/\bar z )^\dagger -\varphi_n^R (e^{i\theta})^\dagger) \, d\mu
(\theta).
\label{a3}
\end{align}

\begin{proposition}
The second kind polynomials obey the recurrence relations
\begin{align}
\psi_{n+1}^L(z)&=(\rho_n^L)^{-1}(z\psi_n^L(z)+\alpha_n^\dagger \psi_
{n}^{R,*}(z)),
\label{a4}
\\
\psi_{n+1}^{R,*}(z)&=(\rho_n^R)^{-1}(z\alpha_n \psi_{n}^L(z)+\psi_{n}^
{R,*}(z))
\label{a5}
\end{align}
and
\begin{align}
\psi_{n+1}^R(z) & = (z \psi_n^R(z) + \psi_n^{L,*}(z) \alpha_n^
\dagger) (\rho_n^R)^{-1}, \label{a4r}
\\
\psi_{n+1}^{L,*}(z) & = (\psi_n^{L,*}(z) + z \psi_n^R(z) \alpha_n)
(\rho_n^L)^{-1} \label{a5r}
\end{align}
for $n\geq0$.
\end{proposition}
\begin{proof}
1. Let us check \eqref{a4} for $n\geq1$. Denote the right-hand side
of \eqref{a4} by
$\tilde\psi_{n+1}^L(z)$. Using the recurrence relations for $
\varphi_n^L$, $\varphi_n^{R,*}$ and
the definition \eqref{a1} of $\psi_n^L$, we find
$$
\psi_{n+1}^L(z)- \tilde \psi_{n+1}^L(z) = \int \frac{e^{i\theta}+z}{e^
{i\theta}-z} \,
A_n(\theta,z) \, d\mu(\theta)
$$
where
\begin{align*}
A_n(\theta,z) & = \varphi_{n+1}^L(e^{i\theta}) -\varphi_{n+1}^L(z)
\\
& \qquad - (\rho_n^L)^{-1}[z \varphi_{n}^L(e^{i\theta}) -z\varphi_{n}
^L(z) +
\alpha_n^\dagger z^n \varphi_n^R(1/\bar z)^\dagger
- \alpha_n^\dagger z^n \varphi_n^R(e^{i\theta})^\dagger]
\\
& = (\rho_n^L)^{-1} [ e^{i\theta}\varphi_{n}^L(e^{i\theta}) -
\alpha_n^\dagger\varphi_{n}^{R,*}(e^{i\theta}) -
z\varphi_{n}^L(z)+\alpha_n^\dagger\varphi_{n}^{R,*}(z)
\\
& \qquad - z \varphi_{n}^L(e^{i\theta}) +z\varphi_{n}^L(z) - \alpha_n^
\dagger
\varphi_n^{R,*}(z) + \alpha_n^\dagger z^n
\varphi_n^R(e^{i\theta})^\dagger ]
\\
& = (\rho_n^L)^{-1}[(e^{i\theta}-z)\varphi_{n}^L(e^{i\theta}) +
\alpha_n^\dagger(z^n
e^{-in\theta}-1)\varphi_n^{R,*}(e^{i\theta})].
\end{align*}
Using the orthogonality relations
\begin{equation}
\int \varphi_n^L(e^{i\theta})\, d\mu(\theta) e^{-im\theta} = \int
\varphi_n^{R,*}(e^{i\theta})\, d\mu(\theta) e^{-i(m+1)\theta} =
\bdzero, \label{a6}
\end{equation}
$m=0,1,\dots,n-1$, and the formula
$$
\frac{e^{in\theta}-z^n}{e^{i\theta}-z} =
e^{i(n-1)\theta}+e^{i(n-2)\theta}z+\dots+z^{n-1}
$$
we obtain
\begin{align*}
\rho_n^L \int \frac{e^{i\theta}+z}{e^{i\theta}-z}\, A_n(\theta,z)\, d
\mu(\theta) & = \int
(e^{i\theta}\varphi_n^L(e^{i\theta}) - \alpha_n^\dagger
\varphi_n^{R,*}(e^{i\theta}))\, d\mu(\theta) \\
& = \rho_n^L\int \varphi_{n+1}^L (e^{i\theta}) \, d\mu(\theta) =
\bdzero.
\end{align*}

\smallskip
2. Let us check \eqref{a5} for $n\geq1$. Denote the right-hand side
of \eqref{a5} by
$\tilde\psi_{n+1}^{R,*}(z)$. Similarly to the argument above, we find
$$
\psi_{n+1}^{R,*}(z)- \tilde \psi_{n+1}^{R,*}(z) = \int
\frac{e^{i\theta}+z}{e^{i\theta}-z} \, B_n(\theta,z) \, d\mu(\theta),
$$
where
\begin{align*}
B_n(\theta,z) & = z^{n+1}\varphi_{n+1}^R(1/\bar z)^\dagger -
z^{n+1}\varphi_{n+1}^R(e^{i\theta})^\dagger \\
& \qquad - (\rho_n^R)^{-1} [ z\alpha_n\varphi_n^L(e^{i\theta}) - z
\alpha_n\varphi_n^L(z)
+
z^n\varphi_n^R(1/\bar z)^\dagger - z^n\varphi_n^R(e^{i\theta})^
\dagger ] \\
& = (\rho_n^R)^{-1} [ \varphi_n^{R,*}(z) - \alpha_n z \varphi_n^L(z) -
z^{n+1}e^{-i(n+1)\theta}
\varphi_n^{R,*}(e^{i\theta}) + z^{n+1}e^{-in\theta}\alpha_n
\varphi_n^L(e^{i\theta}) \\
& \qquad - z \alpha_n \varphi_n^L(e^{i\theta}) + z \alpha_n
\varphi_n^L(z) -
\varphi_n^{R,*}(z) +
z^ne^{-in\theta} \varphi_n^{R,*}(e^{i\theta}) ] \\
& = (\rho_n^R)^{-1} [ z \alpha_n (z^n e^{-in\theta}-1) \varphi_n^L(e^
{i\theta}) + z^n
e^{-in\theta}(1-ze^{-i\theta})\varphi_n^{R,*}(e^{i\theta}) ].
\end{align*}
Using the orthogonality relations \eqref{a6}, we get
\begin{align*}
\rho_n^R \int \frac{e^{i\theta}+z}{e^{i\theta}-z} \, & B_n(\theta,z)
\, d\mu(\theta) = \\
& = z^{n+1}\int ( e^{-i(n+1)\theta}\varphi_n^{R,*}(e^{i\theta}) -
\alpha_n
e^{-in\theta}\varphi_n^{L}(e^{i\theta})) \, d\mu(\theta) \\
& = z^{n+1}\rho_n^R \int e^{-i(n+1)\theta} \varphi_{n+1}^{R,*}(e^{i
\theta}) \, d\mu(\theta) = \bdzero.
\end{align*}

\smallskip
3. Relations \eqref{a4} and \eqref{a5} can be checked for $n=0$ by
a direct substitution of \eqref{firstpoly1} and
\eqref{firstpoly2}.

\smallskip
4. We obtain \eqref{a4r} and \eqref{a5r} from \eqref{a4} and \eqref
{a5} by applying the
$*$-operation.
\end{proof}

Writing the above recursion in matrix form, we get
$$
\begin{pmatrix} \psi_{n+1}^L \\ \psi_{n+1}^{R,*} \end{pmatrix}
= A^L(-\alpha_n,z) \begin{pmatrix} \psi_n^L \\ \psi_n^{R,*} \end
{pmatrix}, \qquad
\begin{pmatrix} \psi_0^L \\ \psi_0^{R,*} \end{pmatrix} =
\begin{pmatrix} \bdone \\ \bdone \end{pmatrix}
$$
for left-orthogonal polynomials and
$$
\begin{pmatrix} \psi_{n+1}^R & \psi_{n+1}^{L,*} \end{pmatrix} =
\begin{pmatrix} \psi_n^R & \psi_n^{L,*} \end{pmatrix} A^R(-\alpha_n,z),
\qquad
\begin{pmatrix} \psi_0^R & \psi_0^{L,*} \end{pmatrix} =
\begin{pmatrix} \bdone & \bdone \end{pmatrix}.
$$
for right-orthogonal polynomials.

Equivalently,
\begin{equation}
\begin{pmatrix} \psi_{n+1}^L \\ - \psi_{n+1}^{R,*} \end{pmatrix}
= A^L(\alpha_n,z) \begin{pmatrix} \psi_n^L \\ - \psi_n^{R,*} \end
{pmatrix}, \qquad
\begin{pmatrix} \psi_0^L \\ - \psi_0^{R,*} \end{pmatrix} =
\begin{pmatrix} \bdone \\ - \bdone \end{pmatrix}
\label{*2}
\end{equation}
and
$$
\begin{pmatrix} \psi_{n+1}^R & - \psi_{n+1}^{L,*} \end{pmatrix} =
\begin{pmatrix} \psi_n^R & - \psi_n^{L,*} \end{pmatrix} A^R
(\alpha_n,z), \qquad
\begin{pmatrix} \psi_0^R & - \psi_0^{L,*} \end{pmatrix} =
\begin{pmatrix} \bdone & - \bdone \end{pmatrix}.
$$

In particular, we see that the second kind polynomials $\psi_n^{L,R}$
correspond to
Verblunsky coefficients $\{ - \alpha_n \}$. We have the following
Wronskian-type relations:

\begin{proposition}
For $n \ge 0$ and $z \in \C$, we have
\begin{align}
2z^n \bdone & = \varphi_n^L(z) \psi_n^{L,*}(z) + \psi_n^L(z)
\varphi_n^{L,*}(z), \label{wronsk1} \\
2z^n \bdone & = \psi_n^{R,*}(z) \varphi_n^R(z) + \varphi_n^{R,*}(z)
\psi_n^R(z). \label{wronsk2} \\
\bdzero & = \varphi_n^L(z) \psi_n^R(z) - \psi_n^L(z) \varphi_n^R(z),
\label{wronsk3} \\
\bdzero & = \psi_n^{R,*}(z) \varphi_n^{L,*}(z) - \varphi_n^{R,*}(z)
\psi_n^{L,*}(z).
\label{wronsk4}
\end{align}
\end{proposition}

\begin{proof}
We prove this by induction. The four identities clearly hold for
$n = 0$. Suppose \eqref{wronsk1}--\eqref{wronsk4} hold for some $n
\ge 0$. Then,
\begin{align*}
\varphi_{n+1}^L \psi_{n+1}^{L,*} & + \psi_{n+1}^L \varphi_{n+1}^{L,*} =
\\
& =
\left(\rho_n^L\right)^{-1} [(z \varphi_n^L - \alpha_n^\dagger
\varphi_n^{R,*})
(\psi_n^{L,*} + z \psi_n^R \alpha_n)
\\
& \qquad + (z\psi_n^L+\alpha_n^\dagger \psi_{n}^{R,*})
(\varphi_n^{L,*} - z \varphi_n^R \alpha_n) ] \left(\rho_n^L\right)^{-1}
\\
& = \left(\rho_n^L\right)^{-1} [ z(\varphi_n^L \psi_n^{L,*} +
\psi_n^L \varphi_n^{L,*}) -
z \alpha_n^\dagger ( \psi_n^{R,*} \varphi_n^R + \varphi_n^{R,*}
\psi_n^R )\alpha_n
\\
& \qquad + \alpha_n^\dagger ( \psi_n^{R,*} \varphi_n^{L,*} -
\varphi_n^{R,*} \psi_n^{L,*}
) + z2 ( \varphi_n^L \psi_n^R - \psi_n^L \varphi_n^R ) \alpha_n ]
\left(\rho_n^L\right)^{-1}
\\
& = \left(\rho_n^L\right)^{-1} [ 2z^{n+1} (\bdone - \alpha_n^\dagger
\alpha_n) ]
\left(\rho_n^L\right)^{-1} = 2z^{n+1} \bdone,
\end{align*}
where we used \eqref{wronsk1}--\eqref{wronsk4} for $n$ in the third
step. Thus,
\eqref{wronsk1} holds for $n+1$.

For \eqref{wronsk2}, we note that
\begin{align*}
\psi_{n+1}^{R,*} \varphi_{n+1}^R & + \varphi_{n+1}^{R,*} \psi_{n+1}^R =
\\
& = (\rho_n^R)^{-1} [
(z\alpha_n \psi_{n}^L + \psi_{n}^{R,*}) (z \varphi_n^R - \varphi_n^{L,*}
\alpha_n^\dagger)
\\
& \qquad + (\varphi_n^{R,*} - z \alpha_n \varphi_n^L)(z \psi_n^R+
\psi_n^{L,*}
\alpha_n^\dagger)]
\left( \rho_n^R \right)^{-1}
\\
& = (\rho_n^R)^{-1} [ z (\psi_n^{R,*} \varphi_n^R + \varphi_n^{R,*}
\psi_n^R) - z
\alpha_n (\psi_n^L \varphi_n^{L,*} + \varphi_n^L \psi_n^{L,*})
\alpha_n^\dagger
\\
& \qquad + z2 \alpha_n (\psi_n^L \varphi_n^R - \varphi_n^L \psi_n^R)
- (\psi_n^{R,*}
\varphi_n^{L,*} - \varphi_n^{R,*} \psi_n^{L,*}) \alpha_n^\dagger ]
\left( \rho_n^R \right)^{-1}
\\
& = (\rho_n^R)^{-1} 2 z^{n+1}( \bdone - \alpha_n \alpha_n^\dagger )
\left( \rho_n^R
\right)^{-1} = 2z^{n+1} \bdone ,
\end{align*}
again using \eqref{wronsk1}--\eqref{wronsk4} for $n$ in the third step.

Next,
\begin{align*}
\varphi_{n+1}^L \psi_{n+1}^R & - \psi_{n+1}^L \varphi_{n+1}^R = \\
& = (\rho_n^L)^{-1} [ (z
\varphi_n^L - \alpha_n^\dagger \varphi_n^{R,*}) (z\psi_n^R + \psi_n^
{L,*} \alpha_n^\dagger) \\
& \qquad - (z\psi_n^L + \alpha_n^\dagger \psi_n^{R,*}) (z\varphi_n^R
- \varphi_n^{L,*}
\alpha_n^\dagger) ] (\rho_n^R)^{-1} \\
& = (\rho_n^L)^{-1} [ z2(\varphi_n^L \psi_n^R - \psi_n^L
\varphi_n^R) - \alpha_n^\dagger
(\varphi_n^{R,*} \psi_n^{L,*} - \psi_n^{R,*} \varphi_n^{L,*})
\alpha_n^\dagger \\
& \qquad - z \alpha_n^\dagger (\varphi_n^{R,*} \psi_n^R + \psi_n^
{R,*} \varphi_n^R) + z
(\varphi_n^L \psi_n^{L,*} + \psi_n^L \varphi_n^{L,*}) \alpha_n^
\dagger ] (\rho_n^R)^{-1}
\\
&= \bdzero
\end{align*}
which implies first \eqref{wronsk3} for $n+1$ and then, by applying
the $*$-operation of
order $2n+2$, also \eqref{wronsk4} for $n+1$. This concludes the
proof of the
proposition.
\end{proof}

\subsection{Christoffel--Darboux Formulas} \lb{s3.3A}

\begin{proposition} \lb{P3.6}
{\rm (a) (CD)-left orthogonal}
\begin{align*}
(1 - \bar \xi z) \sum_{k = 0}^n \varphi_k^L (\xi)^\dagger \varphi_k^L
(z)
& = \varphi_n^{R,*} (\xi)^\dagger \varphi_n^{R,*} (z) - \bar \xi z
\varphi_n^L (\xi)^\dagger \varphi_n^L (z)
\\
& = \varphi_{n+1}^{R,*} (\xi)^\dagger \varphi_{n+1}^{R,*} (z) -
\varphi_{n+1}^L (\xi)^\dagger \varphi_{n+1}^L (z).
\end{align*}
{\rm (b) (CD)-right orthogonal}
\begin{align*}
(1 - \bar \xi z) \sum_{k = 0}^n \varphi_k^R (z) \varphi_k^R (\xi)^
\dagger
& = \varphi_n^{L,*} (z) \varphi_n^{L,*} (\xi)^\dagger - \bar \xi z
\varphi_n^R (z) \varphi_n^R (\xi)^\dagger \\
& = \varphi_{n+1}^{L,*} (z) \varphi_{n+1}^{L,*} (\xi)^\dagger -
\varphi_{n+1}^R (z) \varphi_{n+1}^R (\xi)^\dagger.
\end{align*}
{\rm (c) (Mixed CD)-left orthogonal}
\begin{align*}
(1 - \bar \xi z) \sum_{k = 0}^n \psi_k^L (\xi)^\dagger
\varphi_k^L (z) & = 2\cdot \bdone - \psi_n^{R,*} (\xi)^\dagger
\varphi_n^{R,*} (z) - \bar \xi z
\psi_n^L (\xi)^\dagger \varphi_n^L (z) \\
& = 2\cdot \bdone - \psi_{n+1}^{R,*} (\xi)^\dagger \varphi_{n+1}^
{R,*} (z) -
\psi_{n+1}^L (\xi)^\dagger \varphi_{n+1}^L (z).
\end{align*}
{\rm (d) (Mixed CD)-right orthogonal}
\begin{align*}
(1 - \bar \xi z) \sum_{k = 0}^n \varphi_k^R (z) \psi_k^R (\xi)^
\dagger & =
2\cdot\bdone - \varphi_n^{L,*} (z) \psi_n^{L,*} (\xi)^\dagger - \bar
\xi z
\varphi_n^R (z) \psi_n^R (\xi)^\dagger \\
& = 2\cdot\bdone - \varphi_{n+1}^{L,*} (z) \psi_{n+1}^{L,*} (\xi)^
\dagger -
\varphi_{n+1}^R (z) \psi_{n+1}^R (\xi)^\dagger.
\end{align*}
\end{proposition}

\begin{remark} Since the $\psi$'s are themselves MOPUCs, the analogue of
(a) and (b), with all $\varphi$'s replaced by $\psi$'s, holds.
\end{remark}

\begin{proof}
(a) Write
$$
F^L_n(z) = \begin{pmatrix} \varphi_n^L(z) \\
\varphi_n^{R,*}(z)\end {pmatrix}, \quad J = \left(
\begin{array}{rr} \bdone & \bdzero \\ \bdzero & -\bdone
\end{array} \right), \quad
\tilde J  = \left( \begin{array}{cr} \bar \xi z \bdone & \bdzero \\
\bdzero & -\bdone
\end{array} \right).
$$
Then,
$$
F^L_{n+1}(z) = A^L (\alpha_n,z) F^L_n(z)
$$
and
\begin{align*}
A^L (\alpha,\xi)^\dagger & J A^L(\alpha,z) = \\
& = \begin{pmatrix} \bar \xi (\rho^L)^{-1} & - \bar \xi \alpha^
\dagger (\rho^R)^{-1} \\
-\alpha (\rho^L)^{-1} & (\rho^R)^{-1} \end{pmatrix}
\begin{pmatrix} z (\rho^L)^{-1} & - (\rho^L)^{-1} \alpha^\dagger \\
z (\rho^R)^{-1} \alpha & - (\rho^R)^{-1} \end{pmatrix} \\
& = \begin{pmatrix} \bar \xi z (\rho^L)^{-2} - \bar \xi z
\alpha^\dagger (\rho^R)^{-2} \alpha & - \bar \xi (\rho^L)^{-2}
\alpha^\dagger + \bar \xi \alpha^\dagger (\rho^R)^{-2} \\
- z \alpha (\rho^L)^{-2} + z (\rho^R)^{-2} \alpha & \alpha
(\rho^L)^{-2} \alpha^\dagger - (\rho^R)^{-2} \end{pmatrix} \\
& = \begin{pmatrix} \bar \xi z \bdone & \bdzero \\ \bdzero & - \bdone
\end{pmatrix}
= \tilde J .
\end{align*}
Thus,
\begin{align*}
F^L_{n+1}(\xi)^\dagger J F^L_{n+1}(z) &=  F^L_n(\xi)^\dagger A^L
(\alpha_n,\xi)^\dagger J A^L
(\alpha_n,z) F^L_n(z) \\
&= F^L_n(\xi)^\dagger \tilde J F^L_n(z)
\end{align*}
and hence
$$
\varphi_{n+1}^L(\xi)^\dagger \varphi_{n+1}^L(z) - \varphi_{n+1}^{R,*}
(\xi)^\dagger
\varphi_{n+1}^{R,*}(z) = \bar \xi z \varphi_n^L(\xi)^\dagger
\varphi_n^L(z) -
\varphi_n^{R,*}(\xi)^\dagger \varphi_n^{R,*}(z)
$$
which shows that the last two expressions in (a) are equal. Denote
their common value by
$Q^L_n(z,\xi)$. Then,
\begin{align*}
Q^L_n  (z,\xi) - Q^L_{n-1}(z,\xi) & =
    \varphi_n^{R,*} (\xi)^\dagger \varphi_n^{R,*} (z) - \bar \xi z
\varphi_n^L (\xi)^\dagger \varphi_n^L (z) \\
& \qquad - \varphi_n^{R,*} (\xi)^\dagger \varphi_n^{R,*} (z) +
\varphi_n^L (\xi)^\dagger \varphi_n^L (z) \\
& = (1 - \bar \xi z) \varphi_n^L (\xi)^\dagger \varphi_n^L (z).
\end{align*}
Summing over $n$ completes the proof since $Q^L_{-1}(z,\xi)=0$.

\smallskip
(b) The proof is analogous to (a): Write $F^R_n(z) = \left(
\varphi_n^R(z) \quad  \varphi_n^{L,*}(z) \right)$. Then,
$F^R_{n+1}(z) = F^R_n(z) A^R(\alpha_n,z)$ and $A^R (\alpha,z) J
A^R (\alpha,\xi)^ \dagger = \tilde J$. Thus,
\begin{align*}
F^R_{n+1}(z) J F^R_{n+1}(\xi)^\dagger & = F^R_n(z) A^R(\alpha_n,z)
J A^R (\alpha_n,\xi)^\dagger F^R_n
(\xi)^\dagger \\
&= F^R_n(z) \tilde J F^R_n(\xi)^\dagger
\end{align*}
and hence
$$
\varphi_{n+1}^R(z) \varphi_{n+1}^R(\xi)^\dagger - \varphi_{n+1}^{L,*}(z)
\varphi_{n+1}^{L,*}(\xi)^\dagger = \bar \xi z \varphi_n^R(z)
\varphi_n^R(\xi)^\dagger -
\varphi_n^{L,*}(z) \varphi_n^{L,*}(\xi)^\dagger
$$
which shows that the last two expressions in (b) are equal. Denote
their common value by $Q^R_n(z,\xi)$. Then,
$$
Q^R_n (z,\xi) - Q^R_{n-1}(z,\xi) = (1 - \bar \xi z) \varphi_n^R
(z) \varphi_n^R (\xi)^\dagger
$$
and the assertion follows as before.

\smallskip
(c) Write
$$
\tilde F^L_n(z) =  \begin{pmatrix} \psi_n^L(z) \\ -
\psi_n^{R,*}(z)
\end{pmatrix}
$$
with the second kind polynomials $\psi_n^{L,R}$. As in (a), we see
that
\begin{align*}
\tilde F^L_{n+1}(\xi)^\dagger J F^L_{n+1}(z) &=  \tilde
F^L_n(\xi)^\dagger A^L(\alpha_n,\xi)^\dagger J A^L
(\alpha_n,z) F^L_n(z) \\
& = \tilde F^L_n(\xi)^\dagger \tilde J F^L_n(z)
\end{align*}
and hence
$$
\psi_{n+1}^L(\xi)^\dagger \varphi_{n+1}^L(z) + \psi_{n+1}^{R,*}(\xi)^
\dagger
\varphi_{n+1}^{R,*}(z) = \bar \xi z \psi_n^L(\xi)^\dagger \varphi_n^L
(z) +
\psi_n^{R,*}(\xi)^\dagger \varphi_n^{R,*}(z).
$$
Denote
$$
\tilde Q^L_n(z,\xi) = 2\cdot \bdone - \psi_{n+1}^{R,*}(\xi)^\dagger
\varphi_{n+1}^{R,*}(z) -
\psi_{n+1}^L(\xi)^\dagger \varphi_{n+1}^L(z).
$$
Then,
\begin{align*}
\tilde Q^L_n  (z,\xi) - \tilde Q^L_{n-1}(z,\xi)
& = - \psi_n^{R,*} (\xi)^\dagger \varphi_n^{R,*} (z) - \bar \xi z
\psi_n^L (\xi)^\dagger \varphi_n^L (z)\\
& \qquad  + \psi_n^{R,*} (\xi)^\dagger \varphi_n^{R,*} (z) + \psi_n^L
(\xi)^\dagger \varphi_n^L (z) \\
& = (1 - \bar \xi z) \psi_n^L (\xi)^\dagger \varphi_n^L (z)
\end{align*}
and the assertion follows as before.

\smallskip
(d) Write $\tilde F^R_n(z) = \left(  \psi_n^R(z) \quad -
\psi_n^{L,*}(z) \right)$. As in (b), we see that
\begin{align*}
F^R_{n+1}(z) J \tilde F^R_{n+1}(\xi)^\dagger &= F^R_n(z)
A^R(\alpha_n,z) J A^R(\alpha_n,\xi)^\dagger \tilde
F^R_n(\xi)^\dagger \\
&= F^R_n(z) \tilde J \tilde F^R_n(\xi)^\dagger
\end{align*}
and hence
$$
\varphi_{n+1}^R(z) \psi_{n+1}^R(\xi)^\dagger + \varphi_{n+1}^{L,*}(z)
\psi_{n+1}^{L,*}(\xi)^\dagger = \bar \xi z \varphi_n^R(z) \psi_n^R
(\xi)^\dagger +
\varphi_n^{L,*}(z) \psi_n^{L,*}(\xi)^\dagger.
$$
With $\tilde Q^R_n(z,\xi) = 2\cdot\bdone - \varphi_{n+1}^{L,*}(z)
\psi_{n +1}^{L,*}(\xi)^\dagger - \varphi_{n+1}^R(z)
\psi_{n+1}^R(\xi)^\dagger$, we have
$$
\tilde Q^R_n (z,\xi) - \tilde Q^R_{n-1}(z,\xi) = (1 - \bar \xi z)
\varphi_n^R (z)
\psi_n^R (\xi)^\dagger
$$
and we conclude as in (c).
\end{proof}

\subsection{Zeros of MOPUC} \lb{s3.3B}

Our main result in this section is:

\begin{theorem}\lb{T3.6A} All the zeros of $\det(\varphi_n^R(z))$ lie
in $\bbD=
\{z\, \colon \abs{z}<1\}$.
\end{theorem}

We will also prove:

\begin{theorem}\lb{T3.6B} For each $n$,
\begin{equation}\lb{3.24a}
\det(\varphi_n^R(z))=\det(\varphi_n^L(z)).
\end{equation}
\end{theorem}

The scalar analogue of Theorem~\ref{T3.6A} has seven proofs in \cite
{S}! The simplest is due to
Landau \cite{Landau} and its MOPUC analogue is Theorem~2.13.7 of \cite
{S}. There is also a proof
in Delsarte et al.\ \cite{DGK} who attribute the theorem to Whittle
\cite{Whit}.  We provide two
more proofs here, not only for their intrinsic interest: our first
proof we need because it
depends only on the recursion relation (it is related to the proof of
Delsarte et al.\ \cite{DGK}).
The second proof is here since it relates zeros to eigenvalues of a
cutoff CMV matrix.

\begin{theorem}\lb{T3.6C} We have
\begin{SL}
\item[{\rm{(i)}}] For $z\in\partial\bbD$, all of $\varphi_n^{R,*}(z)
$, $\varphi_n^{L,*}(z)$, $\varphi_n^R(z)$,
$\varphi_n^L(z)$ are invertible.

\item[{\rm{(ii)}}] For $z\in\partial\bbD$, $\varphi_n^L(z) (\varphi_n^
{R,*}(z))^{-1}$ and $(\varphi_n^{*,L}(z))^{-1}
\varphi_n^R(z)$ are unitary.

\item[{\rm{(iii)}}] For $z\in\bbD$, $\varphi_n^{R,*}(z)$ and
$\varphi_n^{L,*}(z)$ are invertible.

\item[{\rm{(iv)}}] For $z\in\bbD$, $\varphi_n^L(z)(\varphi_n^{R,*}(z))
^{-1}$ and $(\varphi_n^{*,L}(z))^{-1}
\varphi_n^R(z)$ are of norm at most $1$ and, for $n\geq 1$, strictly
less than $1$.

\item[{\rm{(v)}}] All zeros of $\det(\varphi_n^{R,*}(z))$ and $\det
(\varphi_n^{L,*}(z))$ lie in $\bbC\setminus
\ol{\bbD}$.

\item[{\rm{(vi)}}] All zeros of $\det(\varphi_n^R(z))$ and $\det
(\varphi_n^L(z))$ lie in $\bbD$.
\end{SL}
\end{theorem}

\begin{remark} (vi) is our first proof of Theorem~\ref{T3.6A}.
\end{remark}

\begin{proof} All these results are trivial for $n=0$, so we can hope
to use an inductive argument. So  suppose we have the result for
$n -1$.

By \eqref{3.7b},
\begin{equation}\lb{3.24b}
\varphi_n^{R,*} = (\rho_{n-1}^R)^{-1} (\bdone
-z\alpha_{n-1}\varphi_{n-1}^L (\varphi_{n-1}^{R,*})^{-1})
\varphi_{n-1}^{R,*}.
\end{equation}
Since $\abs{\alpha_{n-1}}<1$, if $\abs{z}\leq 1$, each factor on the
right of \eqref{3.24b} is invertible.
This proves (i) and (iii) for $\varphi_n^{R,*}$ and a similar
argument works for $\varphi_n^{L,*}$. If
$z=e^{i\theta}$, $\varphi_n^R (e^{i\theta})=e^{in\theta}\varphi_n^
{R,*} (e^{i\theta})^\dagger$ is also
invertible, so we have (i) and (iii) for $n$.

Next, we claim that if $z\in\partial\bbD$, then
\begin{equation}\lb{3.24c}
\varphi_n^{R,*}(z)^\dagger \varphi_n^{R,*}(z) =\varphi_n^L(z)^\dagger
\varphi_n^L(z).
\end{equation}
This follows from taking $z=\xi\in\partial\bbD$ in Proposition~\ref
{P3.6}(a). Given that
$\varphi_n^{R,*}(z)$ is invertible, this implies
\begin{equation}\lb{3.24d}
\bdone =(\varphi_n^L(z) \varphi_n^{R,*}(z)^{-1})^\dagger (\varphi_n^L
(z) \varphi_n^{R,*}(z)^{-1})
\end{equation}
proving the first part of (ii) for $n$. The second part of (ii) is
proven similarly by using Proposition~\ref{P3.6}(b).

For $z\in\ol{\bbD}$, let
\[
F(z)=\varphi_n^L(z) \varphi_n^{R,*}(z)^{-1}.
\]
Then $F$ is analytic in $\bbD$, continuous in $\ol{\bbD}$, and $\norm
{F(z)}=1$ on $\partial\bbD$,
so (iv) follows from the maximum principle.

Since $\varphi_n^{R,*}(z)$ is invertible on $\ol{\bbD}$, its $\det$
is non-zero there, proving (v).
(vi) then follows from
\begin{equation}\lb{3.24e}
\det(\varphi_n^R(z)) =z^{nl}\, \ol{\det(\varphi_n^{R,*}(1/\bar z))}\, .
\end{equation}
\end{proof}

Let $\calV$ be the $\bbC^l$-valued functions on $\partial\bbD$ and
$\calV_n$ the span of the $\bbC^l$-valued polynomials of degree at
most $n$, so
\[
\dim(\calV_n)=\bbC^{l(n+1)}.
\]
Let $\calV_\infty$ be the set $\cup_n \calV_n$ of all
$\bbC^l$-valued polynomials. Let $\pi_n$ be the projection onto
$\calV_n$ in the $\calV$ inner product \eqref {1.33}.

It is easy to see that
\begin{equation}\lb{3.24f}
\calV_n \cap \calV_{n-1}^\perp = \{\Phi_n^R(z) v\, \colon v\in\bbC^l\}
\end{equation}
since $\jap{z^l, \Phi_n^R(z) v} =0$ for $l=0,\dots, n-1$ and the
dimensions  on the left and right of \eqref{3.24f} coincide.
$\calV_n \cap \calV_{n-1}^\perp $ can also be described as the set
of $(v^\dagger \Phi_n^L(z))^\dagger$ for $v\in\bbC^l$.

We define $M_z \colon \calV_{n-1} \to \calV_n$ or $\calV_\infty
\to \calV_\infty$ as the operator of multiplication by $z$.

\begin{theorem}\lb{T3.6D} For all $n$, we have
\begin{equation}\lb{3.24g}
\det_{\bbC^l} (\Phi_n^R(z)) =\det_{\calV_{n-1}} (z\bdone-\pi_{n-1}
M_z \pi_{n-1}).
\end{equation}
\end{theorem}

\begin{remarks} 1. Since $\norm{M_z}\leq 1$, \eqref{3.24g}
immediately implies zeros of $\det(\varphi_n^R(z))$
lie in $\ol{\bbD}$, and a small additional argument proves they lie
in $\bbD$. As we will see, this also
implies Theorem~\ref{T3.6B}.

\smallskip
2. Of course, $\pi_{n-1} M_z \pi_{n-1}$ is a cutoff CMV matrix if
written in a CMV basis.
\end{remarks}

\begin{proof} If $Q\in \calV_{n-k}$, then by \eqref{3.24f},
\begin{equation}\lb{3.24h}
\pi_{n-1}[(z-z_0)^k Q] = 0 \Leftrightarrow (z-z_0)^k Q=\Phi_n^R(z) v
\end{equation}
for some $v\in\bbC^l$. Thus writing $\det(\Phi_n^R(z)) =\Phi_n^R(z)
v_1 \wedge \dots \wedge
\Phi_n^R(z) v_l$ in a Jordan basis for $\Phi_n^R(z)$, we see that the
order of the zeros of
$\det(\Phi_n^R(z))$ at $z_0$ is exactly the order of $z_0$ as an
algebraic eigenvalue of
$\pi_{n-1} M_z \pi_{n-1}$, that is, the order of $z_0$ as a zero of
the right side of \eqref{3.24g}.

Since both sides of \eqref{3.24g} are monic polynomials of degree $nl
$ and their zeros including
multiplicity are the same, we have proven \eqref{3.24g}.
\end{proof}

\begin{proof}[Proof of Theorem~\ref{T3.6B}] On the right side of
\eqref{3.24h}, we can put $(\Phi_n^L(z) v^\dagger)^\dagger$ and so
conclude \eqref{3.24g} holds with $\Phi_n^L(z)$ on the left.

This proves \eqref{3.24a} if $\varphi$ is replaced by $\Phi$. Since
$\alpha_j^*\alpha_j$ and $\alpha_j\alpha_j^*$ are unitarily equivalent,
$\det(\rho_j^L)=\det (\rho_j^R)$. Thus, $\det(\kappa_n^L) = \det
(\kappa_n^R)$,
and we obtain \eqref{3.24a} for $\varphi$.
\end{proof}

It is a basic fact (Theorem~1.7.5 of \cite{S}) that for the scalar
case, any set of $n$ zeros in $\bbD$ are the zeros of a unique
OPUC $\Phi_n$ and any monic polynomial with all its zeros in
$\bbD$ is a monic OPUC. It is easy to see that any set of $nl$
zeros in $\bbD$ is the set of zeros of an OPUC $\Phi_n$, but
certainly not unique. It is an interesting open question to
clarify what matrix monic OPs are monic MOPUCs.

\subsection{Bernstein--Szeg\H{o} Approximation} \lb{s3.5A}

Given $\{\alpha_j\}_{j=0}^{n-1}\in\bbD^n$, we use Szeg\H{o} recursion
to define polynomials
$\varphi_j^R, \varphi_j^L$ for $j=0,1,\dots, n$. We define a measure
$d\mu_n$ on $\partial\bbD$
by
\begin{equation}\lb{3.32a}
d\mu_n(\theta) = [\varphi_n^R(e^{i\theta})\varphi_n^R (e^{i\theta})^
\dagger]^{-1}
\, \f{d\theta}{2\pi}.
\end{equation}
Notice that \eqref{3.24c} can be rewritten
\begin{equation}\lb{3.32b}
\varphi_n^R(e^{i\theta}) \varphi_n^R(e^{i\theta})^\dagger =
\varphi_n^L (e^{i\theta})^\dagger
\varphi_n^L (e^{i\theta}).
\end{equation}

We use here and below the fact that the proof of
Theorem~\ref{T3.6C} only depends on Szeg\H{o} recursion and not on
the a priori existence of a measure. That theorem also shows the
inverse in \eqref{3.32a} exists. Thus,
\begin{equation}\lb{3.32c}
d\mu_n(\theta)= [\varphi_n^L (e^{i\theta})^\dagger \varphi_n^L (e^{i
\theta})]^{-1} \, \f{d\theta}{2\pi}\, .
\end{equation}

\begin{theorem} \lb{T3.10A} The measure $d\mu_n$ is normalized
{\rm{(}}i.e., $\mu_n(\partial \bbD)=\bdone${\rm{)}}
and its right MOPUC for $j=0,\dots, n$ are
$\{\varphi_j^R\}_{j=0}^n$, and for $j>n$,
\begin{equation}\lb{3.32d}
\varphi_j^R(z) = z^{j-n} \varphi_n^R(z).
\end{equation}
The Verblunsky coefficients for $d\mu_n$ are
\begin{equation}\lb{3.32e}
\alpha_j (d\mu_n) = \begin{cases} \alpha_j, & j\leq n, \\
\bdzero, & j\geq n+1.
\end{cases}
\end{equation}
\end{theorem}

\begin{remarks} 1. In the scalar case, one can multiply by a constant
and renormalize, and then
prove the constant is $1$. Because of commutativity issues, we need a
different argument here.

\smallskip
2. Of course, using \eqref{3.32c}, $\varphi_n^L$ are left MOPUC for $d
\mu_n$.

\smallskip
3. Our proof owes something to the scalar proof in \cite{ENZG91}.
\end{remarks}

\begin{proof}
Let $\ang{\cdot,\cdot}_R$ be the inner product associated with $
\mu_n$. By a direct computation, $\ang{\varphi_n^R, \varphi_n^R}_R
=\bdone$, and for $j=0,1,\dots, n-1$,
\begin{align*}
\ang{z^j, \varphi_n^R}_R &= \f{1}{2\pi} \int_0^{2\pi} e^{-ij\theta}
(\varphi_n^R (e^{i\theta})^\dagger)^{-1}\, d\theta \\
&= \f{1}{2\pi i} \oint z^{n-j-1} (\varphi_n^{R,*}(z))^{-1} \, dz =
\bdzero
\end{align*}
by analyticity of $\varphi_n^{R,*}(z)^{-1}$ in $\bbD$ (continuity
in $\ol{\bbD}$).

This proves $\varphi_n^R$ is a MOPUC for $d\mu_n$ (and a similar
calculation works for the right side of \eqref{3.32d} if $j\geq
n$). By the inverse Szeg\H{o} recursion and induction downwards,
$\{\varphi_j^R\}_{j=0}^{n-1}$ are also OPs, and by the Szeg\H{o}
recursion, they are normalized. In particular, since
$\varphi_0^R\equiv\bdone$ is normalized, $d\mu_n$ is normalized.
\end{proof}

\subsection{Verblunsky's Theorem}\lb{s3.5B}

We can now prove the analogue of Favard's theorem for MOPUC; the
scalar version is called Verblunsky's
theorem in \cite{S} after \cite{V35}. A history and other proofs can
be found in \cite{S}. The proof
below is essentially the matrix analogue of that of Geronimus \cite
{Ger46}
(rediscovered in \cite{DGK,ENZG91}).
Delsarte et al. \cite{DGK} presented their proof in the MOPUC case
and they seem to have been the first
with a matrix Verblunsky theorem. One can extend the CMV and the
Geronimus theorem proofs from the scalar
case to get alternate proofs of the theorem below.

\begin{theorem}[Verblunsky's Theorem for MOPUC]\lb{T3.10B}
Any sequence $\{\alpha_j\}_{j=0}^\infty \in\bbD^\infty$ is the sequence of
Verblunsky coefficients of a unique measure.
\end{theorem}

\begin{proof} Uniqueness is easy, since the $\alpha$'s determine the $
\varphi_j^R$'s and so the
$\Phi_j^R$'s which determine the moments.

Given a sequence $\{\alpha_j\}_{j=0}^\infty$, let $d\mu_n$ be the
measures of the last section. By
compactness of $l\times l$ matrix-valued probability measures on
$\partial\bbD$, they have a weak limit. By
using limits, $\{\varphi_j^R\}_{j=0}^\infty$ are the right MOPUC for
$d\mu$ and they determine
the proper Verblunsky coefficients.
\end{proof}

\subsection{Matrix POPUC}\lb{s3.5C}

Analogously to the scalar case (see \cite
{CMV02,Gol02,JNT89,S308,Wong}), given any unitary $\beta$
in $\calM_l$, we define
\begin{equation}\lb{3.32f}
\varphi_n^R(z;\beta) = z\varphi_{n-1}^R(z) -\varphi_{n-1}^{L,*}(z)
\beta^\dagger.
\end{equation}
As in the scalar case, this is related to the secular determinant of
unitary extensions of the
cutoff CMV matrix. Moreover,

\begin{theorem}\lb{T3.10C} Fix $\beta$. All the zeros of $(\varphi_n
(z;\beta))$ lie on $\partial\bbD$.
\end{theorem}

\begin{proof} If $\abs{z}<1$, $\varphi_{n-1}^{L,*}(z)$ is invertible and
\[
\varphi_n^R(z;\beta) = -\varphi_{n-1}^{L,*}(z) \beta^\dagger (\bdone -
z\beta \,
\varphi_{n-1}^{L,*}(z)^{-1} \varphi_{n-1}^R(z))
\]
is invertible since the last factor differs from $\bdone$ by a strict
contraction. A similar argument
shows invertibility if $\abs{z}>1$. Thus, the only zeros of $\det
(\dott)$ lie in $\partial\bbD$.
\end{proof}

\subsection{Matrix-Valued Carath\'eodory and Schur Functions}

An analytic matrix-valued function $F$ defined on $\D$ is called a
(matrix-valued)
Carath\'eodory function if $F(0) = \bdone$ and $\mathrm{Re} \, F(z)
\equiv \f12 (F(z)+F(z)^\dagger)
\ge 0$ for every $z \in\D$. The following result can be found in \cite[Thm.~2.2.2]{DFK}.

\begin{theorem}[Riesz--Herglotz]
If $F$ is a matrix-valued Carath\'eodory function, then there
exists a unique positive semi-definite matrix measure $d\mu$ such
that
\begin{equation}\label{rep}
F(z) = \int \frac{e^{i\theta} + z}{e^{i\theta} - z} \, d\mu(\theta).
\end{equation}
The measure $d\mu$ is given by the unique weak limit of the measures
$d\mu_r(\theta) =
\mathrm{Re} \, F(re^{i\theta}) \f{d\theta}{2\pi}$ as $r \uparrow 1$.
Moreover,
$$
F(z) = c_0 + 2 \sum_{n=1}^\infty c_n z^n
$$
where
$$
c_n = \int e^{-in\theta} \, d\mu(\theta).
$$
Conversely, if $d\mu$ is a positive semi-definite matrix measure,
then \eqref{rep} defines a matrix-valued Carath\'eodory function.
\end{theorem}

An analytic matrix-valued function $f$ defined on $\D$ is called a
(matrix-valued) Schur
function if $f(z)^\dagger f(z) \le \bdone$ for every $z \in \D$. This
condition is equivalent
to $f(z) f(z)^\dagger \le \bdone$ for every $z \in \D$ and to $\norm{f
(z)}\leq 1$ for every $z\in
\bbD$. By the maximum principle, if $f$ is not constant, the
inequalities are strict.
The following can be found in \cite[Prop.~4.5.3]{S}:

\begin{proposition}
The association
\begin{align}
f(z) & = z^{-1} (F(z) - \bdone)(F(z) + \bdone)^{-1}, \label{Ff}\\
F(z) & = (\bdone + zf(z))(\bdone - zf(z))^{-1} \label{fF}
\end{align}
sets up a one-one correspondence between matrix-valued Carath\'eodory
functions and
matrix-valued Schur functions.
\end{proposition}

\begin{proposition}\label{prop.3.16}
For $z \in \D$, we have
\begin{equation}\label{ReFfromf}
\mathrm{Re} \, F(z) = (\bdone - \bar z f(z)^\dagger)^{-1} (\bdone - |
z|^2 f(z)^\dagger f(z))
(\bdone - z f(z))^{-1}
\end{equation}
and the non-tangential boundary values $\mathrm{Re} \, F(e^{i\theta})
$ and
$f(e^{i\theta})$ exist for Lebesgue almost every $\theta$.

Write $d\mu(\theta) = w(\theta) \frac{d\theta}{2\pi} + d\mu_\s$.
Then, for
almost every $\theta$,
\begin{equation}\label{wandReF}
w(\theta) = \mathrm{Re} \, F(e^{i\theta})
\end{equation}
and for a.e.\ $\theta$, $\det (w(\theta)) \not= 0$ if and only if $f
(e^{i\theta})^\dagger
f(e^{i\theta}) < \bdone$.
\end{proposition}

\begin{proof}
The identity \eqref{ReFfromf} follows from \eqref{fF}. The
existence of the boundary values of $f$ follows by application of
the scalar result to the individual entries of $f$. Then
\eqref{fF} gives the boundary values of $F$. We also used the
following fact: Away from a set of zero Lebesgue measure, $\det
(\bdone - z f(z))$ has non-zero boundary values by general
$H^\infty$ theory.

\eqref{wandReF} holds for $\jap{\eta, F(z)\eta}_{\C^l}$ and $\jap
{\eta,d\mu\eta}_{\C^l}$ for any $\eta\in\bbC^l$ by the scalar
result. We get \eqref{wandReF} by polarization. From
$$
w(\theta) = (\bdone - e^{-i\theta} f(e^{i\theta})^\dagger)^{-1}
(\bdone - f(e^{i\theta})^\dagger
f(e^{i\theta})) (\bdone - e^{i\theta} f(e^{i\theta}))^{-1}
$$
it follows immediately that $f(e^{i\theta})^\dagger f(e^{i\theta}) <
\bdone$ implies $\det
(w(\theta)) > 0$. Conversely, if $f(e^{i\theta})^\dagger f(e^{i
\theta}) \le \bdone$ but not
$f(e^{i\theta})^\dagger f(e^{i\theta}) < 1$, then $\det (\bdone - f(e^
{i\theta})^\dagger
f(e^{i\theta})) = 0$ and by our earlier arguments $\det (\bdone - e^{-
i\theta}
f(e^{i\theta})^\dagger)^{-1}$ and $\det (\bdone - e^{i\theta} f(e^{i
\theta}))^{-1}$ exist and
are finite; hence $\det (w(\theta)) = 0$. All previous statements are
true away from
suitable sets of zero Lebesgue measure.
\end{proof}

\subsection{Coefficient Stripping, the Schur Algorithm, and
Geronimus' Theorem} \lb{s3.5} The matrix version of Geronimus'
theorem goes back at least to the book of Bakonyi--Constantinescu
\cite{BakCon}. Let $F(z)$ be the matrix-valued Carath\'eodory
function \eqref{rep} (with the same measure $\mu$ as the one used
in the definition of $ \ang{\cdot,\cdot}_R$). Let us denote
\begin{align*}
u_n^L(z) & = \psi_n^L(z) + \varphi_n^L(z) F(z), \\
u_n^R(z) & = \psi_n^R(z) + F(z) \varphi_n^R(z).
\end{align*}
We also define
\begin{align*}
u_n^{L,*}(z) & = \psi_n^{L,*}(z) - F(z) \varphi_n^{L,*}(z), \\
u_n^{R,*}(z) & = \psi_n^{R,*}(z) - \varphi_n^{R,*}(z) F(z).
\end{align*}

\begin{proposition}\label{prp8}
For any $\abs{z}<1$, the sequences $u_n^L(z)$, $u_n^R(z)$, $u_n^{L,*}
(z)$, $u_n^{R,*}(z)$
are square summable.
\end{proposition}
\begin{proof}
Denote
$$
f(\theta)= \frac{e^{-i\theta}+\bar z}{e^{-i\theta}-\bar z}\, , \qquad
g(\theta)=\frac{e^{i\theta}+z}{e^{i\theta}-z}\, .
$$
By the definitions \eqref{a1}--\eqref{a3}, we have
\begin{align*}
u_n^L(z)&= \ang{f, \varphi_n^L}_L,
\\
-u_n^{R,*}(z)&= z^n \ang{\varphi_n^R,g}_R,
\\
-u_n^{L,*}(z)&= z^n \ang{\varphi_n^L,g}_L,
\\
u_n^R(z)&= \ang{f, \varphi_n^R}_R.
\end{align*}
Using the Bessel inequality and the fact that $\abs{z}<1$, we obtain
the required
statements.
\end{proof}

Next we will consider sequences defined by
\begin{equation}\label{cocy}
\begin{pmatrix} s_n \\ t_n \end{pmatrix} = A^L(\alpha_{n-1},z) \cdots
A^L(\alpha_0,z) \begin{pmatrix} s_0 \\ t_0 \end{pmatrix}
\end{equation}
where $s_n, t_n\in\calM_l$. Similarly, we will consider the sequences
\begin{equation}\label{cocyr}
(s_n,  t_n ) = (s_0, t_0) A^R(\alpha_{0},z) \cdots A^R(\alpha_{n-1},z)
\end{equation}

\begin{theorem}\label{th.solu}
Let $z \in \D$ and let  $f$ be the Schur function associated with $d
\mu$ via \eqref{rep}
and \eqref{Ff}. Then:
\begin{SL}
\item[{\rm{(i)}}] A solution of \eqref{cocy} is square summable if
and only if the initial condition is
of the form
$$
\begin{pmatrix} s_0 \\ t_0 \end{pmatrix} = \begin{pmatrix} c \\
z f(z) c \end{pmatrix}
$$
for some matrix $c$ in $\calM_l$.

\item[{\rm{(ii)}}] A solution of \eqref{cocyr} is square summable if
and only if the initial condition
is of the form
$$
(  s_0, t_0)  = (c, c z f(z))
$$
for some matrix $c$.
\end{SL}
\end{theorem}

\begin{proof}
We shall prove (i); the proof of (ii) is similar.

1. By Proposition~\ref{prp8} and \eqref{*1}, \eqref{*2}, we have
the square summable solution
\begin{equation}
\begin{pmatrix} s_n \\ t_n \end{pmatrix} = \begin{pmatrix} u_n^L(z) \
\ -u_n^{R,*}(z)\end{pmatrix} , \quad
\begin{pmatrix} s_0 \\ t_0 \end{pmatrix} = \begin{pmatrix} d \\ zf(z)
d \end{pmatrix}, \quad
d =(F(z)+\bdone). \label{a9}
\end{equation}
The matrix $d=F(z)+\bdone$ is invertible. Thus, multiplying the above
solution on the right by
$d^{-1}c$ for any given matrix $c$, we get the ``if" part of the
theorem.

\smallskip
2. Let us check that $\varphi_n^{R,*}c$ is not square summable for
any matrix $c\not=0$.
By the CD formula with $\xi=z$, we have
$$
(1-\abs{z}^2) \sum_{k=1}^{n-1} \varphi_k^L(z)^\dagger \varphi_k^L(z) +
\varphi_n^L(z)^\dagger \varphi_n^L(z) = \varphi_n^{R,*}(z)^\dagger
\varphi_n^{R,*}(z)
$$
and so
$$
\varphi_n^{R,*}(z)^\dagger \varphi_n^{R,*}(z) \geq (1-\abs{z}^2)
\varphi_0^L(z)^\dagger \varphi_0^L(z) = (1-\abs{z}^2) \bdone .
$$
Thus, we get
$$
\norm{\varphi_n^{R,*}c}_R2 \ge (1-\abs{z}^2)\Tr c^\dagger c >0.
$$

\smallskip
3. Let $\left(\begin{smallmatrix} s_n  \\ t_n\end{smallmatrix}\right)
$ be any square summable
solution to \eqref{cocy}. Let us write this solution as
\begin{equation}
\begin{pmatrix}s_n \\ t_n \end{pmatrix} = \begin{pmatrix} \varphi_n^L
a \\ \varphi_n^{R,*}a \end{pmatrix}
+ \begin{pmatrix} \psi_n^L b \\ -\psi_n^{R,*}b \end{pmatrix},
\quad a= \frac{s_0+t_0}{2}\, , \quad b= \frac{s_0-t_0}{2}\, .
\label{a10}
\end{equation}
Multiplying the solution \eqref{a9} by $b$ and subtracting from \eqref
{a10}, we get a
square summable solution
$$
\begin{pmatrix} \varphi_n^L (z)(a-F(z)b) \\ \varphi_n^{R,*}(a-F(z)b)
\end{pmatrix}.
$$
It follows that $a=F(z)b$, which proves the ``only if" part.
\end{proof}

The Schur  function $f$ is naturally associated with $d\mu$ and
hence with the Verblunsky coefficients $\alpha_0, \alpha_1,
\alpha_2, \dots$. The Schur functions obtained by coefficient
stripping will be denoted by $f_1, f_2, f_3, \dots$, that is,
$f_n$ corresponds to Verblunsky coefficients $\alpha_n,
\alpha_{n+1}, \alpha_{n+2}, \dots$. We also write $f_0 \equiv f$.

\begin{theorem}[Schur Algorithm and Geronimus' Theorem]\lb{T3.12}
For the Schur functions $f_0, f_1, f_2, \ldots$ associated with
Verblunsky coefficients $\alpha_0, \alpha_1, \alpha_2, \dots$, the
following relations hold:
\begin{align}
f_{n+1}(z) & = z^{-1} (\rho_n^R)^{-1} [f_n(z) - \alpha_n]\,
[\bdone - \alpha_n^\dagger f_n(z)]^{-1} \rho_n^L, \label{salg}\\
f_n(z) & = (\rho_n^R)^{-1} [z f_{n+1}(z) + \alpha_n] \,[\bdone +
z \alpha_n^\dagger f_{n+1}(z)]^{-1} \rho_n^L. \label{salginv}
\end{align}
\end{theorem}

\begin{remarks} 1. See \eqref{1.50} and Theorem~\ref{T1.4} to understand
this result.

\smallskip
2. \eqref{1.82a} provides an alternate way to write \eqref{salg} and \eqref{salginv}.
\end{remarks}

\begin{proof}
It clearly suffices to consider the case $n = 0$. Consider the
solution of \eqref{cocy}
with initial condition
$$
\begin{pmatrix} \bdone \\ z f_0(z)\end{pmatrix}.
$$
By Theorem~\ref{th.solu}, there exists a matrix $c$ such that
\begin{align*}
\begin{pmatrix} c \\ z f_1(z) c \end{pmatrix} & = A^L(\alpha_0,z)
\begin{pmatrix} \bdone \\ z f_0(z) \end{pmatrix} \\
& = \begin{pmatrix} z (\rho_0^L)^{-1} - z (\rho_0^L)^{-1}
\alpha_0^\dagger f_0(z) \\ - z (\rho_0^R)^{-1} \alpha_0 + z (\rho_0^R)
^{-1} f_0(z)
\end{pmatrix}.
\end{align*}
{}From this we can compute $zf_1(z)$:
\begin{align*}
z f_1(z) & = [- z (\rho_0^R)^{-1} \alpha_0 + z (\rho_0^R)^{-1} f_0(z)]
\, [z (\rho_0^L)^{-1} - z
(\rho_0^L)^{-1} \alpha_0^\dagger f_0(z)]^{-1} \\
& = (\rho_0^R)^{-1} [f_0(z) - \alpha_0]\, [\bdone -\alpha_0^\dagger
f_0(z)]^{-1} \rho_0^L
\end{align*}
which is \eqref{salg}.

Similarly, we can express $f_0$ in terms of $f_1$. From
\begin{align*}
\begin{pmatrix} \bdone \\ z f_0(z) \end{pmatrix}
& = A^L(\alpha_0,z)^{-1} \begin{pmatrix} c \\ z f_1(z) c \end
{pmatrix}  \\
& = \begin{pmatrix} z^{-1} (\rho_0^L)^{-1} & z^{-1} (\rho_0^L)^{-1}
\alpha_0^\dagger \\
(\rho_0^R)^{-1} \alpha_0 & (\rho_0^R)^{-1} \end{pmatrix}
\begin{pmatrix} c \\ z f_1(z) c \end{pmatrix}  \\
& = \begin{pmatrix} z^{-1} (\rho_0^L)^{-1} c +
(\rho_0^L)^{-1} \alpha_0^\dagger f_1(z) c \\ (\rho_0^R)^{-1}
\alpha_0 c + (\rho_0^R)^{-1} z f_1(z) c \end{pmatrix}
\end{align*}
we find that
\begin{align*}
z f_0(z) & = [(\rho_0^R)^{-1} \alpha_0 c + (\rho_0^R)^{-1}
z f_1(z) c] \,[z^{-1} (\rho_0^L)^{-1} c + (\rho_0^L)^{-1} \alpha_0^
\dagger f_1(z) c ]^{-1} \\
& = [(\rho_0^R)^{-1} \alpha_0 + (\rho_0^R)^{-1} z f_1(z)]\, [z^{-1}
(\rho_0^L)^{-1} +
(\rho_0^L)^{-1} \alpha_0^\dagger f_1(z)]^{-1} \\
& = (\rho_0^R)^{-1} [\alpha_0 + z f_1(z)]\, [z^{-1} \bdone + \alpha_0^
\dagger f_1(z)]^{-1} \rho_0^L
\end{align*}
which gives \eqref{salginv}.
\end{proof}

\subsection{The CMV Matrix}

In this section and the next, we discuss CMV matrices for MOPUC.
This was discussed first by Simon in \cite{SimonRev}, which also
has the involved history in the scalar case. Most of the results
in this section appear already in \cite{SimonRev}; the results of
the next section are new here---they parallel the discussion in
\cite[Sect.~4.4]{S} where these results first appeared in the
scalar case.

\subsubsection{The CMV basis}

Consider the two sequences $\chi_n$, $x_n\in \h$, defined by
\begin{alignat*}{2}
\chi_{2k}(z) &= z^{-k} \varphi_{2k}^{L,*}(z) , \qquad & \chi_{2k-1}
(z) &= z^{-k+1}
\varphi_{2k-1}^R(z),
\\
x_{2k}(z) &= z^{-k} \varphi_{2k}^R(z) , \qquad & x_{2k-1}(z) &= z^{-k}
\varphi_{2k-1}^{L,*}(z).
\end{alignat*}
For an integer $k\geq0$, let us introduce the following notation:
$i_k$ is the $(k+1)$th term of the sequence
$0,1,-1,2,-2,3,-3,\dots$, and $j_k$ is the $(k+1)$th term of the
sequence $0,-1,1,-2,2,-3,3,\dots$. Thus, for example, $i_1=1$,
$j_1=-1$.

We use the right module structure of $\h$. For a set of functions
$\{f_k(z)\}_{k=0}^n\subset\h$, its module span is the set of all
sums $\sum f_k(z)a_k$ with $a_k\in \calM_l$.

\begin{proposition}
\begin{SL}
\item[{\rm{(i)}}] For any $n\geq 1$, the module span of $\{\chi_k\}_
{k=0}^n$ coincides with the module
span of $\{z^{i_k}\}_{k=0}^n$ and the module span of $\{x_k\}_{k=0}^n
$ coincides with the
module span of $\{z^{j_k}\}_{k=0}^n$.

\item[{\rm{(ii)}}] The sequences $\{\chi_k\}_{k=0}^\infty$ and $\{x_k
\}_{k=0}^\infty$ are orthonormal:
\begin{equation}
\ang{\chi_k, \chi_m}_R=\ang{x_k, x_m}_R=\delta_{km}. \label{cmv1}
\end{equation}
\end{SL}
\end{proposition}

\begin{proof}
(i) Recall that
\begin{align*}
\varphi_n^R(z)&= \kappa_n^R z^n+\lc\{\bdone,\dots, z^{n-1}\},
\\
\varphi_n^{L,*}(z)&= (\kappa_n^L)^\dagger+\lc\{z,\dots, z^{n}\},
\end{align*}
where both $ \kappa_n^R$ and $(\kappa_n^L)^\dagger$ are invertible
matrices. It follows
that
\begin{align*}
\chi_n(z)&=\gamma_n z^{i_n}+\lc\{z^{i_0},\dots, z^{i_{n-1}}\},
\\
x_n(z)&=\delta_n z^{j_n}+\lc\{z^{j_0},\dots, z^{j_{n-1}}\},
\end{align*}
where  $\gamma_n$, $\delta_n$ are invertible matrices. This proves (i).

\smallskip
(ii) By the definition of $\varphi_n^L$ and $\varphi_n^R$, we have
\begin{gather}
\ang{\varphi_n^R, \varphi_m^R}_R=\ang{\varphi_n^{L,*}, \varphi_m^
{L,*}}_R =\delta_{nm},
\label{cmv2}
\\
\ang{\varphi_n^R,z^m}_R=\bdzero, \; m=0,\dots,n-1; \quad \ang
{\varphi_n^{L,*},z^m}_R=\bdzero, \;
m=1,\dots,n. \label{cmv3}
\end{gather}
{}From \eqref{cmv2} with $n=m$, we get
$$
\ang{\chi_n,\chi_n}_R=\ang{x_n,x_n}_R=\bdone.
$$
Considering separately the cases of even and odd $n$, it is easy
to prove that
\begin{alignat}{2}
\ang{\chi_n, z^m}_R &=\bdzero, \quad &  m &= i_0,i_1,\dots,i_{n-1},
\label{cmv4}
\\
\ang{x_n, z^m}_R&=\bdzero, \quad & m &= j_0,j_1,\dots,j_{n-1}. \label
{cmv5}
\end{alignat}
For example, for $n=2k$, $m\in\{i_0,\dots, i_{2k-1}\}$ we have
$m+k\in\{1,2,\dots, 2k\}$ and so, by \eqref{cmv3},
$$
\ang{\chi_n,z^m}_R = \ang{z^{-k}\varphi_{2k}^{L,*},z^m}_R =
\ang{\varphi_{2k}^{L,*},z^{m+k}}_R =\bdzero.
$$
The other three cases are considered similarly. From \eqref{cmv4},
\eqref{cmv5}, and (i), we get \eqref{cmv1} for $k\not=m$.
\end{proof}

{}From the above proposition and the $\norm{\cdot}_\infty$-density of
Laurent
polynomials in $C(\partial\bbD)$, it follows that  $\{\chi_k\}_{k=0}^
\infty$ and $
\{x_k\}_{k=0}^\infty$ are right orthonormal modula bases in $\h$,
that is,
any element $f\in \h$ can be represented as a
\begin{equation}
f=\sum_{k=0}^\infty \chi_k\ang{\chi_k,f}_R =\sum_{k=0}^\infty x_k\ang
{x_k,f}_R.
\label{cmv6}
\end{equation}

\subsubsection{The CMV matrix}

Consider the matrix of the right homomorphism $f(z) \mapsto zf(z)$
with respect to the
basis $\{\chi_k\}$. Denote ${\mathcal C}_{nm}=\ang{\chi_n, z \chi_m}_
{R}$. The matrix
${\mathcal C}$ is unitary in the following sense:
$$
\sum_{k=0}^\infty {\mathcal C}_{kn}^\dagger {\mathcal C}_{km} = \sum_
{k=0}^\infty
{\mathcal C}_{nk} {\mathcal C}_{mk}^\dagger=\delta_{nm} \bdone.
$$
The proof follows from \eqref{cmv6}:
\begin{align*}
\delta_{nm}\bdone & = \ang{z\chi_n,z\chi_m}_R = \biggl< \!\!\!
\biggl< \sum_{k=0}^\infty \chi_k
\ang{\chi_k,z\chi_n}_R,z\chi_m \biggr>\!\!\! \biggr>_R =  \sum_{k=0}^
\infty
{\mathcal C}_{kn}^\dagger {\mathcal C}_{km},
\\
\delta_{nm}\bdone &= \ang{\bar z\chi_n,\bar z\chi_m}_R = \biggl< \!\!
\! \biggl< \sum_{k=0}^\infty
\chi_k \ang{\chi_k,\bar z\chi_n}_R, \bar z\chi_m\biggr>\!\!\!
\biggr>_R = \sum_{k=0}^\infty
{\mathcal C}_{nk} {\mathcal C}_{mk}^\dagger.
\end{align*}

We note an immediate consequence:

\begin{lemma}\label{sol}
Let $\abs{z}\leq1$. Then, for every $m \ge 0$,
$$
\sum_{n = 0}^\infty \chi_n(z) \mathcal{C}_{nm} = z \chi_m(z), \qquad
\sum_{n = 0}^\infty
\mathcal{C}_{mn} \chi_n(1/\bar z)^\dagger = z \chi_m(1/\bar z)^\dagger.
$$
\end{lemma}

\begin{proof}
First note that the above series contains only finitely many non-zero
terms. Expanding
$f(z) = z \chi_n$ according to \eqref{cmv6}, we see that
$$
z \chi_n(z) = \sum_{k=0}^\infty \chi_k(z) \ang{\chi_k,z \chi_n}_R =
\sum_{k=0}^\infty
\chi_k(z) \mathcal{C}_{kn}
$$
which is the first identity. Next, taking adjoints, we get
$$
\bar z \chi_n(z)^\dagger = \sum_{k=0}^\infty \mathcal{C}_{kn}^\dagger
\chi_k(z)^\dagger
$$
which yields
\begin{align*}
\bar z \sum_{n = 0}^\infty \mathcal{C}_{mn} \chi_n(z)^\dagger & =
\sum_{n =
0}^\infty
\mathcal{C}_{mn} \sum_{k=0}^\infty \mathcal{C}_{kn}^\dagger \chi_k(z)^
\dagger \\
& = \sum_{k=0}^\infty \biggl(\, \sum_{n = 0}^\infty \mathcal{C}_{mn} \,
\mathcal{C}_{kn}^\dagger \biggr) \chi_k(z)^\dagger \\
& = \sum_{k=0}^\infty \delta_{mk} \chi_k(z)^\dagger = \chi_m(z)^\dagger.
\end{align*}
Replacing $z$ by $1/\bar z$, we get the required statement.
\end{proof}

\subsubsection{The $\mathcal{LM}$-representation}

Using \eqref{cmv6} for $f=\chi_m$, we obtain:
\begin{equation}
\mathcal{C}_{nm}=\ang{\chi_n, z \chi_m}_{R} = \sum_{k=0}^\infty \ang
{\chi_n,z x_k}_R
\ang{x_k,\chi_m}_R = \sum_{k=0}^\infty \mathcal{L}_{nk} \mathcal{M}_
{km}. \label{cmv7}
\end{equation}

Denote by $\Theta(\alpha)$ the $2l \times 2l$ unitary matrix
$$
\Theta(\alpha) = \begin{pmatrix} \alpha^\dagger & \rho^L \\ \rho^R &
- \alpha
\end{pmatrix}.
$$
Using the Szeg\H o recursion formulas \eqref{szright1} and \eqref
{szright4}, we get
\begin{align}
z \varphi_n^R &= \varphi_{n+1}^R\rho_n^R+\varphi_n^{L,*} \alpha_n^
\dagger,
\label{szright5}
\\
\varphi_{n+1}^{L,*} &= \varphi_{n}^{L,*}\rho_n^L-\varphi_{n+1}^{R}
\alpha_n.
\label{szright6}
\end{align}
Taking $n=2k$ and multiplying by $z^{-k}$, we get
\begin{align*}
z x_{2k} &= \chi_{2k} \alpha_{2k}^\dagger+\chi_{2k+1} \rho_{2k}^R,
\\
z x_{2k+1} &= \chi_{2k} \rho_{2k}^L-\chi_{2k+1} \alpha_{2k}.
\end{align*}
It follows that the matrix $\mathcal{L}$ has the structure
$$
\mathcal{L}  = \Theta(\alpha_0) \oplus \Theta(\alpha_2) \oplus \Theta
(\alpha_4) \oplus
\cdots .
$$
Taking $n=2k-1$ in \eqref{szright5}, \eqref{szright6} and multiplying
by $z^{-k}$, we get
\begin{align*}
\chi_{2k-1} &= x_{2k-1} \alpha_{2k-1}^\dagger+x_{2k} \rho_{2k-1}^R,
\\
\chi_{2k} &= x_{2k-1} \rho_{2k-1}^L-x_{2k} \alpha_{2k-1}.
\end{align*}
It follows that the matrix $\mathcal{M}$ has the structure
\begin{equation}
\mathcal{M} = \bdone \oplus \Theta(\alpha_1) \oplus \Theta(\alpha_3)
\oplus \cdots.
\label{*3}
\end{equation}
Substituting this into \eqref{cmv7}, we obtain:
\begin{equation}\label{cmvmatrix}
\mathcal{C} =
\begin{pmatrix} {}& \alpha_0^\dagger & \rho_0^L \alpha_1^\dagger  &
\rho_0^L \rho_1^L & \bdzero & \bdzero & \cdots & {}
\\
{}& \rho_0^R & -\alpha_0 \alpha_1^\dagger & - \alpha_0 \rho_1^L &
\bdzero & \bdzero & \cdots & {}
\\
{}& \bdzero &  \alpha_2^\dagger \rho_1^R & -\alpha_2^\dagger \alpha_1
    & \rho_2^L \alpha_3^\dagger  & \rho_2^L \rho_3^L & \cdots & {}
\\
{}& \bdzero & \rho_2^R \rho_1^R  & -\rho_2^R \alpha_1
    & -\alpha_2 \alpha_3^\dagger & -\alpha_2 \rho_3^L & \cdots & {}
\\
{}& \bdzero & \bdzero & \bdzero & \alpha_4^\dagger \rho_3^R  & -
\alpha_4^\dagger \alpha_3 & \cdots & {}
\\
{}& \vdots & \vdots & \vdots & \vdots & \vdots & \ddots
\end{pmatrix}.
\end{equation}
We note that the analogous formula to this in \cite{SimonRev},
namely, (4.30), is incorrect!
The order of the factors below the diagonal is wrong there.

\subsection{The Resolvent of the CMV Matrix}

We begin by studying solutions to the equations
\begin{alignat}{2}
\sum_{k=0}^\infty {\mathcal C}_{mk} w_k &= zw_m, \quad & m &\geq 2,
\label{r1}
\\
\sum_{k=0}^\infty \tilde w_k {\mathcal C}_{km} &= z\tilde w_m, \quad
&  m &\geq 1. \label{r2}
\end{alignat}
Let us introduce the following functions:
\begin{align*}
\tilde x_n(z) & =\chi_n(1/\bar{z})^\dagger, \\
\Upsilon_{2n}(z) & =-z^{-n}\psi_{2n}^{L,*}(z), \\
\Upsilon_{2n-1}(z)& = z^{-n+1}\psi_{2n-1}^{R}(z), \\
y_{2n}(z) & = - \Upsilon_{2n}(1/\bar{z})^\dagger =z^{-n}\psi_{2n}^{L}
(z), \\
y_{2n-1}(z) & = -\Upsilon_{2n-1}(1/\bar{z})^\dagger=-z^{-n}\psi_{2n-1}
^{R,*}(z), \\
p_n(z) & = y_n(z)+\tilde x_n(z)F(z), \\
\pi_n(z) & = \Upsilon_n(z)+F(z)\chi_n(z).
\end{align*}

\begin{proposition}\label{sol1}
Let $z\in\D\setminus \{0\}$.
\begin{SL}
\item[{\rm{(i)}}] For each $n\geq 0$, a pair of values $(\tilde w_
{2n}, \tilde w_{2n+1})$
uniquely determines a solution $\tilde w_n$ to \eqref{r2}. Also, for
any pair of values
$(\tilde w_{2n},\tilde w_{2n+1})$  in $\calM_l$, there exists a
solution $\tilde w_{n}$ to
\eqref{r2} with these values at $(2n, 2n+1)$.

\item[{\rm{(ii)}}] The set of solutions $\tilde w_n$ to \eqref{r2}
coincides with the set of
sequences
\begin{equation}
\tilde w_n (z)= a \chi_n(z) + b\pi_n(z) \label{r3}
\end{equation}
where $a,b$ range over $\calM_l$.

\item[{\rm{(iii)}}] A solution \eqref{r3} is in $\ell^2$ if and only
if $a=0$.

\item[{\rm{(iv)}}] A solution \eqref{r3} obeys \eqref{r2} for all $m
\geq0$ if and only if $b=0$.
\end{SL}
\end{proposition}

\begin{proposition}\label{sol2}
Let $z\in\D\setminus \{0\}$.
\begin{SL}
\item[{\rm{(i)}}] For each $n\geq 1$, a pair of values $(w_{2n-1},
w_ {2n})$ uniquely determines a solution $w_n$ to \eqref{r1}.
Also, for any pair of values $(w_ {2n-1}, w_{2n})$ in $\calM_l$,
there exists a solution $w_{n}$ to \eqref{r1} with these values at
$(2n-1, 2n)$.

\item[{\rm{(ii)}}] The set of solutions $w_n$ to \eqref{r1} coincides
with the set of sequences
\begin{equation}
w_n (z)= \tilde x_n(z)a + p_n(z)b \label{r4}
\end{equation}
where $a,b$ range over $\calM_l$.

\item[{\rm{(iii)}}] A solution \eqref{r4} is in $\ell^2$ if and only
if $a=0$.

\item[{\rm{(iv)}}] A solution \eqref{r4} obeys \eqref{r1} for all $m
\geq0$ if and only if $b=0$.
\end{SL}
\end{proposition}

\begin{proof}[Proof of Proposition~\ref{sol1}]
(i) The matrix ${\mathcal C}-z$ can be written in the form
\begin{equation}
{\mathcal C}-z=
\begin{pmatrix}
A_0 & B_0 & \bdzero & \bdzero & \cdots \\
\bdzero & A_1 & B_1 & \bdzero & \cdots \\
\bdzero & \bdzero & A_2 & B_2 & \cdots \\
\vdots & \vdots & \vdots & \vdots & \ddots\\
\end{pmatrix}
\label{r5}
\end{equation}
where
$$
A_0=
\begin{pmatrix}
\alpha_0^\dagger-z\\ \rho_0^R
\end{pmatrix}
\qquad A_n=
\begin{pmatrix}
\alpha^\dagger_{2n} \rho_{2n-1}^R & -\alpha_{2n}^\dagger \alpha_
{2n-1}-z \\
\rho_{2n}^R\rho_{2n-1}^R & -\rho_{2n}^R \alpha_{2n-1}
\end{pmatrix}
$$
$$
B_n=
\begin{pmatrix}
\rho_{2n}^L\alpha_{2n+1}^\dagger & \rho_{2n}^L \rho_{2n+1}^L \\
-\alpha_{2n}\alpha_{2n+1}^\dagger-z & -\alpha_{2n}\rho_{2n+1}^L
\end{pmatrix}.
$$
Define $\wti W_n=(\tilde w_{2n}, \tilde w_{2n+1})$ for
$n=0,1,2,\dots$. Then \eqref{r2} for $m=2n+1$, $2n+2$ is
equivalent to
$$
\wti W_{n} B_n +\wti W_{n+1} A_{n+1}=\bdzero.
$$
It remains to prove that the $2l\times 2l$ matrices $A_j$, $B_j$
are invertible. Suppose that for some $x,y\in \C^l$, $A_n
\binom{x}{y}=\binom{0}{0}$. This is equivalent to the system
\begin{align*}
\alpha_{2n}^\dagger \rho_{2n-1}^R x - \alpha_{2n}^\dagger \alpha_
{2n-1}y-zy&=\bdzero,
\\
\rho_{2n}^R\rho_{2n-1}^Rx - \rho_{2n}^R\alpha_{2n-1}y&=\bdzero.
\end{align*}
The second equation of this system yields
$\rho_{2n-1}^Rx=\alpha_{2n-1} y$ (since $\rho_{2n}^R$ is
invertible), and upon substitution into the first equation, we get
$y=x=0$. Thus, $\ker(A_n)=\{0\}$. In a similar way, one proves
that $\ker(B_n)=\{0\}$.

\smallskip
(ii) First note that $\tilde w_n=\chi_n$ is a solution to \eqref{r2}
by Lemma~\ref{sol}.
Let us check that $\tilde w_n=\Upsilon_n$ is also a solution. If $U_
{km}=(-1)^k
\delta_{km}$, then $(U{\mathcal C}U)_{km}$ for $m\geq1$ coincides
with the CMV matrix
corresponding to the coefficients $\{-\alpha_n\}$. Recall that $
\psi_n^{L,R}$ are the
orthogonal polynomials $\varphi_n^{L,R}$, corresponding to the
coefficients
$\{-\alpha_n\}$. Taking into account the minus signs in the
definition of $\Upsilon_n$,
we see that $\tilde w_n=\Upsilon_n$ solves \eqref{r2} for $m\geq1$.
It follows that any
$\tilde w_n$ of the form \eqref{r3} is a solution to \eqref{r2}.

Let us check that any solution to \eqref{r2} can be represented as
\eqref{r3}. By (i), it
suffices to show that for any $\tilde w_0, \tilde w_1$, there exist
$a,b\in\calM_l$ such
that
\begin{align*}
a\chi_0(z)+b\pi_0(z)&=\tilde w_0,
\\
a\chi_1(z)+b\pi_1(z)&=\tilde w_1.
\end{align*}
Recalling that $\chi_0=1$, $\Upsilon_0=-1$,
$\Upsilon_1(z)=(z+\alpha_0^\dagger)(\rho_0^R)^{-1}$,
$\chi_1(z)=(z-\alpha_0^\dagger)(\rho_0^R)^{-1}$, we see that the
above system can be
easily solved for $a,b$ if $z\not=0$.

\smallskip
(iii) Let us prove that the solution $\pi_n$ is square integrable. We
will consider
separately the sequences $\pi_{2n}$ and $\pi_{2n-1}$ and prove that
they both belong to
$\ell^2$. By \eqref{a2} and \eqref{a2a}, we have
\begin{align}
\psi_n^R(z)+F(z)\varphi_n^R(z) &= \int
\frac{e^{i\theta}+z}{e^{i\theta}-z} \, d\mu(\theta)\varphi_n^R(e^{i
\theta}), \label{r6}
\\
\psi_n^{L,*}(z)-F(z)\varphi_n^{L,*}(z) &= -z^n\int
\frac{e^{i\theta}+z}{e^{i\theta}-z}\, d\mu(\theta)\varphi_n^L(e^{i
\theta})^\dagger.
\label{r7}
\end{align}
Taking $n=2k$ in \eqref{r7} and $n=2k-1$ in \eqref{r6}, we
get
\begin{align}
\pi_{2k}(z) &= z^k \int  \frac{e^{i\theta}+z}{e^{i\theta}-z}\,
d\mu(\theta)\varphi_{2k}^L (e^{i\theta})^\dagger, \label{r8}
\\
\pi_{2k-1}(z) &=
z^{-k+1} \int  \frac{e^{i\theta}+z}{e^{i\theta}-z} \, d\mu(\theta)
\varphi_{2k-1}^R(e^{i\theta}).
\label{r9}
\end{align}
As $\varphi_{2k}^L$ is an orthonormal sequence, using the Bessel
inequality, from \eqref{r8} we immediately get that $\pi_{2k}$ is
in $\ell^2$.

Consider the odd terms $\pi_{2k-1}$. We claim that
\begin{equation}
z^{-k+1} \int
\frac{e^{i\theta}+z}{e^{i\theta}-z}\, d\mu(\theta)\varphi_{2k-1}^R(e^
{i\theta}) = \int
\frac{e^{i\theta}+z}{e^{i\theta}-z}\, d\mu(\theta)e^{i(-k+1)\theta}
\varphi_{2k-1}^R(e^{i\theta}).
\label{r10}
\end{equation}
Indeed, using the right orthogonality of $\varphi_{2k-1}^R$ to $e^{im
\theta}$,
$m=0,1,\dots,2k-2$, we get
\begin{align*}
\int  \frac{e^{i\theta}+z}{e^{i\theta}-z}\, d\mu(\theta) \varphi_
{2k-1}^R (e^{i\theta})
& = \biggl<\!\!\! \biggl< 1+2\sum_{m=1}^\infty \bar z^m e^{im\theta},
\varphi_{2k-1}^R\biggr>\!\!\! \biggr>_R \\
& = \biggl<\!\!\! \biggl<  2 \!\!\! \sum_{m=2k-1}^\infty \bar z^m
e^{im\theta},\varphi_{2k-1}^R\biggr>\!\!\! \biggr>_R
\end{align*}
and
\begin{align*}
\int \frac{e^{i\theta}+z}{e^{i\theta}-z}z^{k-1}\, &e^{i(-k+1)\theta} \,
d\mu(\theta)\varphi_{2k-1}^R(e^{i\theta}) = \\
& = \biggl<\!\!\! \biggl< \bar z^{k-1}e^{i(k-1)\theta} \biggl(1+2\sum_
{m=1}^\infty \bar z^m
e^{im\theta}\biggr),\varphi_{2k-1}^R\biggr>\!\!\! \biggr>_R \\
& = \biggl<\!\!\! \biggl< 2 \!\!\! \sum_{m=2k-1}^\infty \bar z^m
e^{im\theta},\varphi_{2k-1}^R\biggr>\!\!\! \biggr>_R
\end{align*}
which proves \eqref{r10}. The identities \eqref{r10} and \eqref{r9}
yield
$$
\pi_{2k-1}(z) = \biggl<\!\!\! \biggl<
\frac{e^{-i\theta}+\bar z}{e^{-i\theta}-\bar z},\chi_{2k-1}\biggr>\!\!
\! \biggr>_R
$$
and, since $\chi_{2k-1}$ is a right orthogonal sequence, the Bessel
inequality ensures
that $\pi_{2k-1}(z)$ is in $\ell^2$. Thus, $\pi_k(z)$ is in $\ell^2$.

Next, as in the proof  of Theorem~\ref{th.solu}, using the CD
formula, we check that the sequence $\norm{\varphi_n^{L,*}(z)}_R$
is bounded below and therefore the sequence $\chi_{2n}(z)$ is not
in $\ell^2$. This proves the statement (iii).

\smallskip
(iv) By Lemma~\ref{sol}, the solution $\chi_n(z)$ obeys \eqref{r2}
for all $m\geq0$. It is easy to check directly that the solution
$\pi_n(z)$ does not obey \eqref{r2} for $m=0$ if $z\not=0$. This
proves the required statement.
\end{proof}

\begin{proof}[Proof of Proposition~\ref{sol2}]
(i) For $j=1,2,\dots$, define $W_j=(w_{2j-1}, w_{2j})$. Then, using
the block structure
\eqref{r5}, we can rewrite \eqref{r1} for $m=2j, 2j+1$ as $A_j  W_j
+B_{j}W_{j+1}=\bdzero$. By
the proof of Proposition~\ref{sol1}, the matrices $A_j$ and $B_j$ are
invertible, which
proves (i).

\smallskip
(ii) Lemma~\ref{sol} ensures that $\tilde x_n(z)$ is a solution of
\eqref{r1}. As in the
proof of Proposition~\ref{sol1}, by considering the matrix $(U
{\mathcal C}U)_{km}$, one
checks that $y_n(z)$ is also a solution to \eqref{r1}.

Let us prove that any solution to \eqref{r1} can be represented in
the form \eqref{r4}.
By (i), it suffices to show that for any $w_1$, $w_2$, there exist
$a,b\in\calM_l$ such that
\begin{align*}
\tilde x_1(z) a + p_1(z) b & = w_1,
\\
\tilde x_2(z) a + p_2(z) b & = w_2.
\end{align*}
We claim that this system of equations can be solved for $a$, $b$.
Here are the main
steps. Substituting the definitions of $\tilde x_n(z)$ and $p_n(z)$,
we rewrite this
system as
\begin{align*}
\varphi_1^{R,*} (a+F(z)b) - \psi_1^{R,*}(z) b & = z w_1,
\\
\varphi_2^{L} (a+F(z)b) + \psi_2^{L}(z) b & = z w_2.
\end{align*}
Using Szeg\H o recurrence, we can substitute the expressions for
$\varphi_2^{L}$, $\psi_2^{L}$, which helps rewrite our system as
\begin{align*}
\varphi_1^{R,*} (a+F(z)b) - \psi_1^{R,*}(z) b & = z w_1,
\\
\varphi_1^{L} (a+F(z)b) + \psi_1^{L}(z) b & = \rho_1^L
w_2+\alpha_1^\dagger  w_1.
\end{align*}
Substituting explicit formulas for $\varphi_1^{R,*}$,
$\varphi_1^{L}$, $ \psi_1^{R,*}$, $\psi_1^{L}$, and expressing
$F(z)$ in terms of $f(z)$, we can rewrite this as
\begin{align*}
(\rho_0^R)^{-1}(1-\alpha_0z)a+2z(\rho_0^R)^{-1}(f(z)-\alpha_0)(\bdone-
zf(z))^{-1}b
& = z w_1,
\\
(\rho_0^L)^{-1}(z-\alpha_0^\dagger)a+2z(\rho_0^L)^{-1}(\bdone-
\alpha_0^\dagger
f(z))(\bdone-zf(z))^{-1}b & = \rho_1^L w_2+\alpha_1^\dagger  w_1.
\end{align*}
Denote
\begin{align*}
a_1 &=(\rho_0^R)^{-1}(1-\alpha_0z)a,
\\
b_1 &=2z(\rho_0^L)^{-1}(\bdone-\alpha_0^\dagger
f(z))(\bdone-zf(z))^{-1}b.
\end{align*}
Then in terms of $a_1$, $b_1$, our system can be rewritten as
\begin{align*}
a_1+X_1b_1 & = z w_1,
\\
X_2a_1+b_1 & = \rho_1^L w_2+\alpha_1^\dagger  w_1,
\end{align*}
where
\begin{align*}
X_1 &=(\rho_0^R)^{-1}(f(z)-\alpha_0)(\bdone-\alpha_0
f(z))^{-1}\rho_0^L,
\\
X_2 &=(\rho_0^L)^{-1}(z-\alpha_0^\dagger)(\bdone-\alpha_0
z)^{-1}\rho_0^R.
\end{align*}
Since $\norm{f(z)}<1$ and $\abs{z}<1$, we can apply
Corollary~\ref{C1.3.5}, which yields $\norm{X_1}<1$ and
$\norm{X_2}<1$. It follows that our system can be solved for
$a_1$, $b_1$.

\smallskip
(iii) As $p_n(z)=-\pi_n(1/\bar z)^\dagger$, by Proposition~\ref
{sol1}, we get that
$p_n(z)$ is in $\ell^2$. In the same way, as $\tilde
x_n(z)=\chi_n(1/\bar z)^\dagger$, we get that $\tilde x_n(z)$ is not
in $\ell^2$.

\smallskip
(iv) By Lemma~\ref{sol}, the solution $\tilde x_n(z)$ obeys \eqref
{r1} for all $m\geq0$.
Using the explicit formula for $y_n(z)$, one easily checks that the
solution $y_n(z)$
does not obey \eqref{r1} for $m=0,1$.
\end{proof}

\begin{theorem} \lb{T3.21}
We have for $z \in \D$,
$$
[(\mathcal{C} - z)^{-1}]_{k,l} =
\begin{cases} (2z)^{-1} \tilde x_k(z) \pi_l(z), & l > k \text{
or } k = l \text{ even}, \\
(2z)^{-1} p_k(z) \chi_l(z), & k > l \text{ or } k = l \text{
odd}.
\end{cases}
$$
\end{theorem}

\begin{proof}
Fix $z \in \D$. Write $G_{k,l}(z) = [(\mathcal{C} -
z)^{-1}]_{k,l}$. Then $G_{\bddot , l}(z)$ is equal to
$(\mathcal{C} - z)^{-1} \delta_{l}$, which means that $G_{k,l}(z)$
solves \eqref{r1} for $m \not= l$. Since $G_{\bddot , l}(z)$ is
$\ell^2$ at infinity and obeys the equation at $m = 0$, we see
that it is a right-multiple of $p$ for large $k$ and a
right-multiple of $\tilde x$ for small $k$. Thus,
$$
G_{k,l}(z) =
\begin{cases} \tilde x_k(z) a_l(z), & k < l \text{ or } k = l  \text
{ even}, \\
p_k(z) b_l(z), & k > l \text{ or } k = l \text{ odd}.
\end{cases}
$$
Similarly,
$$
G_{k,l}(z) = \begin{cases} \tilde b_k(z) \pi_l(z), & k < l \text
{ or } k = l \text{ even}, \\
\tilde a_k(z) \chi_l(z), & k > l \text{ or } k = l \text{ odd}.
\end{cases}
$$
Equating the two expressions, we find
\begin{alignat}{2}
\label{rfp1} \tilde x_k(z) a_l(z) &= \tilde b_k(z) \pi_l(z)  \quad &
k &< l \text{ or } k = l \text{ even}, \\
\label{rfp2} p_k(z) b_l(z) &= \tilde a_k(z) \chi_l(z)  \quad& k &> l
\text{ or } k = l
\text{ odd}.
\end{alignat}
Putting $k = 0$ in \eqref{rfp1} and setting $\tilde b_0(z) = c_1(z)$,
we find $a_l(z) =
c_1(z) \pi_l(z)$. Putting $l = 0$ in \eqref{rfp2} and setting $b_0(z)
= c_2(z)$, we find
$p_k(z) c_2(z) = \tilde a_k(z)$. Thus,
$$
G_{k,l}(z) =
\begin{cases} \tilde x_k(z) c_1(z) \pi_l(z), & k < l \text{ or } k =
l  \text{ even}, \\
p_k(z) c_2(z) \chi_l(z), & k > l \text{ or } k = l \text{ odd}.
\end{cases}
$$
We claim that $c_1(z) = c_2(z) = (2z)^{-1} \bdone$. Consider the
case $k = l = 0$. Then, on the one hand, by the definition,
\begin{align}
G_{0,0}(z) & = \int \frac{1}{e^{i\theta} - z} \, d\mu (e^{i\theta})
\notag
\\
& = \int (2z)^{-1} \biggl[ \frac{e^{i\theta} + z}{e^{i\theta} - z} -
1 \biggr] d\mu (e^{i\theta}) \notag
\\
&= (2z)^{-1} (F(z) - \bdone)
\label{spthm}
\end{align}
and on the other hand,
$$
G_{0,0}(z) = \tilde x_0(z) c_1(z) \pi_0(z) = c_1(z) (F(z) - \bdone).
$$
This shows $c_1(z) = (2z)^{-1} \bdone$. Next, consider the case $k=1
$, $l=0$. Then, on the one
hand, by the definition,
$$
G_{1,0}(z) = \ang{\chi_1, (e^{i\theta}-z)^{-1}\chi_0}_{R}
$$
and on the other hand,
$$
G_{1,0}(z) = p_1(z) c_2(z) \chi_0(z).
$$
Let us calculate the expressions on the right-hand side. We have
\begin{equation}
p_1(z) c_2(z) \chi_0(z) = (\rho_0^R)^{-1}(-z^{-1}-\alpha_0+(z^{-1}-
\alpha_0)F(z))c_2(z)
\label{r11}
\end{equation}
and
\begin{align*}
\ang{\chi_1, &(e^{i\theta}-z)^{-1}\chi_0}_{R} = \\
& = (\rho_0^R)^{-1} \int (e^{-i\theta}-\alpha_0)d\mu(\theta)(e^{i
\theta}-z)^{-1}
\\
& = (\rho_0^R)^{-1}\int [z^{-1}(e^{i\theta}-z)^{-1}-z^{-1}e^{-i
\theta}-\alpha_0
(e^{i\theta}-z)^{-1}] \, d\mu(\theta)
\\
& = (\rho_0^R)^{-1}\biggl[\frac{1}{2z2}\, (F(z)-\bdone)-\frac{1}{2z}
\, \alpha_0 (F(z)-\bdone) -
\frac{1}{z}\int e^{-i\theta}\, d\mu(\theta)\biggr].
\end{align*}
Taking into account the identity
$$
\int e^{-i\theta}d\mu(\theta)=\alpha_0
$$
(which can be obtained, e.g., by expanding $\ang{\varphi_1^R,
\varphi_0^R}_R=0$),
we get
$$
\ang{\chi_1, (e^{i\theta}-z)^{-1}\chi_0}_{R} = \frac{1}{2z}\,
(\rho_0^R)^{-1}
(-z^{-1}-\alpha_0+(z^{-1}-\alpha_0)F(z)).
$$
Comparing this with \eqref{r11}, we get $c_2(z)=(2z)^{-1}\bdone$.
\end{proof}

As an immediate corollary, evaluating the kernel on the diagonal for
even and odd
indices, we obtain the formulas
\begin{align}
\int \varphi_{2n}^L(e^{i\theta})\, \frac{d\mu(\theta)}{e^{i\theta}-z}\,
\varphi_{2n}^L(e^{i\theta})^\dagger &= -\frac1{2 z^{2n+1}}\, \varphi_
{2n}^L(z)
u_{2n}^{L,*}(z), \label{r12}
\\
\int \varphi_{2n-1}^R(e^{i\theta})^\dagger\, \frac{d\mu(\theta)}{e^{i
\theta}-z}\,
\varphi_{2n-1}^R(e^{i\theta}) & = -\frac1{2 z^{2n}}\, u_{2n-1}^{R,*}(z)
\varphi_{2n-1}^R(z). \label{r13}
\end{align}
Combining this with \eqref{wronsk1} and \eqref{wronsk2}, we find
\begin{align}
u_n^L(z) \varphi_n^{L,*}(z) + \varphi_n^L(z) u_n^{L,*}(z) & = 2 z^n,
\label{wronsk5} \\
\varphi_n^{R,*}(z) u_n^R(z) + u_n^{R,*}(z) \varphi_n^R(z) & = 2 z^n.
\label{wronsk6}
\end{align}

\subsection{Khrushchev Theory} \lb{s3.13}

Among the deepest and most elegant methods in OPUC are those of
Khrushchev
\cite{Kh2000,Khr,KhGo}. We have not been able to extend them to
MOPUC! We regard
their extension as an important open question; we present the first
very partial
steps here.

Let
$$
\Omega = \{ \theta \,\colon \det w(\theta) > 0 \}.
$$

\begin{theorem}
For every $n \ge 0$,
$$
\{ \theta : f_n(e^{i\theta})^\dagger f_n(e^{i\theta}) < \bdone \} =
\Omega
$$
up to a set of zero Lebesgue measure.

Consequently,
\begin{equation}\label{fnintlower}
\int \| f_n (e^{i\theta}) \| \, \frac{d\theta}{2\pi} \ge 1 - \frac{|
\Omega|}{2\pi}\, .
\end{equation}
\end{theorem}

\begin{proof}
Recall that, by Proposition~\ref{prop.3.16}, up to a set of zero
Lebesgue measure,
$$
\{ \theta\, \colon f_0(e^{i\theta})^\dagger f_0(e^{i\theta}) < \bdone
\} =
\{ \theta\, \colon \det w(\theta) > 0 \}
$$
so, by induction, it suffices to show that, up to a set of zero
Lebesgue measure,
$$
\{ \theta : f_0(e^{i\theta})^\dagger f_0(e^{i\theta}) < \bdone \}
= \{ \theta : f_1(e^{i\theta})^\dagger f_1(e^{i\theta}) < \bdone
\}.
$$
This in turn follows from the fact that the Schur algorithm, which
relates the two functions, preserves the property $g^\dagger g <
\bdone$.

Notice that away from $\Omega$, $f_n (e^{i\theta})$ has norm one and
therefore,
$$
\int \| f_n (e^{i\theta}) \| \, \frac{d\theta}{2\pi} \ge \int_
{\Omega^c} \| f_n
(e^{i\theta}) \| \, \frac{d\theta}{2\pi} =  \int_{\Omega^c} 1 \, \frac
{d\theta}{2\pi}
$$
which yields \eqref{fnintlower}.
\end{proof}

Define
$$
b_n(z;d\mu) = \varphi_n^L(z;d\mu) \varphi_n^{R,*}(z;d\mu)^{-1}.
$$

\begin{proposition}
\begin{SL}
\item[{\rm{(a)}}] $b_{n+1} = (\rho_n^L)^{-1} (z b_n - \alpha_n^
\dagger) (\bdone - z
\alpha_n b_n)^{-1} \rho_n^R$.

\item[{\rm{(b)}}] The Verblunsky coefficients of $b_n$ are $(-
\alpha_{n-1}^ \dagger, - \alpha_{n-2}^\dagger, \dots , -
\alpha_0^\dagger, \bdone)$.
\end{SL}
\end{proposition}

\begin{proof}
(a) By the Szeg\H{o} recursion, we have that
\begin{align*}
b_{n+1} & = \varphi_{n+1}^L (\varphi_{n+1}^{R,*})^{-1} \\
& = ((\rho_n^L)^{-1} z \varphi_n^L - (\rho_n^L)^{-1}
\alpha_n^\dagger \varphi_n^{R,*}) ((\rho_n^R)^{-1}
\varphi_n^{R,*} - z (\rho_n^R)^{-1} \alpha_n \varphi_n^L)^{-1} \\
& = (\rho_n^L)^{-1} (z \varphi_n^L - \alpha_n^\dagger \varphi_n^{R,*})
(\varphi_n^{R,*} - z \alpha_n \varphi_n^L)^{-1} \rho_n^R \\
& = (\rho_n^L)^{-1} (z \varphi_n^L (\varphi_n^{R,*})^{-1} -
\alpha_n^\dagger \varphi_n^{R,*} (\varphi_n^{R,*})^{-1})  \\
& \qquad \qquad \qquad (\varphi_n^{R,*} (\varphi_n^{R,*})^{-1} - z
\alpha_n \varphi_n^L
(\varphi_n^{R,*})^{-1})^{-1} \rho_n^R \\
& = (\rho_n^L)^{-1} ( z b_n - \alpha_n^\dagger) (\bdone - z
\alpha_n b_n )^{-1} \rho_n^R.
\end{align*}

\smallskip
(b) It follows from part (a) that the first Verblunsky coefficient
of $b_n$ is $- \alpha_{n-1}^\dagger$ and that its first Schur
iterate is $b_{n-1}$; compare Theorem~\ref{T3.12}. This gives the
claim by induction and the fact that $b_0 = \bdone$.
\end{proof}

\section{The Szeg\H{o} Mapping and the Geronimus Relations} \lb{s4}

In this chapter, we present the matrix analogue of the Szeg\H{o}
mapping and the resulting Geronimus relations. This establishes a
correspondence between certain matrix-valued measures on the unit
circle and matrix-valued measures on the interval $[-2,2]$ and,
consequently, a correspondence between Verblunsky coefficients and
Jacobi parameters. Throughout this chapter, we will denote
measures on the circle by $d\mu_C$ and measures on the interval by
$d\mu_I$.

The scalar versions of these objects are due to Szeg\H{o} \cite
{Sz22a} and Geronimus \cite{Ger46}. There are four proofs that we
know of: the original argument of Geronimus \cite{Ger46} based on
Szeg\H{o}'s formula in \cite{Sz22a}, a proof of Damanik--Killip
\cite{DKacta} using Schur functions, a proof of Killip--Nenciu
\cite{KN} using CMV matrices, and a proof of Faybusovich--Gekhtman
\cite{FG99} using canonical moments.

The matrix version of these objects was studied by
Yakhlef--Marcell\'an \cite{YM} who proved Theorem~\ref{T4.2} below
using the Geronimus--Szeg\H{o} approach. Our proof uses the
Killip--Nenciu--CMV approach. In comparing our formula with \cite
{YM}, one needs the following dictionary (their objects on the
left of the equal sign and ours on the right):
\begin{align*}
H_n &= -\alpha_{n+1}^\dagger ,\\
D_n &= A_n, \\
E_n &= B_{n+1}.
\end{align*}

Dette--Studden \cite{DS} have extended the theory of canonical
moments from OPRL to MOPRL.
It would be illuminating to use this to extend the proof that
Faybusovich--Gekhtman
\cite{FG99} gave of Geronimus relations for scalar OPUC to MOPUC.

Suppose $d\mu_C$ is a non-trivial positive semi-definite Hermitian
matrix measure on the unit circle that is invariant under $\theta
\mapsto - \theta$ (i.e., $z \mapsto \bar z =z^{-1}$). Then we
define the measure $d\mu_I$ on the interval $[-2,2]$ by
$$
\int f(x) \, d\mu_I(x) = \int f(2 \cos \theta) \, d\mu_C(\theta)
$$
for $f$ measurable on $[-2,2]$. The map
$$
\mathrm{Sz} : d\mu_C \mapsto d\mu_I
$$
is called the Szeg\H{o} mapping.

The Szeg\H{o} mapping can be inverted as follows. Suppose $d\mu_I$
is a non-degenerate positive semi-definite matrix measure on
$[-2,2]$. Then we define the measure $d\mu_C$ on the unit circle
which is invariant under $\theta \mapsto - \theta$ by
$$
\int g(\theta) \, d\mu_C(\theta) = \int g \left( \arccos (x/2)
\right) \, d\mu_I(x)
$$
for $g$ measurable on $\partial \D$ with $g(\theta) = g(-\theta)$.

We first show that for the measures on the circle of interest in this
section, the
Verblunsky coefficients are always Hermitian.

\begin{lemma}\label{alphaherm}
Suppose $d\mu_C$ is a non-trivial positive semi-definite
Hermitian matrix measure on the unit circle. Denote the associated
Verblunsky coefficients by $ \{ \alpha_n \}$. Then, $d\mu_C$ is
invariant under $\theta \mapsto - \theta$ if and only if
$\alpha_n^ \dagger = \alpha_n$ for every $n$.
\end{lemma}

\begin{proof}
For a polynomial $P$, denote $\tilde P(z)=P(\bar z)^\dagger$.

1. Suppose that $d\mu_C$ is invariant under $\theta \mapsto - \theta$.
Then we have
$$
\ang{f,g}_L=\ang{\tilde g,\tilde f}_R
$$
for all $f$, $g$. Inspecting the orthogonality conditions which
define $\Phi_n^L$ and $\Phi_n^R$, we see that
\begin{equation}
\tilde\Phi_n^L=\Phi^R_n \text{ and }
\ang{\Phi_n^L,\Phi_n^L}_L=\ang{\Phi_n^R,\Phi_n^R}_R. \label{4.1}
\end{equation}
Next, we claim that
\begin{equation}
\kappa_n^L=\kappa_n^{R,\dagger}.
\label{4.2}
\end{equation}
Indeed, recall the definition of $\kappa_n^L$, $\kappa_n^R$:
\begin{alignat*}{2}
\kappa_n^L&=u_n\ang{\Phi_n^L,\Phi_n^L}_L^{-1/2},
\text{ $u_n$ is unitary, } \kappa_{n+1}^L(\kappa_n^L)^{-1}>0, \quad &
\kappa_0^L &=1,
\\
\kappa_n^R&=\ang{\Phi_n^R,\Phi_n^R}_R^{-1/2}v_n,
\text{ $v_n$ is unitary, } (\kappa_n^R)^{-1}\kappa_{n+1}^R>0, \quad &
\kappa_0^R &=1.
\end{alignat*}
Using this definition and \eqref{4.1}, one can easily prove by induction
that $v_n=u_n^\dagger$ and therefore \eqref{4.2} holds true.

Next, taking $z=0$ in \eqref{szmonic}, we get
\begin{align*}
\alpha_n=-(\kappa_n^R)^{-1}\Phi_{n+1}^L(0)^\dagger(\kappa_n^L)^\dagger,
\\
\alpha_n=-(\kappa_n^R)^\dagger\Phi_{n+1}^R(0)^\dagger (\kappa_n^L)^{-1}.
\end{align*}
{}From here and \eqref{4.1}, \eqref{4.2}, we get $\alpha_n=\alpha_n^
\dagger$.

\smallskip
2. Assume $\alpha_n^\dagger=\alpha_n$ for all $n$. Then, by
Theorem~\ref{szethm}(c), we have $\rho_n^L=\rho_n^R$. It follows
that in this case the Szeg\H{o} recurrence relation is invariant
with respect to the change $\varphi_n^L\mapsto \tilde
\varphi_n^R$, $\varphi_n^R\mapsto \tilde \varphi_n^L$. It follows
that $\varphi_n^L= \tilde \varphi_n^R$, $\varphi_n^R=\tilde
\varphi_n^L$. In particular, we get
\begin{equation}
\ang{\varphi_n^L,\varphi_m^L}_L=\ang{\tilde\varphi_m^R, \tilde
\varphi_n^R}_R.
\label{4*}
\end{equation}
Now let $f$ and $g$ be any polynomials; we have
$$
f(z)=\sum_n f_n \varphi^L_n(z),
\quad
\tilde f(z)=\sum_n \tilde \varphi_n^L(z) f_n^\dagger,
$$
and a similar expansion for $g$. Using these expansions and \eqref
{4*}, we get
$$
\ang{f,g}_L=\ang{\tilde g,\tilde f}_R \text{ for all polynomials $f$,
$g$.}
$$
{}From here it follows that the measure $d\mu_C$ is invariant
under $\theta \mapsto - \theta$.
\end{proof}

Now consider two measures $d\mu_C$ and $d\mu_I =
\mathrm{Sz}(d\mu_C)$ and the associated CMV and Jacobi matrices.
What are the relations between the parameters of these matrices?

\begin{theorem} \lb{T4.2}
Given $d\mu_C$ and $d\mu_I = \mathrm{Sz}(d\mu_C)$ as above, the
coefficients of the associated CMV and Jacobi matrices satisfy the
Geronimus relations:
\begin{align}
B_{k+1} & = \sqrt{\bdone - \alpha_{2k-1}} \; \alpha_{2k}
\sqrt{\bdone-\alpha_{2k-1}} -
\sqrt{\bdone+\alpha_{2k-1}} \; \alpha_{2k-2} \sqrt{\bdone +\alpha_
{2k-1}}\, ,  \lb{4.3}\\
A_{k+1} & = \sqrt{\bdone - \alpha_{2k-1}} \sqrt{\bdone-\alpha_{2k}^2}
\sqrt{\bdone+\alpha_{2k+1}}\, . \lb{4.4}
\end{align}
\end{theorem}

\begin{remarks} 1. For these formulas to hold for $k=0$,
we set $\alpha_{-1} = - {\boldsymbol 1}$.

\smallskip
2. There are several proofs of the Geronimus relations in the
scalar case. We follow the proof given by Killip and Nenciu in \cite
{KN}.

\smallskip
3. These $A$'s are, in general, not type~1 or 2 or 3.
\end{remarks}

\begin{proof}
For a Hermitian $l \times l$ matrix $\alpha$ with $\| \alpha \| < 1$,
define the unitary
$2l \times 2l$ matrix $S(\alpha)$ by
$$
S(\alpha) = \frac{1}{\sqrt{2}} \begin{pmatrix} \sqrt{\bdone - \alpha}
& - \sqrt{\bdone + \alpha} \\
\sqrt{\bdone + \alpha} & \sqrt{\bdone - \alpha} \end{pmatrix}.
$$
Since $\alpha^\dagger = \alpha$, the associated $\rho^L$ and $\rho^R$
coincide and will
be denoted by $\rho$. We therefore have
$$
\Theta(\alpha) = \begin{pmatrix} \alpha & \rho \\ \rho & - \alpha
\end{pmatrix}
$$
and hence, by a straightforward calculation,
\begin{align*}
S(\alpha) \Theta(\alpha) S(\alpha)^{-1} & = \frac{1}{2} \begin
{pmatrix} \sqrt{\bdone - \alpha} & - \sqrt{\bdone + \alpha} \\
\sqrt{\bdone + \alpha} & \sqrt{\bdone - \alpha} \end{pmatrix}  \begin
{pmatrix} \alpha & \rho \\
\rho & - \alpha
\end{pmatrix}  \begin{pmatrix} \sqrt{\bdone - \alpha} & \sqrt{\bdone
+ \alpha} \\
- \sqrt{\bdone + \alpha} & \sqrt{\bdone - \alpha} \end{pmatrix} \\
& = \begin{pmatrix} - {\boldsymbol 1} & {\boldsymbol 0} \\
{\boldsymbol 0} & {\boldsymbol
1} \end{pmatrix}.
\end{align*}

Thus, if we define
$$
\mathcal{S} = \bdone\oplus S(\alpha_1)\oplus S(\alpha_3)\oplus \dots
$$
and
$$
\mathcal{R} = \bdone\oplus(-\bdone)\oplus \bdone\oplus(-\bdone)\oplus
\dots,
$$
it follows that (see \eqref{*3})
$$
\mathcal{S} \mathcal{M} \mathcal{S}^\dagger = \mathcal{R}.
$$

The matrix $\mathcal{L} \mathcal{M} + \mathcal{M} \mathcal{L}$ is
unitarily equivalent to
$$
\mathcal{A} = \mathcal{S} ( \mathcal{L} \mathcal{M} +
\mathcal{M} \mathcal{L} ) \mathcal{S}^\dagger
=
\mathcal{S}
\mathcal{L} \mathcal{S}^\dagger \mathcal{R} + \mathcal{R}
\mathcal{S} \mathcal{L} \mathcal{S}^\dagger.
$$
Observe that $\mathcal{A}$ is the direct sum of two block Jacobi
matrices. Indeed, it follows quickly from the explicit form of
$\mathcal{R}$ that the even-odd and odd-even block entries of
$\mathcal{A}$ vanish. Consequently, $\mathcal{A}$ is the direct sum
of its odd-odd and even-even block entries. We will call these two
block Jacobi matrices $J$ and $\tilde J$, respectively.

Consider $\mathcal C$ as an operator on $\h_v$.
Then $d\mu_C$ is the spectral measure of $\mathcal{C}$
in the following sense:
$$
[\mathcal{C} ^m]_{0,0}=\int e^{im\theta} d\mu_C(\theta);
$$
see \eqref{spthm}. Then, by the spectral theorem, $d\mu_C(-\theta)$
is the spectral measure of $\mathcal{C}^{-1}$ and so
$d\mu_I$ is the spectral measure of
$\mathcal{C}+\mathcal{C}^{-1}=
\mathcal{L} \mathcal{M} + (\mathcal{L} \mathcal{M})^{-1}=
\mathcal{L} \mathcal{M} +\mathcal{M} \mathcal{L}$. Since  $\mathcal{S}
$ leaves $\begin{pmatrix}
{\boldsymbol 1} & {\boldsymbol 0} & {\boldsymbol 0} & {\boldsymbol
0} & \cdots
\end{pmatrix}^t$ invariant, we see that $d\mu_I$ is the spectral
measure of $J$.

To determine the block entries of $J$, we only need to compute the
odd-odd block entries of $\mathcal{A}$. For $k \ge 0$, we have
\begin{align*}
\mathcal{A}_{2k+1,2k+1} & = \begin{pmatrix} \sqrt{\bdone+\alpha_
{2k-1}} &
\sqrt{\bdone - \alpha_{2k-1}} \end{pmatrix}  \begin{pmatrix} -
\alpha_{2k-2} & {\boldsymbol 0} \\ {\boldsymbol 0} & \alpha_{2k}
\end{pmatrix}
\begin{pmatrix} \sqrt{\bdone+\alpha_{2k-1}} \\ \sqrt{\bdone-\alpha_
{2k-1}}
\end{pmatrix} \\
& = \sqrt{\bdone - \alpha_{2k-1}} \; \alpha_{2k} \sqrt{\bdone-\alpha_
{2k-1}} -
\sqrt{\bdone+\alpha_{2k-1}} \; \alpha_{2k-2} \sqrt{\bdone+\alpha_{2k-1}}
\end{align*}
and
\begin{align*}
\mathcal{A}_{2k+1,2k+3} & = \begin{pmatrix} \sqrt{\bdone+\alpha_
{2k-1}} &
\sqrt{\bdone - \alpha_{2k-1}} \end{pmatrix}  \begin{pmatrix}
{\boldsymbol
0} & {\boldsymbol 0} \\ \rho_{2k} & {\boldsymbol 0}
\end{pmatrix}
\begin{pmatrix} \sqrt{\bdone+\alpha_{2k+1}} \\ \sqrt{\bdone-\alpha_{2k
+1}}
\end{pmatrix} \\
& = \sqrt{\bdone - \alpha_{2k-1}} \sqrt{\bdone-\alpha_{2k}^2}
\sqrt{\bdone+\alpha_{2k+1}}\, .
\end{align*}
The result follows.
\end{proof}

As in \cite{KN,S2}, one can also use this to get Geronimus relations
for the
second Szeg\H{o} map.

\section{Regular MOPRL} \lb{s5}

The theory of regular (scalar) OPs was developed by Stahl--Totik \cite
{StT} generalizing
a definition of Ullman \cite{Ull} for $[-1,1]$. (See Simon \cite
{EqMC} for a review and
more about the history.) Here we develop the basics for MOPRL; it is
not hard to do the same
for MOPUC.

\subsection{Upper Bound and Definition} \lb{s5.1}

\begin{theorem}\lb{T5.1} Let $d\mu$ be a non-trivial $l\times l$
matrix-valued measure
on $\bbR$ with $E=\supp(d\mu)$ compact. Then {\rm{(}}with $C(E)=$
logarithmic capacity
of $E${\rm{)}}
\begin{equation} \lb{5.1}
\limsup_{n\to\infty}\, \abs{\det(A_1\cdots A_n)}^{1/n} \leq C(E)^l.
\end{equation}
\end{theorem}

\begin{remarks} 1. $\abs{\det(A_1\cdots A_n)}$ is constant over
equivalent Jacobi parameters.

\smallskip
2. For the scalar case, this is a result of Widom \cite{Wid} (it
might be older) whose proof
extends to the matrix case.
\end{remarks}

\begin{proof} Let $T_n$ be the Chebyshev polynomials for $E$ (see
\cite[Appendix~B]{EqMC} for a
definition) and let $T_n^{(l)}$ be $T_n\otimes\bdone$, that is, the $l
\times l$ matrix polynomial
obtained by multiplying $T_n(x)$ by $\bdone$. $T_n^{(l)}$ is monic
so, by \eqref{2.7d} and \eqref{2.32a},
\begin{align*}
\abs{\det(A_1 \cdots A_n)}^{1/n}
& \leq \biggl| \det\biggl( \int \abs{T_n^{(l)}(x)}^2\, d\mu(x)\biggr)
\biggr|^{1/2n} \\
&\leq \sup_n\, \abs{T_n(x)}^{l/n}.
\end{align*}
By a theorem of Szeg\H{o} \cite{Sz24}, $\sup_n \abs{T_n(x)}^{1/n}\to C
(E)$, so \eqref{5.1} follows.
\end{proof}

\begin{definition} Let $d\mu$ be a non-trivial $l\times l$ matrix-valued measure with $E=\supp(d\mu)$
compact. We say $\mu$ is {\it regular\/} if equality holds in \eqref
{5.1}.
\end{definition}

\subsection{Density of Zeros} \lb{s5.2}

The following is a simple extension of the scalar results
(see \cite{DLS} or \cite[Sect.~2]{EqMC}):

\begin{theorem}\lb{T5.2} Let $d\mu$ be a regular measure with $E=\supp
(d\mu)$. Let $d\nu_n$
be the zero counting measure for $\det(p_n^L(x))$, that is, if $\{x_j^
{(n)}\}_{j=1}^{nl}$ are
the zeros of this determinant {\rm{(}}counting degenerate zeros
multiply{\rm{)}}, then
\begin{equation} \lb{5.2}
d\nu_n = \f{1}{nl} \sum_{j=1}^{nl} \delta_{x_j^{(n)}}.
\end{equation}
Then $d\nu_n$ converges weakly to $d\rho_E$, the equilibrium measure
for $E$.
\end{theorem}

\begin{remark} For a discussion of capacity, equilibrium measure,
quasi-every (q.e.), etc.,
see \cite[Appendix~A]{EqMC} and \cite{Helms,Land,Ran}.
\end{remark}

\begin{proof} By \eqref{2.54L} and \eqref{2.54x},
\begin{equation} \lb{5.3}
\abs{\det(p_n^R(x))} \geq \bigl(\tfrac{d}{D}\bigr)^l
\bigl(1+\bigl(\tfrac{d}{D}\bigr)^2\bigr)^{(n-1)l/2}
\end{equation}
so, in particular,
\begin{equation} \lb{5.4}
\liminf_n\, \abs{\det(p_n^R(x))}^{1/nl} \geq
\bigl(1+\bigl(\tfrac{d}{D}\bigr)^2 \bigr)^{1/2} \geq 1.
\end{equation}
But
\begin{equation} \lb{5.5}
\abs{\det(p_n^R(x))} = \abs{\det(A_1\cdots A_n)}^{-1} \exp (-nl\Phi_
{\nu_n}(x))
\end{equation}
where $\Phi_\nu$ is the potential of the measure $\nu$. Let $\nu_
\infty$ be a limit point of
$\nu_n$ and use equality in \eqref{5.1} and \eqref{5.4} to conclude,
for $x\notin\cvh(E)$,
\[
\exp(-\Phi_{\nu_\infty}(x))\geq C(E)
\]
which, as in the proof of Theorem~2.4 of \cite{EqMC}, implies that $
\nu_\infty=\rho_e$.
\end{proof}

The analogue of the almost converse of this last theorem has an extra
subtlety relative to
the scalar case:

\begin{theorem}\lb{T5.3} Let $d\mu$ be a non-trivial $l\times l$
matrix-valued measure on
$\bbR$ with $E=\supp(d\mu)$ compact. If $d\nu_n\to d\rho_E$, then
either $\mu$ is regular
or else, with $d\mu=M(x)\, d\mu_\tr(x)$, there is a set $S$ of
capacity zero, so $\det(M(x))=0$
for $d\mu_\tr$-a.e.\ $x\notin S$.
\end{theorem}

\begin{remark} By taking direct sums of regular point measures, it is
easy to find regular
measures where $\det(M(x))=0$ for $d\mu_\tr$-a.e.\ $x$.
\end{remark}

\begin{proof} For a.e.\ $x$ with $\det(M(x))\neq 0$, we have (see
Lemma~\ref{L5.5} below)
\[
p_n^R(x)\leq C(n+1)\bdone.
\]
The theorem then follows from the proof of Theorem~2.5 of \cite{EqMC}.
\end{proof}

\subsection{General Asymptotics} \lb{s5.3}

The following generalizes Theorem~1.10 of \cite{EqMC} from OPRL to
MOPRL---it is the matrix analogue
of a basic result of Stahl--Totik \cite{StT}. Its proof is
essentially the same as in \cite{EqMC}. By $\sigma_\ess(\mu)$, we
mean the
essential spectrum of the block Jacobi matrix associated to $\mu$.

\begin{theorem}\lb{T5.4} Let $E\subset\bbR$ be compact and let $\mu$
be an $l\times l$ matrix-valued
measure of compact support with $\sigma_\ess(\mu)=E$.
Then the following are equivalent:
\begin{SL}
\item[{\rm{(i)}}] $\mu$ is regular, that is, $\lim_{n\to\infty} \abs
{\det(A_1 \cdots A_n)}^{1/n}
=C(E)^l$.

\item[{\rm{(ii)}}] For all $z$ in $\bbC$, uniformly on compacts,
\begin{equation}\lb{5.6}
\limsup \abs{\det(p_n^R(z))}^{1/n} \leq e^{G_E(z)}.
\end{equation}

\item[{\rm{(iii)}}] For q.e.\ $z$ in $E$, we have
\begin{equation}\lb{5.7}
\limsup \abs{\det(p_n^R(z))}^{1/n}\leq 1.
\end{equation}
\end{SL}

\smallskip
Moreover, if {\rm{(i)--(iii)}} hold, then
\begin{SL}
\item[{\rm{(iv)}}] For every $z\in\bbC\setminus\cvh(\supp(d\mu))$, we
have
\begin{equation}\lb{5.8}
\lim_{n\to\infty}\, \abs{\det(p_n^R(z))}^{1/n} = e^{G_E(z)}.
\end{equation}

\item[{\rm{(v)}}] For q.e.\ $z\in\partial\Omega$,
\begin{equation}\lb{5.9}
\limsup_{n\to\infty} \, \abs{\det(p_n^R(z))}^{1/n} =1.
\end{equation}
\end{SL}
\end{theorem}

\begin{remarks} 1. $G_E$, the potential theorists' Green's function
for $E$, is defined by $G_E(z)=
-\log(C(E))-\Phi_{\rho_E}(z)$.

\smallskip
2. There is missing here one condition from Theorem~1.10 of \cite
{EqMC} involving general polynomials.
Since the determinant of a sum can be much larger than the sum of the
determinants, it is not obvious
how to extend this result.
\end{remarks}

\subsection{Weak Convergence of the CD Kernel and Consequences} \lb
{s5.3A}

The results of Simon in \cite{Weak-cd} extend to the matrix case. The
basic result is:

\begin{theorem}\lb{T5.4A} The measures $d\nu_n$ and $\f{1}{(n+1)l}\Tr
(K_n(x,x))\, d\mu(x)$
have the same weak limits. In particular, if $d\mu$ is regular,
\begin{equation}\lb{5.9a}
\f{1}{(n+1)l}\, \Tr (K_n(x,x)d\mu(x)) \overset{w}{\longrightarrow} d
\rho_E.
\end{equation}
\end{theorem}

As in \cite{Weak-cd}, $(\pi_n M_x \pi_n)^j$ and $(\pi_n M_x^j \pi_n)$
have a difference of
traces which is bounded as $n\to\infty$, and this implies the result.
Once one has this,
combining it with Theorem~\ref{T2.18B} leads to:

\begin{theorem}\lb{T5.4B} Let $I=[a,b]\subset E\subset\bbR$ with $E$
compact. Let $\sigma_\ess
(d\mu)=E$ for an $l\times l$ matrix-valued measure, and suppose $d\mu$
is regular for $E$ and
\begin{equation}\lb{5.9b}
d\mu = W(x)\, dx + d\mu_\s
\end{equation}
where $d\mu_\s$ is singular and $\det(W(x))>0$ for a.e.\ $x\in I$. Then,
\begin{alignat}{2}
& \text{\rm{(1)}} \qquad && \lim_{n\to\infty} \int_I p_n^R(x)^\dagger
\, d\mu_\s (x)
p_n^R(x)\to \boldsymbol{0}, \lb{5.9c} \\
& \text{\rm{(2)}} \qquad && \int_I\, \biggl\| \f{1}{n+1} \sum_{j=0}^n
p_j^R(x)^\dagger
W(x) p_j(x) -\rho_E(x)\bdone \biggr\|\, dx \to 0. \lb{5.9d}
\end{alignat}
\end{theorem}

\subsection{Widom's Theorem} \lb{s5.4}

\begin{lemma}\lb{L5.5} Let $d\mu$ be an $l\times l$ matrix-valued
measure supported on a compact
$E\subset\bbR$ and let $d\eta$ be a scalar measure on $\bbR$ so
\begin{equation} \lb{5.10}
d\mu(x) =W(x)\, d\eta(x) + d\mu_\s (x)
\end{equation}
where $d\mu_\s$ is $d\eta$ singular. Suppose for $d\eta$-a.e.\ $x$,
\begin{equation} \lb{5.11}
\det(W(x)) >0.
\end{equation}
Then for $d\eta$-a.e.\ $x$, there is a positive real function $C(x)$
so that
\begin{equation} \lb{5.12}
\norm{p_n^R(x)} \leq C(x) (n+1)\bdone.
\end{equation}
In particular,
\begin{equation} \lb{5.13}
\abs{\det(p_n^R(x))} \leq C(x)^l (n+1)^l.
\end{equation}
\end{lemma}

\begin{proof} Since $\norm{p_n^R}_R^2=1$, we have
\begin{equation} \lb{5.14}
\sum_{n=0}^\infty (n+1)^{-2} \norm{p_n^R}_R^2 <\infty
\end{equation}
so
\[
\sum_{n=0}^\infty (n+1)^{-2} \Tr (p_n^R(x)^\dagger W(x) p_n^R(x)) <
\infty
\]
for $d\eta$-a.e.\ $x$. Since \eqref{5.11} holds, for a.e.\ $x$,
\[
W(x)\geq b(x)\bdone
\]
for some scalar function $b(x)$. Thus, for a.e.\ $x$,
\[
\sum_{n=0}^\infty (n+1)^{-2} \Tr(p_n^R(x)^\dagger p_n^R(x)) \leq C(x)^2.
\]

Since $\norm{A}^2\leq \Tr(A^\dagger A)$, we find \eqref{5.12}, which
in turn implies \eqref{5.13}.
\end{proof}

This lemma replaces Lemma~4.1 of \cite{EqMC} and then the proof there
of Theorem~1.12 extends to give
(a matrix version of the theorem of Widom \cite{Wid}):

\begin{theorem}\lb{T5.6} Let $d\mu$ be an $l\times l$ matrix-valued
measure with $\sigma_\ess (d\mu)
=E\subset\bbR$ compact. Suppose
\begin{equation} \lb{5.15}
d\mu(x) =W(x)\, d\rho_E(x) + d\mu_\s (x)
\end{equation}
with $d\mu_\s$ singular with respect to $d\rho_E$. Suppose for $d
\rho_E$-a.e.\ $x$, $\det(W(x))>0$.
Then $\mu$ is regular.
\end{theorem}

\subsection{A Conjecture} \lb{s5.5}

We end our discussion of regular MOPRL with a conjecture---an analog
of a theorem of Stahl--Totik
\cite{StT}; see also Theorem~1.13 of \cite{EqMC} for a proof and
references. We expect the key will be
some kind of matrix Remez inequality. For direct sums, this
conjecture follows from Theorem~1.13 of
\cite{EqMC}.

\begin{conjecture}\lb{Con5.7} Let $E$ be a finite union of disjoint
closed intervals in $\bbR$. Suppose
$\mu$ is an $l\times l$ matrix-valued measure on $\bbR$ with $\sigma_
\ess (d\mu) =E$. For each $\eta>0$
and $m=1,2,\dots$, define
\begin{equation} \lb{5.16}
S_{m,\eta} = \{x \,\colon \mu([x-\tfrac{1}{m}, x+\tfrac{1}{m}])\geq e^
{-\eta m}\bdone\}.
\end{equation}
Suppose that for each $\eta$ (with $\abs{\cdot}=$ Lebesgue measure)
\[
\lim_{m\to\infty} \, \abs{E\setminus S_{m,\eta}} =0.
\]
Then $\mu$ is regular.
\end{conjecture}

\medskip

\noindent{\bf Acknowledgments.} It is a pleasure to thank Alexander Aptekarev, Christian Berg,
Antonio Dur\'an, Jeff Geronimo, Fritz~Gesztesy, Alberto Gr\"unbaum, Paco Marcell\'an,
Ken McLaughlin, Hermann Schulz-Baldes, and Walter Van Assche for useful correspondence.
D.D.\ was supported in part by NSF grants DMS--0500910 and
DMS-0653720. B.S.\ was supported in part by NSF grants DMS--0140592
and DMS-0652919 and U.S.--Israel Binational Science Foundation
(BSF) Grant No.\ 2002068.

\bigskip

\bigskip
\noindent David Damanik \\
\noindent Department of Mathematics \\
\noindent Rice University \\
\noindent Houston, TX 77005, USA \\
\noindent damanik@rice.edu

\bigskip
\noindent Alexander Pushnitski \\
\noindent Department of Mathematics \\
\noindent King's College London \\
\noindent Strand, London WC2R 2LS, UK \\
\noindent alexander.pushnitski@kcl.ac.uk

\bigskip
\noindent Barry Simon \\
\noindent Mathematics 253--37 \\
\noindent California Institute of Technology \\
\noindent Pasadena, CA 91125, USA \\
\noindent bsimon@caltech.edu

\endddoc
